\documentclass[11pt]{amsart}

\usepackage[margin=1in]{geometry}
\usepackage{amsmath,amssymb,amsfonts,amsthm}
\usepackage{mathrsfs}
\usepackage{enumitem}
\usepackage{graphicx}

\usepackage{color}

\newcommand{\hol}{\ensuremath{\mathcal{O}}}

\newcommand\Ga{\Gamma}

\newcommand\ga{\gamma}

\newcommand{\FF}{\ensuremath{\mathbb{F}}}

\newcommand{\PP}{\mathbb{P}}
\newcommand{\QQ}{\mathbb{Q}}

\newcommand{\ZZ}{\mathbb{Z}}

\newcommand{\rank}{\mathrm{rank\,}}

\usepackage{array,booktabs,multirow}
\usepackage{tikz}
\usetikzlibrary{arrows.meta,calc,intersections,decorations.pathmorphing}

\usepackage{tikz-cd}
\usepackage[colorlinks=true,linkcolor=blue,citecolor=blue,urlcolor=blue]{hyperref}

\theoremstyle{plain}
\newtheorem{theorem}{Theorem}[section]
\newtheorem{proposition}[theorem]{Proposition}
\newtheorem{lemma}[theorem]{Lemma}

\theoremstyle{definition}
\newtheorem{definition}[theorem]{Definition}
\newtheorem{example}[theorem]{Example}

\theoremstyle{remark}
\newtheorem{remark}[theorem]{Remark}

\title[Odd and even fake quadrics isogenous to a product]{Odd and even fake quadrics isogenous to a product}

\author{Fabrizio Catanese}
\address{Accademia Nazionale dei Lincei, Palazzo Corsini, Via della Lungara 10, 00165 Roma, Italia; Mathematisches Institut der Universit\"at Bayreuth, NW II, Universit\"atsstr. 30, 95447 Bayreuth, Germany}
\email{Fabrizio.Catanese@gmail.com, Fabrizio.Catanese@lincei.it, fabrizio.catanese@uni-bayreuth.de}

\author{JongHae Keum}
\address{School of Mathematics, Korea Institute for Advanced Study, Seoul 02455, Republic of Korea}
\email{jhkeum@kias.re.kr}

\author{Do Geun Kim}
\address{Department of Mathematics, POSTECH, 77 Cheongam-ro, Nam-gu, Pohang, Gyeongbuk, 37673, Korea}
\email{rlaehrms@postech.ac.kr}

\author{Kyoung-Seog Lee}
\address{Department of Mathematics, POSTECH, 77 Cheongam-ro, Nam-gu, Pohang, Gyeongbuk, 37673, Korea}
\email{kyoungseog@postech.ac.kr}

\date{\today}

\begin{document}
\raggedbottom

\begin{abstract}
A $\mathbb{Q}$-homology quadric is a smooth compact complex surface which has the same $\mathbb{Q}$-homology as a quadric surface, that is, the product of two $2$-dimensional spheres. It is called fake if it is minimal and not rational. Among them are the fake quadrics which are isogenous to a product of complex curves. In this paper, we determine the parity of the intersection forms of all fake quadrics isogenous to a product of unmixed type. In particular, we show that there are  fake quadrics with odd intersection forms.
\end{abstract}
\maketitle

\tableofcontents

\section{Introduction}

In the paper \cite{CasciniCataneseChenKeum25}, the authors studied $\mathbb{Q}$-homology quadrics.

\begin{definition}
A $\QQ$-homology quadric is a smooth compact complex surface $S$ with Betti numbers $b_1(S)=0, b_2(S)=2$: it is algebraic and has $K^2_S=8$. The surface $S$ is called a {\bf  fake quadric} if moreover it is minimal and of general type (equivalently,  the  canonical divisor $K_S$ is nef).
\end{definition}

For example,  the Segre-Hirzebruch surfaces $$\FF_m = \mathrm{Proj} (\hol_{\PP^1} \oplus \hol_{\PP^1} (m))$$
  are $\mathbb{Q}$-homology quadrics which are not fake quadrics.

\begin{definition}
A $\mathbb{Q}$-homology quadric $S$ is such that its Severi lattice $\mathrm{Num}(S)$, the group of divisor classes modulo numerical equivalence, has rank 2 and coincides with $ H^2(S, \ZZ)/\mathrm{Tors}(S)$. A $\mathbb{Q}$-homology quadric $S$ is
called an {\bf even}, resp. {\bf odd}, homology quadric if the intersection form on $\mathrm{Num}(S)$, which is unimodular  by Poincar\'e duality,
is even, resp. odd.
\end{definition}

For example,  the parity of  $\FF_m$   is the same as the parity of $m$. \\

Let us also mention the following problems posed in \cite{CasciniCataneseChenKeum25}.

\medskip

\noindent\textbf{Problem 1.}\cite{CasciniCataneseChenKeum25} Does  there exist an odd fake quadric? \\
(a) Determine which fake quadrics are even. (b) Determine which fake quadrics are odd.

\medskip

\noindent\textbf{Problem 2.}\cite{CasciniCataneseChenKeum25} Can an odd fake quadric contain a  curve of negative self-intersection?

\medskip

\noindent\textbf{Problem 3.}\cite{CasciniCataneseChenKeum25} Can an odd fake quadric have two fibrations?

\medskip

A surface $S$ is said to be isogenous to a higher product of unmixed type,
if
$$S\cong (C_1 \times C_2)/G,$$
where $C_i$ is a smooth compact curve with $g(C_i)\ge 2$ and $G$ is a finite group acting faithfully on $C_1$ and on $C_2$,
 such that  the diagonal action $ g(x,y) : = (gx, gy)$  of $G$ on $C_1\times C_2$ is free.  If moreover $S$ has geometric genus $p_g(S)=0$ and irregularity $q(S)=0$, then $S$ is a fake quadric.

The fake quadrics isogenous to a product of unmixed type were classified in \cite{BCG08} into 12 families and the parity of 5 families among the 12 was determined to be even in \cite{CasciniCataneseChenKeum25}. 
Our main result is the determination of the parity for the remaining 7 families, giving a positive answer to Problem 1 for such surfaces. It turns out that three families among the 12 are odd. The existence of an odd fake quadric was announced in \cite{odd}, but\footnote{due to a mistake concerning linearized 2-torsion bundles} it was claimed that the first surface in the  list below was odd.

All fake quadrics isogenous to a product are Mori dream surfaces having two fibrations and no negative curves (cf. \cite{FrapportiLee22, KeumLee19, KeumLee26}). Therefore, our result also answers positively Problem 3.

\begin{theorem}\label{main} In the table below, the parity of all the 12 families are given, new results in capital letters. Also given are the numbers $\#\Theta(C_i)^G$ of $G$-invariant theta characteristics for the $G$ action on $C_i$ with signature $T_i$ (this means: $(0; T_i)$). In particular, three families among the 12 are odd fake quadrics.

\begin{table}[h!]
\centering
\resizebox{\textwidth}{!}{%
\begin{tabular}{c|c|c|c|c|c|c|c|c}
$G$ & $g(C_1)$ & $g(C_2)$ & $T_1$ & $T_2$ & $H_1(S,\mathbb{Z})$ & $\#\Theta(C_1)^G$ & $\#\Theta(C_2)^G$ & parity \\
\hline
$\mathfrak{A}_5$ & 4 & 21 & $[2,5,5]$ & $[3,3,3,3]$ & $(\mathbb{Z}/3)^2 \times \mathbb{Z}/{15}$ & 1 & 2 & EVEN \\

$\mathfrak{A}_5$ & 13 & 6 & $[5,5,5]$ & $[2,2,2,3]$ & $(\mathbb{Z}/{10})^2$ & 1 & 4 & ODD \\

$\mathfrak{A}_5$ & 5 & 16 & $[3,3,5]$ & $[2,2,2,2,2]$ & $(\mathbb{Z}/2)^3 \times \mathbb{Z}/6$ & 1 & 16 & ODD \\

$\mathfrak{S}_4 \times \mathbb{Z}/2$ & 3 & 25 & $[2,4,6]$ & $[2,2,2,2,2,2]$ & $(\mathbb{Z}/2)^4 \times \mathbb{Z}/4$ & 2 & 32 & EVEN \\

$G(32)$ & 9 & 5 & $[2,2,4,4]$ & $[2,2,2,4]$ & $(\mathbb{Z}/2)^2 \times \mathbb{Z}/4 \times \mathbb{Z}/8$ & 8 & 8 & ODD \\

$(\mathbb{Z}/5)^2$ & 6 & 6 & $[5,5,5]$ & $[5,5,5]$ & $(\mathbb{Z}/5)^3$ & 1 & 1 & even \\

$\mathfrak{S}_4$ & 3 & 13 & $[3,4,4]$ & $[2,2,2,2,2,2]$ & $(\mathbb{Z}/2)^4 \times \mathbb{Z}/8$ & 2 & 32 & even \\

$G(16)$ & 5 & 5 & $[2,2,4,4]$ & $[2,2,4,4]$ & $(\mathbb{Z}/2)^2 \times \mathbb{Z}/4 \times \mathbb{Z}/8$ & 8 & 8 & even \\

$D_4 \times \mathbb{Z}/2$ & 3 & 9 & $[2,2,2,4]$ & $[2,2,2,2,2,2]$ & $(\mathbb{Z}/2)^3 \times (\mathbb{Z}/4)^2$ & 4 & 32 & EVEN \\

$(\mathbb{Z}/2)^4$ & 5 & 5 & $[2,2,2,2,2]$ & $[2,2,2,2,2]$ & $(\mathbb{Z}/4)^4$ & 32 & 32 & even \\

$(\mathbb{Z}/3)^2$ & 4 & 4 & $[3,3,3,3]$ & $[3,3,3,3]$ & $(\mathbb{Z}/3)^5$ & 1 & 1 & even \\

$(\mathbb{Z}/2)^3$ & 3 & 5 & $[2,2,2,2,2]$ & $[2,2,2,2,2,2]$ & $(\mathbb{Z}/2)^4 \times (\mathbb{Z}/4)^2$ & 8 & 32 & EVEN \\

\end{tabular}%
}
\end{table}

\end{theorem}

\begin{remark} Our method also applies to recover the previous determination for the 5 families in \cite{CasciniCataneseChenKeum25}. The parity of the case with $G=D_4 \times \mathbb{Z}/2$ has been determined by a different method in \cite{CataneseFrapporti26}.
\end{remark}

Let us explain how to determine the parity of the intersection form  of $S = (C_1 \times C_2) / G$ for each family.  First we recall the following exact sequence.

\begin{proposition}[{\cite[Proposition~2.2]{Dolgachev99}}, {\cite[Remark~7.2]{Dolgachev03}}]\label{prop:dolgachev_exact}
Let $G$ be a finite group and let $X$ be a smooth projective variety with $G$-action. Then we have the following exact sequence of abelian groups
\begin{equation}\label{eq:exact}
0 \longrightarrow \operatorname{Hom}(G,\mathbb{C}^*)
\longrightarrow \operatorname{Pic}(G;X)
\longrightarrow \operatorname{Pic}(X)^{G}
\longrightarrow H^{2}(G,\mathbb{C}^*).
\end{equation}

Here, $\operatorname{Pic}(G;X)$ denotes the group of $G$-linearized line bundles, and $\operatorname{Pic}(X)^{G}$ denotes the group of divisor classes which are fixed by $G$. Moreover, the map $\operatorname{Pic}(X)^G \longrightarrow H^{2}(G,\mathbb{C}^*)$ is surjective when $X$ is a smooth projective curve with genus $g\geq 2$.
\end{proposition}

We analyse the above exact sequences when $X = C_1 \times C_2$ or $X=C$.

Using the results of \cite{Dolgachev99} and group theory, we are able to compute all the terms in \eqref{eq:exact} except $\operatorname{Pic}(C)^G$.
Indeed, the groups $\operatorname{Hom}(G,\mathbb{C}^*)$ and $H^{2}(G,\mathbb{C}^*)$ depend only on  the group $G$,  and their values can be read off from the webpage:
\cite[\href{https://www.lmfdb.org/Groups/Abstract}{LMFDB: Abstract group}]{LMFDB}.
Moreover, we can compute $\operatorname{Pic}(G;C)$ using the following results.

\begin{proposition}[{\cite[Proposition~5.5]{BradenDisneyHogg25}, \cite{Dolgachev99}}]\label{Pic(G;C)}
Suppose that the action of $G$ on $C$ has $g(C/G)=0$ and  signature $T=[m_1,\dots,m_r]$, with $m_1\leq\cdots\leq m_r$, where $m_i$   is the order of the cyclic stabilizer of a ramification point  over the $i$-th branch point $p_i$ of the cover $\pi: C\to C/G$. Then, $m_i p'_i$ being  the
pull-back divisor of the point $p_i$, we have
\[
\operatorname{Pic}(G;C)\ \cong\  (\bigoplus_1^r \mathbb{Z} p'_i )/ \langle m_1 p'_1 - m_2 p'_2 , m_2 p'_2-m_3 p'_3, \dots , m_{r-1} p'_{r-1}  - m_r p'_r \rangle ,
\]
which in turn is isomorphic to
\[
\operatorname{Pic}(G;C)\ \cong\ \mathbb{Z}\ \oplus\ (\bigoplus_{i=1}^{r-1}\ \mathbb{Z}/(d_i/d_{i-1})\mathbb{Z}\ ),
\]
where
\[
d_i:= \gcd\Bigl\{\ \prod_{j\in S} m_j\ \Bigm|\ S\subset\{1,\dots,r\},\ |S|=i\ \Bigr\}
~~~ \text{and we set } d_0=1.
\]

Moreover, the free part of $\operatorname{Pic}(G;C)$ is generated by a line bundle of degree
\[
\frac{\deg(\pi)}{\operatorname{lcm}(m_1, \dots, m_r)}=\frac{|G|}{\operatorname{lcm}(m_1, \dots, m_r)}.
\]
\end{proposition}

It remains to determine $\operatorname{Pic}(C)^G$. First, we compute  its 2-torsion subgroup $\operatorname{Pic}(C)^G[2]$ and  the number of $G$-invariant theta characteristics on $C$.
Let
\[
\Theta(C)
=
\{
\kappa\in \operatorname{Pic}^{g-1}(C)
\mid
\kappa^{\otimes 2}\simeq K_C
\}
\]
be the set of theta characteristics of $C$, and let
\[
\Theta(C)^G
=
\{
\kappa\in\Theta(C)
\mid
g^*\kappa\simeq\kappa
\text{ for every }g\in G
\}.
\]
be the set of $G$-invariant theta characteristics of $C$.

\begin{lemma}\label{2torsions} Suppose that there is a $G$-invariant theta characteristic on $C$. Then the number of $G$-invariant theta characteristics on $C$ is equal to the order of $\operatorname{Pic}(C)^G[2]$, that is,
$$ \#\Theta(C)^G= \#\operatorname{Pic}(C)^G[2].$$
\end{lemma}

In the next sections we will compute $\#\Theta(C_1)^G$ and $\#\Theta(C_2)^G$ for each of the 12 families.
Using this information together with  the exact sequence \eqref{eq:exact} for $C_1$ and $ C_2$, we determine $\operatorname{Pic}(C_i)^G$ (observing that here the group of Schur multipliers in the new cases is of exponent  2. Then, from these exact sequences and the ideas from \cite{KimKimLee17,LeeThesis14,Lee15,Lee16,LeeShabalin18}, we can explicitly construct a basis of $H^2(S,\mathbb{Z})_{\mathrm{free}}$ and compute the intersection form. More precisely, we use the following strategy.

\bigskip

\noindent\textbf{STRATEGY}

\medskip

Since the $G$-action on $C_1\times C_2$ is free, giving a line bundle on
$S=(C_1\times C_2)/G$ is equivalent to giving a $G$-linearized line bundle on $C_1\times C_2$.  The class $L'$ on $S=(C_1\times C_2)/G$ induced by a $G$-linearized line bundle $L$ on $C_1\times C_2$ is called the {\bf descent } of $L$.

\begin{enumerate}
\item If we take  two $G$-linearized line bundles on $C_1\times C_2$ of bidegree $(a, \, b)$ and $(c, \, d)$, respectively,  then their descent line bundles have  intersection number
\[
\frac{ (a, \, b)\cdot(c, \, d)}{|G|}=\frac{ad+bc}{|G|}=\frac{ad+bc}{(g(C_1)-1)(g(C_2)-1)}.
\]

\item Let $L_i\in \operatorname{Pic}(C_i)^G$, $i=1,2$. If both $L_1$ and $L_2$ are $G$-linearized, then the external tensor product
$L_1\boxtimes L_2$ is a $G$-linearized line bundle on $C_1\times C_2$. If neither is $G$-linearized, but their obstruction classes in the  Schur multiplier coincide with each other, then $L_1\boxtimes L_2^{-1}$ has obstruction $0$, hence is a $G$-linearized line bundle on $C_1\times C_2$.
This is explained by the following ladder diagram:
\[
\begin{tikzcd}[row sep=2.8em, column sep=2em]
0 \arrow[r] & \operatorname{Hom}(G,\mathbb C^*) \arrow[r] \arrow[d, equal] & \operatorname{Pic}(G;C_{1}) \arrow[r] \arrow[d] & \operatorname{Pic}(C_{1})^{G} \arrow[r, "\delta_1"] \arrow[d, hook] & H^{2}(G;\mathbb{C}^*) \arrow[r] \arrow[d, equal] & 0 \\
0 \arrow[r] & \operatorname{Hom}(G,\mathbb C^*) \arrow[r] & \operatorname{Pic}(G;C_{1}\times C_{2}) \arrow[r] & \operatorname{Pic}(C_{1}\times C_{2})^{G} \arrow[r] & H^{2}(G;\mathbb{C}^*) \arrow[r] & 0
\\
0 \arrow[r] & \operatorname{Hom}(G,\mathbb C^*) \arrow[r] \arrow[u, equal] & \operatorname{Pic}(G;C_{2}) \arrow[r] \arrow[u] & \operatorname{Pic}(C_{2})^{G} \arrow[r,"\delta_2"] \arrow[u, hook] & H^{2}(G;\mathbb{C}^*) \arrow[r] \arrow[u, equal] & 0
\end{tikzcd}
\]
 Moreover, if $L\in \operatorname{Pic}(C_1)^G$, then, since the map $\delta_2$ is surjective, there exists $M \in \operatorname{Pic}(C_2)^G$ such that $\delta_2(M)=\delta_1(L)$, hence the external tensor product
$L\boxtimes M^{-1}$ has obstruction 0, giving a $G$-linearized line bundle on $C_1\times C_2$.

\item If we can construct two $G$-linearized line bundles on $C_1\times C_2$ of bidegree $(g(C_1)-1, \,0)$ and $(0, \, g(C_2)-1)$, respectively, then the intersection numbers of their descents can be computed by (1) and the intersection matrix of $H^2(S,\mathbb{Z})_{\mathrm{free}}$ is
$
\begin{pmatrix}
0 & 1 \\
1 & 0
\end{pmatrix},
$
which is even.

\item If we can construct two $G$-linearized line bundles on $C_1\times C_2$ of bidegree $(g(C_1)-1, \, 0)$ and $(\frac{g(C_1)-1}{2}, \, g(C_2)-1)$, respectively, then again by (1) their descents have intersection matrix
$
\begin{pmatrix}
0 & 1 \\
1 & 1
\end{pmatrix},
$
and the intersection form of $H^2(S,\mathbb{Z})_{\mathrm{free}}$ is odd. To prove just oddness it is enough to construct a $G$-linearized line bundle on $C_1\times C_2$ whose descent line bundle has self-intersection odd. Similarly, a $G$-linearized line bundle on $C_1\times C_2$ of bidegree $(g(C_1)-1, \, \frac{g(C_2)-1}{2})$, if exists, also determines a  descent line bundle of  self-intersection 1.
\end{enumerate}

\bigskip
\noindent{\bf NOTATION}

\medskip

We will work over $\mathbb{C}.$ For a smooth projective curve $C$ and  a surface $S$ on which a finite group $G$ acts, we will use the following notations.

\smallskip

 \noindent $T=[g(C/G); m_1,\ldots,m_r]$: the signature of the quotient map $\pi\colon C\to C/G$, where $m_j$ are the ramification indices with $m_1\leq\cdots\leq m_r$. When $C/G\simeq\mathbb{P}^1$, we often write $T=[m_1,\ldots,m_r]$.

\smallskip

\noindent $\Gamma:=\pi_1^{\mathrm{orb}}\bigl(C/G\setminus\{p_1,\cdots,p_r\} ,(m_1,\ldots,m_r)\bigr)$: the orbifold fundamental group associated with $T=[g(C/G); m_1,\ldots,m_r]$, where $p_1,\ldots,p_r$ are the branch points. $\Gamma$ fits in the exact sequence $$1\longrightarrow \pi_1(C) \longrightarrow \Gamma \xrightarrow{\ \phi\ }G \longrightarrow 1$$ and admits the presentation
(cf. \cite[Definition~5.5, Proposition~5.8, Example~5.10]{Cat08})
\[
\Gamma=\left\langle\alpha_1,\beta_1,\ldots,\alpha_g,\beta_g,\gamma_1,\ldots,\gamma_r \mid\prod_{j=1}^{g}[\alpha_j,\beta_j]\,\gamma_1\cdots\gamma_r=1,\,\,\gamma_i^{m_i}=1,\,\, 1\leq i\leq r \right\rangle
\]
where $g=g(C/G).$ If $g=0$, this reduces to $\Gamma=\left\langle\gamma_1,\ldots,\gamma_r \mid \gamma_1\cdots\gamma_r=1,\;\gamma_i^{m_i}=1,\,\,1\leq i\leq r \right\rangle.$

\smallskip

\noindent $\Theta(C)$: the set of theta characteristics of $C$.

\noindent $\Theta(C)^G$: the set of $G$-invariant theta characteristics of $C$.

\noindent $\operatorname{Num}(S)$: the group of divisors on $S$ modulo numerical equivalence.

\noindent $\operatorname{Pic}(C)$: the group of isomorphism classes of line bundles on $C$.

\noindent $\operatorname{Pic}(C)^G$: the group of isomorphism classes of $G$-invariant line bundles on $C$.

\noindent $\operatorname{Pic}(G;C)$: the group of $G$-linearized line bundles on $C$.

\noindent $\delta\colon\operatorname{Pic}(C)^G\to H^2(G;\mathbb{C}^*)$: the obstruction map for  $G$-linearization.

\noindent $L_1\boxtimes L_2$: the external tensor product on $C_1\times C_2$  of line bundles $L_i$ on $C_i$.

\noindent $D_1\equiv D_2$: the divisors $D_1$ and $D_2$ are linearly equivalent.

\noindent $G(n,k)$: the $k$-th group of order $n$ in GAP Small Groups Library. $G(32):=G(32,27)$. 

\noindent $H^k(G;\mathbb{Z}/2)$ : the group cohomology with coefficients in $\mathbb{Z}/2$ (which is  a trivial $G$-module).

\noindent $\operatorname{res}_H$: the restriction homomorphism in group cohomology for a subgroup $H$ in $G$.

\noindent $\beta_G$: the Bockstein homomorphism for  a group $G$.

\bigskip

\section{Generalities on $G$-invariant divisors and group cohomology}

 In this section, we collect some general results about $G$-linearized sheaves and invariant divisors which will be used to analyse  $G$-actions on  curves and their products.

\begin{proposition}\label{admit theta}
Let $G$ be a finite group and $C$ be a smooth projective curve with a faithful $G$-action.  Suppose that the quotient map $\pi: C\to C/G$ has signature $T=[0;m_1,\ldots,m_r]$. Let $p_1,\ldots,p_r$ be the branch points and define  $E_i:=\pi^*(p_i)_{\mathrm{red}}$,  $m_i E_i:=\pi^*(p_i)$ and $F:=\pi^*(p)$ for a general point $p\in C/G\cong \mathbb P^1$. Assume that every $m_i$ is odd. Then
$$
K_C\equiv 2(-F+\sum_{i=1}^r \frac{m_i-1}{2} E_i)
$$
and $\kappa:=\mathcal{O}_C(-F+\sum_{i=1}^r \frac{m_i-1}{2} E_i)$ is a $G$-invariant theta
characteristic.
\end{proposition}
\begin{proof}
It follows from the Hurwitz formula.
\end{proof}

\begin{remark}
By the above proposition, the cases in the table with signatures $[3,3,3,3]$, $[5,5,5]$, and $[3,3,5]$ admit a $G$-invariant theta characteristic.
\end{remark}

\begin{lemma}
If there exists a $G$-invariant theta characteristic $\kappa_0\in\Theta(C)^G$, then $$\Theta(C)^G=\kappa_0\otimes H^1(C,\mathbb{Z}/2)^G$$
and $\#\Theta(C)^G=2^d$, $d=\dim_{\mathbb{Z}/2}H^1(C,\mathbb{Z}/2)^G$.
\end{lemma}

Assume that $C\longrightarrow C/G\cong\mathbb{P}^1$ has signature $[0;m_1,\ldots,m_r]$. Let
$\Gamma$ be the polyhedral group
\[
\Gamma : = P (m_1, \dots,  m_r) : = \langle x_1,\ldots,x_r \mid x_i^{m_i}=1,\;x_1\cdots x_r=1,\; i=1,\cdots,r \rangle.
\]
A surjective homomorphism $\varphi:\Gamma\twoheadrightarrow G, \quad x_i\longmapsto g_i,$
is determined by a spherical generating  vector $(g_1,\ldots,g_r) \in G^r $ such that  $\operatorname{ord}(g_i)=m_i$, $g_1\cdots g_r=1$, $\langle g_1,\ldots, g_r \rangle=G$. The  kernel of $\varphi$ is the ordinary fundamental group of $C$ and we have the exact sequence
$$1\longrightarrow\pi_1(C)\longrightarrow\Gamma\overset{\varphi}{\longrightarrow}G\longrightarrow1.$$

\medskip

We recall a consequence of the {\bf Lyndon--Hochschild--Serre spectral sequences Theorem}

\begin{theorem}
 Let $\Ga$  be any group, and
$K$ a normal subgroup, with  factor group $G:  = \Ga/K$.
Assume that $M$  is a $\Ga$-module, such that the respective modules of covariants and invariants,
$$M_{\Ga} : = M / \langle \{ \ga(a)-a| \ga \in \Ga, a \in M\} \rangle \ \ {\rm  and } \ M^{\Ga} : = \{ a | \ga(a) = a, \forall \ga \in \Ga\}$$
 are $G$-modules. Then there
are converging first quadrant spectral sequences
$$E^2_{p,q}
= H_p(G; H_q (K, M)) \Rightarrow H_{p+q} (\Ga, M)$$
$$E ^{p,q}_2
= H^p(G; H^q (K, M)) \Rightarrow H^{p+q} (\Ga, M).$$
 \end{theorem}

namely,

\begin{theorem}\cite[Theorem~2]{HochschildSerre53}
Let $\Gamma$ be a group, let $K$ be a normal subgroup of $\Gamma$, and let $M$ be a $\Gamma$-module. Let $m\geq 1$, and assume that $H^n(K,M)=0$, for $0<n<m$.
 Then we have the following exact sequence:
\[
0
\longrightarrow H^m(\Gamma/K,M^K)
\longrightarrow H^m(\Gamma,M)
\longrightarrow H^m(K,M)^\Gamma
\longrightarrow H^{m+1}(\Gamma/K,M^K)
\longrightarrow H^{m+1}(\Gamma,M).
\]
\end{theorem}

\begin{remark}\cite[Remark, p.130]{HochschildSerre53}
In the case $m=1$, the hypothesis of the preceding theorem is vacuous, and hence we have always the following exact sequence:
\[
0
\longrightarrow H^1(\Gamma/K,M^K)
\longrightarrow H^1(\Gamma,M)
\longrightarrow H^1(K,M)^\Gamma
\longrightarrow H^2(\Gamma/K,M^K)
\longrightarrow H^2(\Gamma,M).
\]
\end{remark}

We apply this to $K=\pi_1(C)$ and $M=\mathbb{Z}/2$, where $M$ is regarded as a trivial $\Gamma$-module. Since the action of $K=\pi_1(C)$ on $\mathbb{Z}/2$ is trivial,
\[
M^K=\{m\in\mathbb{Z}/2 \mid k\cdot m=m \text{ for every }k\in K\}=\mathbb{Z}/2.
\]
Hence the five-term exact sequence becomes
\[
0
\longrightarrow H^1(G,\mathbb{Z}/2)
\longrightarrow H^1(\Gamma,\mathbb{Z}/2)
\longrightarrow H^1(\pi_1(C),\mathbb{Z}/2)^\Gamma
\longrightarrow H^2(G,\mathbb{Z}/2)
\overset{\varphi^*}{\longrightarrow} H^2(\Gamma,\mathbb{Z}/2).
\]

The following results are well-known.

\begin{lemma}  \cite[Exercise 6.1.5]{Weibel94}\label{Hom} For each group $G$ and each abelian group $A$ with trivial $G$-action $A$, $H^1(G,A)=\operatorname{Hom}(G,A)\cong\operatorname{Hom}(G^{\mathrm{ab}},A)$. If the above applies to $ A = \mathbb{Z}/2$, and to the groups $G=\pi_1(C)$, we obtain $H^1(\pi_1(C),\mathbb{Z}/2)\cong H^1(C,\mathbb{Z}/2)$.
\end{lemma}

\begin{lemma}
$H^1(\pi_1(C),\mathbb{Z}/2)^\Gamma=H^1( C,\mathbb{Z}/2)^G$.
\end{lemma}
\begin{proof}
Since $K : =\pi_1(C)\triangleleft\Gamma$, every $\gamma\in\Gamma$ induces an automorphism of $K$ by conjugation: $k\longmapsto\gamma k\gamma^{-1}$. This induces an action of $\Gamma$ on
\[
H^1(K,\mathbb{Z}/2)=\operatorname{Hom}(K,\mathbb{Z}/2).
\]
Explicitly, for $\psi\in\operatorname{Hom}(K,\mathbb{Z}/2)$ and $\gamma\in\Gamma$, the action is given by $(\gamma\cdot\psi)(k)=\psi(\gamma^{-1}k\gamma)$.

Since these homomorphisms factor through the abelianization of $K$, clearly inner automorphisms act trivially, whence
 the action of $\Gamma$ on $H^1(K,\mathbb{Z}/2)$ is trivial on
$K$, and therefore factors through $\Gamma/K\cong G$. Consequently,
\[
H^1(K,\mathbb{Z}/2)^\Gamma=H^1(K,\mathbb{Z}/2)^G \cong H^1(C,\mathbb{Z}/2)^G.
\]
\end{proof}

By above lemma the five-term sequence becomes
\[
0
\longrightarrow H^1(G,\mathbb{Z}/2)
\longrightarrow H^1(\Gamma,\mathbb{Z}/2)
\longrightarrow H^1(C,\mathbb{Z}/2)^G
\longrightarrow H^2(G,\mathbb{Z}/2)
\overset{\varphi^*}{\longrightarrow} H^2(\Gamma,\mathbb{Z}/2).
\]

And we have the following.

\begin{proposition}\label{dim}
$\dim H^1(C,\mathbb{Z}/2)^G=\dim H^1(\Gamma,\mathbb{Z}/2)-\dim H^1(G,\mathbb{Z}/2)+\dim\ker\varphi^*$.
\end{proposition}

We now explain how to compute these terms.

\begin{proposition}\label{dim H1 Gamma}
Let $[m_1,\cdots,m_r]$ be the signature for  $\pi: C \to C/G \cong \PP^1$. Let $e:=\#\{i \mid m_i \text{ : even }\}$. Then $\dim_{\mathbb{Z}/2}H^1(\Gamma;\mathbb{Z}/2)=e-1$, $e>0$.
\end{proposition}

\begin{proof}
Let $\chi :\Gamma \longrightarrow \mathbb{Z}/2$ be a homomorphism. If $m_i$ is odd, then $m_i\chi(x_i)=0$ implies $\chi(x_i)=0$. If $m_i$ is even, the product relation gives $\sum_{m_i \text{ : even} }\chi(x_i)=0.$ Therefore, we obtain $e-1$ free variables for $\chi(x_i)$. This implies $\dim_{\mathbb{Z}/2}H^1(\Gamma;\mathbb{Z}/2)=e-1.$
\end{proof}

\begin{remark}\label{Schur}
The abelianization and Schur multiplier for some finite groups are well known  and can be readily found in \cite[\href{https://www.lmfdb.org/Groups/Abstract}{LMFDB: Abstract group}]{LMFDB} as follows. \\
\resizebox{\textwidth}{!}{
\renewcommand{\arraystretch}{1}
\begin{tabular}{c|c|c|c|c|c|c|c|c|c|c|c|c}
\toprule
Group & $\mathfrak{A}_5$& $\mathfrak{A}_5$ & $\mathfrak{A}_5$ & $\mathfrak{S}_4\times \mathbb{Z}/2$& $G(32)$  &  $(\mathbb{Z}/5)^2$& $\mathfrak{S}_4$&   $G(16)$  &   $D_4\times \mathbb{Z}/2$     &   $(\mathbb{Z}/2)^4$   &   $(\mathbb{Z}/3)^2$     &  $(\mathbb{Z}/2)^3$ \\
\midrule
$H_1(G;\mathbb{Z})$\text{(Abelianization)} &   0    &   0   &  0    &  $(\mathbb{Z}/2)^2$    &   $(\mathbb{Z}/2)^3$   &   $(\mathbb{Z}/5)^2$   &   $\mathbb{Z}/2$   &   $\mathbb{Z}/4\oplus \mathbb{Z}/2$   &  $(\mathbb{Z}/2)^3$    &   $(\mathbb{Z}/2)^4$   &  $(\mathbb{Z}/3)^2$  & $(\mathbb{Z}/2)^3$  \\
\midrule
$H_2(G;\mathbb{Z})$\text{(Schur multiplier)} &  $\mathbb{Z}/2$     &   $\mathbb{Z}/2$   &   $\mathbb{Z}/2$   &   $(\mathbb{Z}/2)^2$   &   $(\mathbb{Z}/2)^4$   &  $\mathbb{Z}/5$    &   $\mathbb{Z}/2$   &   $(\mathbb{Z}/2)^2$   &   $(\mathbb{Z}/2)^3$   &   $(\mathbb{Z}/2)^6$   &  $\mathbb{Z}/3$   &  $(\mathbb{Z}/2)^3$ \\
\bottomrule
\end{tabular}}\\

In particular, the above table allows us to determine $\operatorname{Hom}(G,\mathbb{C}^*)$ and $H^2(G;\mathbb{C}^*)$, which will be used later to analyse the exact sequence \eqref{eq:exact}. Indeed, we have $\operatorname{Hom}(G,\mathbb{C}^*)\cong H^1(G;\mathbb{C}^*)\cong \operatorname{Hom}(H_1(G;\mathbb{Z}),\mathbb{C}^*)$. Moreover, since $H^2(G;\mathbb{C}^*)\cong \operatorname{Hom}(H_2(G;\mathbb{Z}),\mathbb{C}^*)$, $H^2(G;\mathbb{C}^*)\cong H_2(G;\mathbb{Z})$.
\end{remark}
Using above remark, we compute $Ext_{\mathbb Z}^1(\mathbb Z/n,\mathbb{Z}/2)$, $\operatorname{Hom}_{\mathbb Z}(\mathbb Z/n,A)$.

\begin{lemma}\cite[Calculation~3.3.2]{Weibel94}
For a cyclic group $\mathbb Z/n$ and every abelian group $A$, we have $\operatorname{Hom}_{\mathbb Z}(\mathbb Z/n,A) \cong \{a\in A\mid na=0\}$, $Ext^1_{\mathbb Z}(\mathbb Z/n,A)\cong A/nA$. Hence, for $A=\mathbb{Z}/2$,\\
\[
Ext_{\mathbb Z}^1(\mathbb Z/n,\mathbb{Z}/2)\cong(\mathbb{Z}/2)/n(\mathbb{Z}/2)=
\begin{cases}
\mathbb{Z}/2,& n \text{ even}\\
0,& n \text{ odd}
\end{cases}
\quad \text{and} \quad
\operatorname{Hom}_{\mathbb Z}(\mathbb Z/n,\mathbb{Z}/2)=
\begin{cases}
\mathbb{Z}/2,& n \text{ even}\\
0,& n \text{ odd}.
\end{cases}
\]
\end{lemma}

By the above lemma and the universal coefficient theorem, we obtain $H^2(G;\mathbb{Z}/2)$.

\begin{proposition}\label{dim H2 G}
$$
\resizebox{\textwidth}{!}{
\renewcommand{\arraystretch}{1}
\begin{tabular}{c|c|c|c|c|c|c|c|c|c|c|c|c}
\toprule
Group& $\mathfrak{A}_5$& $\mathfrak{A}_5$& $\mathfrak{A}_5$& $\mathfrak{S}_4\times \mathbb{Z}/2$& $G(32)$& $(\mathbb{Z}/5)^2$& $\mathfrak{S}_4$& $G(16)$& $D_4\times \mathbb{Z}/2$& $(\mathbb{Z}/2)^4$& $(\mathbb{Z}/3)^2$& $(\mathbb{Z}/2)^3$\\
\midrule
$Ext^1_{\mathbb{Z}}(H_1(G;\mathbb{Z}),\mathbb{Z}/2)$& $0$ & $0$& $0$& $(\mathbb{Z}/2)^2$& $(\mathbb{Z}/2)^3$& $0$& $\mathbb{Z}/2$& $(\mathbb{Z}/2)^2$& $(\mathbb{Z}/2)^3$& $(\mathbb{Z}/2)^4$& $0$& $(\mathbb{Z}/2)^3$\\
\midrule
$\operatorname{Hom}(H_2(G;\mathbb{Z}),\mathbb{Z}/2)$& $\mathbb{Z}/2$& $\mathbb{Z}/2$& $\mathbb{Z}/2$& $(\mathbb{Z}/2)^2$& $(\mathbb{Z}/2)^4$& $0$& $\mathbb{Z}/2$& $(\mathbb{Z}/2)^2$& $(\mathbb{Z}/2)^3$& $(\mathbb{Z}/2)^6$& $0$& $(\mathbb{Z}/2)^3$\\
\midrule
$H^2(G;\mathbb{Z}/2)$& $\mathbb{Z}/2$& $\mathbb{Z}/2$& $\mathbb{Z}/2$& $(\mathbb{Z}/2)^4$& $(\mathbb{Z}/2)^7$& $0$& $(\mathbb{Z}/2)^2$& $(\mathbb{Z}/2)^4$& $(\mathbb{Z}/2)^6$& $(\mathbb{Z}/2)^{10}$& $0$& $(\mathbb{Z}/2)^6$\\
\bottomrule
\end{tabular}
}
$$
\end{proposition}

By Proposition~\ref{dim H1 Gamma} and ~\ref{dim H2 G}, we have the following table.

\begin{center}\label{dim table}
\resizebox{0.66\textwidth}{!}{
\renewcommand{\arraystretch}{1}
\begin{tabular}{c|c|c|c|c|c|c}
\toprule
Family& Curve& Group & Signature & $\dim H^1(\Gamma,\mathbb{Z}/2)$&$\dim H^1(G,\mathbb{Z}/2)$ & $\dim H^2(G,\mathbb{Z}/2)$ \\
\midrule
1 & $C_1$ & $\mathfrak{A}_5$ &$[2,5,5]$& $0$&  $0$  & $1$    \\
1 & $C_2$ & $\mathfrak{A}_5$ &$[3,3,3,3]$& $0$&  $0$  & $1$   \\
\midrule
2 & $C_1$ & $\mathfrak{A}_5$ &$[5,5,5]$& $0$&  $0$  & $1$    \\
2 & $C_2$ & $\mathfrak{A}_5$ &$[2,2,2,3]$& $2$&  $0$  & $1$    \\
\midrule
3 & $C_1$ & $\mathfrak{A}_5$ &$[3,3,5]$& $0$&  $0$  & $1$    \\
3 & $C_2$ & $\mathfrak{A}_5$ &$[2,2,2,2,2]$ & $4$ & $0$ &  $1$    \\
\midrule
4 & $C_1$ & $\mathfrak{S}_4\times \mathbb{Z}/2$ &$[2,4,6]$& $2$ &  $2$  & $4$  \\
4 & $C_2$ & $\mathfrak{S}_4\times \mathbb{Z}/2$ &$[2,2,2,2,2,2]$ & $5$ &  $2$  &$4$   \\
\midrule
5 & $C_1$ & $G(32)$ &$[2,2,4,4]$& $3$ &  $3$  &$7$  \\
5 & $C_2$ & $G(32)$ &$[2,2,2,4]$ & $3$ &  $3$  &$7$  \\
\midrule
6 & $C_1$ & $(\mathbb{Z}/5)^2$ &$[5,5,5]$& $0$ &  $0$  &$0$  \\
6 & $C_2$ & $(\mathbb{Z}/5)^2$ &$[5,5,5]$ & $0$ &  $0$  &$0$  \\
\midrule
7 & $C_1$ & $\mathfrak{S}_4$ &$[3,4,4]$ &  $1$    &  $1$   & $2$  \\
7 & $C_2$ & $\mathfrak{S}_4$ &$[2,2,2,2,2,2]$ & $5$ &  $1$   & $2$  \\
\midrule
8 & $C_1$ & $G(16)$ &$[2,2,4,4]$& $3$ &  $2$   &$4$ \\
8 & $C_2$ & $G(16)$ &$[2,2,4,4]$ & $3$ &  $2$  & $4$  \\
\midrule
9 & $C_1$ & $D_4\times \mathbb{Z}/2$ &$[2,2,2,4]$& $3$ & $3$   &$6$  \\
9 & $C_2$ & $D_4\times \mathbb{Z}/2$ &$[2,2,2,2,2,2]$ & $5$ & $3$   &$6$  \\
\midrule
10 & $C_1$ & $(\mathbb{Z}/2)^4$ &$[2,2,2,2,2]$& $4$ &  $4$  &$10$ \\
10 & $C_2$ & $(\mathbb{Z}/2)^4$ &$[2,2,2,2,2]$ & $4$ &  $4$  &$10$ \\
\midrule
11 & $C_1$ & $(\mathbb{Z}/3)^2$ &$[3,3,3,3]$& $0$ &  $0$   &$0$  \\
11 & $C_2$ & $(\mathbb{Z}/3)^2$ &$[3,3,3,3]$ & $0$ & $0$   &$0$ \\
\midrule
12 & $C_1$ & $(\mathbb{Z}/2)^3$ & $[2,2,2,2,2]$& $4$ &  $3$  &$6$  \\
12 & $C_2$ & $(\mathbb{Z}/2)^3$ &$[2,2,2,2,2,2]$ & $5$ &  $3$  &$6$  \\
\bottomrule
\end{tabular}}
\end{center} 

There remains to understand the map
\[
\varphi^*: H^2(G,\mathbb{Z}/2) \longrightarrow H^2(\Gamma,\mathbb{Z}/2).
\]
 This can be understood observing that  $\mathbb{Z}/2$ is a trivial
$G, \Ga $-module.
 Hence the extensions classified by
$H^2(G,\mathbb{Z}/2)$ are central extensions. More precisely, $u\in H^2(G,\mathbb{Z}/2)$ can be interpreted as equivalence classes of central extensions
\[
1
\longrightarrow \mathbb{Z}/2
\longrightarrow \widetilde G
\longrightarrow G
\longrightarrow 1.
\]

\begin{lemma}\label{pullback}\cite[Thm 6.6.3, Exercise 6.6.5]{Weibel94}
Let $u\in H^2(G,\mathbb{Z}/2)$ and let $1\longrightarrow \mathbb{Z}/2 \longrightarrow \widetilde G\longrightarrow G\longrightarrow 1$ be the central extension corresponding to $u$. Then
\[
\begin{aligned}
\varphi^*(u)=0 &\quad\Longleftrightarrow\quad \text{the pullback extension } 1\to\mathbb{Z}/2\to \Gamma\times_G\widetilde G \to\Gamma\to1 \text{ splits} \\
&\quad\Longleftrightarrow\quad \exists\, s:\Gamma\to \Gamma\times_G\widetilde G \text{ such that } q\circ s=\operatorname{id}_\Gamma \\
&\quad\Longleftrightarrow\quad \exists\, \widetilde\varphi: \Gamma\to\widetilde G \text{ such that } p\circ\widetilde\varphi=\varphi.
\end{aligned}
\]
\end{lemma}
\begin{proof}
Let $u\in H^2(G,\mathbb{Z}/2)$ be represented by
\[
1
\longrightarrow \mathbb{Z}/2
\longrightarrow \widetilde G
\overset{p}{\longrightarrow} G
\longrightarrow 1.
\]
Pulling this extension back via $\varphi:\Gamma\longrightarrow G$ gives the fibre product
\[
\varphi^*\widetilde G=\Gamma\times_G\widetilde G=\{(\gamma,\widetilde g)\in\Gamma\times\widetilde G\mid\varphi(\gamma)=p(\widetilde g)\}.
\]
It fits into a central extension
\[
1
\longrightarrow\mathbb{Z}/2
\longrightarrow\varphi^*\widetilde G
\overset{q}\longrightarrow\Gamma
\longrightarrow 1.
\]
 Clearly  $\varphi^*(u)=0$ if and only if the pull-back extension splits.

\medskip

\medskip

We now explain why such a splitting is equivalent to a lift of $\varphi:\Gamma\to G$ to $\widetilde G$. Suppose first that a splitting $s:\Gamma\longrightarrow\Gamma\times_G\widetilde G$ exists. Since $q\circ s=\operatorname{id}_\Gamma$, the first coordinate of $s(\gamma)$ must be $\gamma$. Hence $s$ has the form $s(\gamma)=\bigl(\gamma,\widetilde\varphi(\gamma)\bigr)$ for some map $\widetilde\varphi:\Gamma\longrightarrow\widetilde G$. Since $s(\gamma)\in\Gamma\times_G\widetilde G$, the defining condition of the fiber product gives $\varphi(\gamma)=p(\widetilde\varphi(\gamma))$. Hence $p\circ\widetilde\varphi=\varphi$. Moreover, since $s$ is a group homomorphism, we have $s(\gamma_1\gamma_2)=s(\gamma_1)s(\gamma_2)$. Writing the two sides explicitly,
\[
\bigl(\gamma_1\gamma_2,\widetilde\varphi(\gamma_1\gamma_2)\bigr)=\bigl(\gamma_1\gamma_2,\widetilde\varphi(\gamma_1)\widetilde\varphi(\gamma_2)\bigr).
\]
Therefore
\[
\widetilde\varphi(\gamma_1\gamma_2)=\widetilde\varphi(\gamma_1)\widetilde\varphi(\gamma_2).
\]
Thus $\widetilde\varphi:\Gamma\longrightarrow\widetilde G$ is itself a group homomorphism. Hence a splitting of the pullback extension determines a lift $\widetilde\varphi$ satisfying $p\circ\widetilde\varphi=\varphi$.

\medskip

Conversely, suppose that there exists a group homomorphism $\widetilde\varphi:\Gamma\longrightarrow\widetilde G$ such that $p\circ\widetilde\varphi=\varphi$. Define $s:\Gamma\longrightarrow\Gamma\times_G\widetilde G$ by $ $ $s(\gamma)=\bigl(\gamma,\widetilde\varphi(\gamma)\bigr)$. Since $p(\widetilde\varphi(\gamma))=\varphi(\gamma)$, the pair $(\gamma,\widetilde\varphi(\gamma))$ indeed lies in the fiber product $\Gamma\times_G\widetilde G$. Moreover,
\[
s(\gamma_1\gamma_2)=\bigl(\gamma_1\gamma_2,\widetilde\varphi(\gamma_1\gamma_2)\bigr)=\bigl(\gamma_1\gamma_2,\widetilde\varphi(\gamma_1)\widetilde\varphi(\gamma_2)\bigr)\\
=s(\gamma_1)s(\gamma_2),
\]
so $s$ is a group homomorphism.
Finally,
\[
q(s(\gamma))=q\bigl(\gamma,\widetilde\varphi(\gamma)\bigr)=\gamma.
\]
Therefore
\[
q\circ s=\operatorname{id}_\Gamma.
\]

\end{proof}

\begin{lemma}\label{even-lifting}
Let
$\Gamma=\langle x_1,\ldots,x_r \mid x_i^{m_i}=1,\;x_1x_2\cdots x_r=1\rangle$ and let $\varphi:\Gamma\twoheadrightarrow G$ be
 an orbifold  surjective homomorphism, i.e., one
 satisfying $\varphi(x_i)=g_i$, $\operatorname{ord}(g_i)=m_i$. Assume that at least  one $m_i$ is even. Then
\[
\ker(\varphi^*:H^2(G;\mathbb{Z}/2)\longrightarrow H^2(\Gamma;\mathbb{Z}/2))=\bigcap_{i=1}^r\ker(\operatorname{res}_{\langle g_i\rangle}:H^2(G;\mathbb{Z}/2)\longrightarrow H^2(\langle g_i\rangle;\mathbb{Z}/2)) =   H^2(\langle x_i\rangle;\mathbb{Z}/2)).
\]
\end{lemma}

\begin{proof}
Let $u\in H^2(G;\mathbb{Z}/2)$:  as explained,  $u$ corresponds to an equivalence class of central extensions
\[
1
\longrightarrow\mathbb{Z}/2
\longrightarrow\widetilde G
\overset{p}{\longrightarrow}G
\longrightarrow1.
\]
 Hence we have an obvious inclusion of the left hand side in the right hand side, because if a central extension
of $\Ga$ splits, also the restriction to a subgroup of $\Ga$ splits.

Conversely, assume that we have a splitting of each restriction to the subgroup spanned by $x_i$, hence a mapping which sends $x_i$ to an element $ x'_i \in \tilde{\Ga}$ of order precisely
$m_i$.

We write the kernel as $\ker p=\langle z\rangle=\{1,z\}$, $z^2=1$ for some $z\in \widetilde G$. Thus we shall write the extension as
\[
1
\longrightarrow\langle z\rangle
\longrightarrow\widetilde G
\overset{p}{\longrightarrow}G
\longrightarrow1.
\]

Correspondingly, we have the extension
\[
1
\longrightarrow\langle z\rangle
\longrightarrow\widetilde {\Ga}
\overset{p}{\longrightarrow}\Ga
\longrightarrow1.
\]

The product $x'_1 \dots x'_r \in \langle z\rangle$, since its image in $\Ga$ is the identity.

If the product equals $1$, then we have found the desired splitting of $\Ga$, and the extension class is trivial, as we wanted to show.

If instead  $x'_1 \dots x'_r =z$, observe that, by assumption, at least one of the $m_i$ is even. Without loss of generality, suppose that $m_1$ is even. Since $p(z)=1$, the element $z x'_1$ is still a lift of $x_1$.

We now replace $ x'_1$  by $z x'_1$ and denote $x'_j$ by $y_j$.

We must verify that this modification does not destroy the power relation $(zy_1)^{m_1}=1$. Since the extension is central, $zy_1=x'_1z$. Therefore $(zy_1)^{m_1}=z^{m_1}y_1^{m_1}$. By construction, $y_1^{m_1}=1$ and  $m_1$ is even, hence  $z^{m_1}=1$. Consequently, $(zy_1)^{m_1}=1$. Thus replacing $x'_1 =y_1$ by $zx'_1$ preserves the
power relation. On the other hand, the product changes by exactly one factor of $z$:
\[
(zx'_1)x'_2\cdots x'_r=z(x'_1x'_2\cdots x'_r)=z\cdot z=z^2=1.
\]
Hence the modified lifts $zy_1,y_2,\ldots,y_r$ satisfy $(z y_1)^{m_1}=1$, $y_i^{\,m_i}=1$, $(i\ge2)$, $(zy_1)y_2\cdots y_r=1$.
Therefore they satisfy all thebdefining relations of $\Gamma$.

Hence  there exists a unique
homomorphism $\widetilde\varphi:\Gamma\longrightarrow\widetilde {\Ga}$
and we  have therefore proved the other inclusion
\[
\bigcap_{i=1}^r\ker(\operatorname{res}_{\langle g_i\rangle})\subseteq\ker\varphi^*.
\]
whence
\[
\ker\varphi^*=\bigcap_{i=1}^r\ker(\operatorname{res}_{\langle g_i\rangle}:H^2(G;\mathbb{Z}/2)\longrightarrow H^2(\langle g_i\rangle;\mathbb{Z}/2)).
\]

\end{proof}

\begin{proposition}\label{kervarphi}
Define
\[
R:H^2(G;\mathbb{Z}/2)\longrightarrow\bigoplus_{m_i \text{ : even }} H^2(\langle g_i \rangle; \mathbb{Z}/2),
\quad
R(u)=(res_{\langle g_i \rangle}(u))_{m_i \text{ : even }}.
\]
Then $\dim ker\varphi^*=\dim H^2(G,\mathbb{Z}/2) - \operatorname{rank} R$.

\end{proposition}
\begin{proof}
 Observe that $H^2 (\ZZ/m, \ZZ/2)=0$ if $m$ is odd, and $\ZZ/2$ if $m$ is even,
because any such central extension splits for $m$ odd, while for $m$ even there are two cases:
$ \ZZ/2m$ and $\ZZ/2 \oplus \ZZ/m$. By Lemma~\ref{even-lifting} and definition of $R$,  if there is some $m_i$ which is even, then
\[
\ker \varphi^* =\bigcap_{i=1}^r\ker(\operatorname{res}_{\langle g_i\rangle}:H^2(G;\mathbb{Z}/2)\longrightarrow H^2(\langle g_i\rangle;\mathbb{Z}/2))= \ker R.
\]

If instead all the $m_i$ are odd, then the right hand side is zero, and also in this case we have,
by the rank-nullity theorem, $\operatorname{rank} R + \dim\ker R = \dim H^2(G,\mathbb{Z}/2)$.
\end{proof}

Thus, we need only compute $\operatorname{rank} R.$

\medskip

In our situation  $H_j$ is a cyclic group of even order. Since $H^2(H_i;\mathbb{Z}/2)\cong\mathbb{Z}/2$, $\bigoplus_{j=1}^r H^2(H_i;\mathbb{Z}/2)\cong(\mathbb{Z}/2)^r$, and we may regard $R:H^2(G;\mathbb{Z}/2)\longrightarrow(\mathbb{Z}/2)^r$ as an ordinary $\mathbb{Z}/2$-linear map. Set $d:=\dim_{\mathbb{Z}/2}H^2(G;\mathbb{Z}/2)$, and choose a basis $u_1,\ldots,u_d$ of $H^2(G;\mathbb{Z}/2)$. Define $$M=\bigl(\operatorname{res}_{H_i}(u_j)\bigr)_{\substack{1\le i\le r\\1\le j\le d}}.$$

Let $u\in H^2(G;\mathbb{Z}/2)$. Write $u=a_1u_1+\cdots+a_du_d$, $a_i\in\mathbb{Z}/2$. By linearity, $M$ is exactly the matrix representing $R$ with respect to the chosen basis of $H^2(G;\mathbb{Z}/2)$.

\medskip

Let us explain how to compute it.

\begin{definition}[Bockstein homomorphism]
Consider the short exact sequence of trivial $G$-modules $0\longrightarrow\mathbb{Z}/2\overset{\iota}{\longrightarrow}\mathbb Z/4\overset{\rho}{\longrightarrow}\mathbb{Z}/2\longrightarrow 0$, where $\iota(1)=2$ and $\rho$ is reduction modulo $2$. The associated long exact sequence in group cohomology contains the connecting homomorphism $\beta_G:H^1(G;\mathbb{Z}/2)\longrightarrow H^2(G;\mathbb{Z}/2)$, which  is called the Bockstein homomorphism.
\end{definition}

\begin{remark}
In the  long exact sequence of  group cohomology
\[
\begin{aligned}
\cdots
\longrightarrow
H^1(G;\mathbb{Z}/2)
\overset{\iota_*}{\longrightarrow}
H^1(G;\mathbb Z/4)
\overset{\rho_*}{\longrightarrow}
H^1(G;\mathbb{Z}/2)
\overset{\beta}{\longrightarrow}
H^2(G;\mathbb{Z}/2)
\longrightarrow
\cdots,
\end{aligned}
\]
since the coefficient modules are trivial, we have $H^1(G;\mathbb{Z}/2)\cong\operatorname{Hom}(G,\mathbb{Z}/2)$ and $H^1(G;\mathbb Z/4)\cong\operatorname{Hom}(G,\mathbb Z/4)$. Under these identifications, the map $\rho_*:H^1(G;\mathbb Z/4)\longrightarrow H^1(G;\mathbb{Z}/2)$ is given by composition with $\rho$:
 $$\rho: \widetilde\chi\longmapsto\rho\circ\widetilde\chi.$$

Therefore, by the exactness of the long exact sequence, $\ker\beta=\operatorname{Im}\rho_*$.

Consequently, for $\chi\in H^1(G;\mathbb{Z}/2)=\operatorname{Hom}(G,\mathbb{Z}/2)$, we have

\[
\beta(\chi)=0 \iff \chi:G\to\mathbb{Z}/2 \text{ lifts to a homomorphism } \widetilde\chi:G\to\mathbb Z/4,
\]
Equivalently,
\[
\beta(\chi)\neq0 \iff \chi \text{ does not admit a }\mathbb Z/4\text{-valued lift}.
\]

\end{remark}

\begin{remark}
Naturality of the connecting homomorphism gives a commutative diagram
\[
\begin{tikzcd}
H^1(G;\mathbb{Z}/2) \arrow[r, "\beta_G"] \arrow[d, "\operatorname{res}_H"'] & H^2(G;\mathbb{Z}/2) \arrow[d, "\operatorname{res}_H"] \\
H^1(H;\mathbb{Z}/2) \arrow[r, "\beta_H"'] & H^2(H;\mathbb{Z}/2).
\end{tikzcd}
\]
Hence, for every $\chi\in H^1(G;\mathbb{Z}/2)$, we have
\[
\operatorname{res}_H(\beta_G(\chi))=\beta_H(\chi|_H).
\]
\end{remark}

As observed in Proposition \ref{kervarphi} $H^2(\mathbb{Z}/m;\mathbb{Z}/2)\cong\mathbb{Z}/2$ for $m$ even, and is zero for $m$ odd.

\begin{lemma}\label{restrict}
 Let $\beta_G(\chi)\in H^2(G;\mathbb{Z}/2)$ for $\chi \in H^1(G,\mathbb{Z}/2)$. Let $H=\langle g\rangle\cong \mathbb{Z}/m$, where $m$ is even. Then
\[
\operatorname{res}_{H}\beta_G(\chi)=
\begin{cases}
\chi(g),&m\equiv 2\pmod 4,\\[2mm]
0,&4\mid m,
\end{cases}
\]
after identifying $H^2(H;\mathbb{Z}/2)\cong\mathbb{Z}/2$.
\end{lemma}

\begin{proof}
Let $H=\mathbb{Z}/m=\langle g\mid g^m=1\rangle$, where $m$ is even, and let $\psi:=\chi\mid_H :H\longrightarrow\mathbb{Z}/2$ be a character. If $\psi(g)=0$, then $\psi=0$, and hence $\beta_H(\psi)=0$. Assume now that $\psi(g)=1$. Then $\psi$ is the unique nonzero character of $\mathbb{Z}/m$. Recall $\beta_H(\psi)=0$ if and only if there exists a homomorphism $\widetilde\psi:\mathbb{Z}/m\longrightarrow\mathbb Z/4$ such that $\rho \circ \widetilde\psi=\psi$. Since $\psi(g)=1$, the element $\widetilde\psi(g)$ must lie in $\rho^{-1}(1)=\{1,3\}\subset\mathbb Z/4$. Both $1$ and $3$ have order $4$ in $\mathbb Z/4$. On the other hand, since $g^m=1$, any homomorphism $\widetilde\psi:\mathbb{Z}/m\to\mathbb Z/4$ must satisfy $m\,\widetilde\psi(g)=0 \in \mathbb Z/4.$ If $\widetilde\psi(g)$ is odd, then it is invertible modulo $4$. Hence $m\,\widetilde\psi(g)\equiv0\pmod4$ holds if and only if $m\equiv0\pmod4$. Therefore,
\[
\beta_H(\psi)=0 \iff 4\mid m.
\]
Since $m$ is even, the only alternative is $m\equiv2\pmod4$, and in this case $\beta_H(\psi)\neq0$. Thus
\[
\beta_H(\psi)=
\begin{cases}
1,&m\equiv2\pmod4,\\
0,&4\mid m.
\end{cases}
\]
Combining the cases $\psi(g)=0$ and $\psi(g)=1$, we obtain
\[
\beta_H(\psi)=
\begin{cases}
\psi(g),&m\equiv2\pmod4,\\
0,&4\mid m.
\end{cases}
\]
Since $\psi=\chi|_H$, we have $\psi(g)=\chi(g)$. Since $\operatorname{res}_H(\beta_G(\chi))=\beta_H(\chi|_H)$, we obtain
\[
\operatorname{res}_{H}\beta_G(\chi)=
\begin{cases}
\chi(g),&m\equiv2\pmod4,\\[2mm]
0,&4\mid m.
\end{cases}
\]
\end{proof}

 Consider a basis $\chi_1, \cdots, \chi_b$ of $H^1(G;\mathbb{Z}/2)$ and $H_i=\langle g_i\rangle\cong \mathbb{Z}/m_i$, $i=1,\cdots r$. Hence the part of the restriction matrix is obtained directly from the character values:
\[
B_{ij}=
\begin{cases}
\chi_j(g_i),&m_i\equiv2\pmod4,\\[2mm]
0,&4\mid m_i.
\end{cases}
\]

\begin{remark}
The classes $\beta_G(\chi)$ alone are not sufficient to determine the full rank of $R$. For example, if $H=\langle g\rangle\cong \mathbb{Z}/4$, then for every $\chi\in H^1(G;\mathbb{Z}/2)$ we have $\operatorname{res}_{H}\beta_G(\chi)=0$. Therefore, additional cohomology classes are needed to determine the remaining part of the matrix.
\end{remark}

\bigskip

\section{Clifford algebra, Binary octahedral group}\label{Clifford algebra}

\medskip

\noindent\textbf{3.1. Clifford algebra.}\label{clifford} \\ \\
 We collect some results about the real Clifford algebras which are needed in this paper. See \cite[\S 1--3]{AtiyhaBottShaprio64} for more details about Clifford algebras. Note that the sign convention used there differs from ours. 
 
\begin{definition}
Let $V=\mathbb R^n$ with its standard Euclidean inner product $\langle\ ,\ \rangle$, and let $e_1,\ldots,e_n$ be the standard orthonormal basis. The Clifford algebra $Cl(V)$ is defined as the quotient
\[
Cl(V)=T(V)\Big/\langle v\otimes v-\langle v,v\rangle\cdot 1\mid v\in V\rangle,
\]
of the tensor algebra of $V$
\[
T(V)=\mathbb R\oplus V\oplus (V\otimes V)\oplus (V^{\otimes 3})\oplus\cdots
\]

\end{definition}
\begin{remark}
There are several important relations that are frequently used in Clifford algebras.
\begin{enumerate}
\item If $u$ is a unit vector, then $u^2=1$, and hence $u^{-1}=u$.
\item If $u\neq0$, then $u^{-1}=\frac{u}{\langle u,u\rangle}$.
\item $uv+vu=2\langle u,v\rangle$
\end{enumerate}
\end{remark}

\begin{definition}
The Clifford algebra is naturally $\mathbb Z/2$-graded:
\[
Cl(V)=Cl^0(V)\oplus Cl^1(V),
\]
where $Cl^0(V)$ is generated by products of an even number of vectors.

\end{definition}

\begin{definition}
The Pin group is generated by unit vectors in $Cl(V)^\times$, it is  the  Clifford multiplicative subgroup   of all elements of the form $v_1v_2\cdots v_k$ where each $v_i\in V$ is of unit length $<v_i,v_i>=1$
\end{definition}

\begin{definition}
The Spin group is the even part in $Pin(V)$, denoted, $\operatorname{Spin}(V):=Pin(V)\cap Cl^0(V)$.
\end{definition}
For example, for $V=\mathbb R^5$, we write $\operatorname{Spin}(5):=\operatorname{Spin}(\mathbb R^5)$. Also, there is an exact sequence $$1\longrightarrow\{\pm1\}\longrightarrow\operatorname{Spin}(5)\overset{p}{\longrightarrow}SO(5)\longrightarrow1.$$

The alternating group  $\frak{A}_5$
acts on $\mathbb R^5$ by permuting the standard basis: $\sigma(e_i)=e_{\sigma(i)}$. Since permutations preserve the Euclidean inner product, this gives an embedding $\mathfrak{A}_5\subset O(5)$. Moreover, the determinant of the permutation matrix associated with $\sigma\in \mathfrak{S}_5$ is $\det(\sigma)=\operatorname{sgn}(\sigma)$. Since every element of $\mathfrak{A}_5$ is an even permutation, $\operatorname{sgn}(\sigma)=1$ for every $\sigma\in \mathfrak{A}_5$. Hence $\mathfrak{A}_5\subset SO(5)$.

\medskip

Let $\widetilde{\mathfrak{A}_5}:=p^{-1}(\mathfrak{A}_5)\subset\operatorname{Spin}(5)$. Restricting the double covering $p:\operatorname{Spin}(5)\to SO(5)$ to this inverse image gives a central extension
\[
1
\longrightarrow\{\pm1\}
\longrightarrow\widetilde{\mathfrak{A}_5}
\longrightarrow \mathfrak{A}_5
\longrightarrow1.
\]
Let $u\in V$ be a unit vector. The reflection through the
hyperplane $u^\perp=\{v\in V\mid \langle u,v\rangle=0\}$ is given by $s_u(v)=v-2\langle v,u\rangle u$.
\begin{lemma}
This reflection can be written in the Clifford algebra as $$s_u(v)=-uvu^{-1}=-uvu.$$
\end{lemma}
\begin{proof}
$uv+vu=2\langle u,v\rangle$, and therefore $vu=2\langle u,v\rangle-uv$. Multiplying on the left by $u$, we obtain $uvu=u\bigl(2\langle u,v\rangle-uv\bigr)=2\langle u,v\rangle u-u^2v$. Since $u$ is a unit vector, $u^2=1$, and hence $uvu=2\langle u,v\rangle u-v$.
Since $u^{-1}=u$, we  get the equality.
\end{proof}

For $i\neq j$, define $\alpha_{ij}:=\frac{e_i-e_j}{\sqrt2}$. First observe that
\[
\begin{aligned}
\|\alpha_{ij}\|^2
=\langle\frac{e_i-e_j}{\sqrt2},\frac{e_i-e_j}{\sqrt2}\rangle=\frac12(\langle e_i,e_i\rangle-2\langle e_i,e_j\rangle+\langle e_j,e_j\rangle)=1.
\end{aligned}
\]
Thus $\alpha_{ij}$ is a unit vector, and $\alpha_{ij}^2=1$. We now compute explicitly the reflection $s_{\alpha_{ij}}(v)=-\alpha_{ij}v\alpha_{ij}$.

\begin{lemma}
\[
s_{\alpha_{ij}}:
\begin{cases}
e_i\longmapsto e_j,\\
e_j\longmapsto e_i,\\
e_k\longmapsto e_k,
& k\neq i,j.
\end{cases}
\]
We may regard the transposition $(ij)$ as the reflection $s_{\alpha_{ij}}$ in the Clifford algebra.
\end{lemma}
\begin{proof}

First, for $e_i$,
\[
\begin{aligned}
s_{\alpha_{ij}}(e_i)
=-\alpha_{ij}e_i\alpha_{ij}=-\frac12(e_i-e_j)e_i(e_i-e_j).
\end{aligned}
\]
Since $e_i^2=1$ and, for $i\neq j$, $e_ie_j=-e_je_i$, we have
\[
\begin{aligned}
(e_i-e_j)e_i
=e_i^2-e_je_i=1+e_ie_j,
\end{aligned}
\]
\[
\begin{aligned}
(e_i-e_j)e_i(e_i-e_j)
=(1+e_ie_j)(e_i-e_j)=e_i-e_j+e_ie_je_i-e_ie_je_j,
\end{aligned}
\]
\[
\begin{aligned}
e_ie_je_i=e_i(e_je_i)=e_i(-e_ie_j)=-e_i^2e_j=-e_j,
\end{aligned}
\]
\[
e_ie_je_j=e_ie_j^2=e_i.
\]
Thus
\[
\begin{aligned}
(e_i-e_j)e_i(e_i-e_j)
=e_i-e_j-e_j-e_i=-2e_j.
\end{aligned}
\]
Hence
$\alpha_{ij}e_i\alpha_{ij}=-e_j$
and consequently
$
s_{\alpha_{ij}}(e_i)=e_j.
$

Similarly,
\[
\begin{aligned}
s_{\alpha_{ij}}(e_j)
=-\alpha_{ij}e_j\alpha_{ij}=-\frac12(e_i-e_j)e_j(e_i-e_j).
\end{aligned}
\]
Since
\[
(e_i-e_j)e_j=e_ie_j-e_j^2=e_ie_j-1,
\]
we obtain
\[
\begin{aligned}
(e_i-e_j)e_j(e_i-e_j)
=(e_ie_j-1)(e_i-e_j)=e_ie_je_i-e_ie_je_j-e_i+e_j=-e_j-e_i-e_i+e_j=-2e_i.
\end{aligned}
\]
Therefore $s_{\alpha_{ij}}(e_j)=e_i.$

Finally, if $k\neq i,j$, then
\[
\langle\alpha_{ij},e_k\rangle=0.
\]
Hence, by the Clifford relation,
\[
\alpha_{ij}e_k=-e_k\alpha_{ij},
\]
and therefore $s_{\alpha_{ij}}(e_k)
=-\alpha_{ij}e_k\alpha_{ij}=-(-e_k\alpha_{ij})\alpha_{ij}=e_k\alpha_{ij}^2=e_k.$
\end{proof}

Let $(ijk)$ be a $3$-cycle. Since $\alpha_{ij}$ lifts the reflection $(ij)$ and $\alpha_{jk}$ lifts the reflection $(jk)$, their product $\alpha_{ij}\alpha_{jk}$ belongs to $\operatorname{Spin}(5)$ that lifts the reflection $(ij)(jk)=(ijk)$. Thus we may define $\widetilde{(ijk)}:=\alpha_{ij}\alpha_{jk}$.

\medskip

It remains to check that this lift has order exactly $3$, not $6$. Set $p:=\alpha_{ij}$, $q:=\alpha_{jk}$. We compute
\[
\begin{aligned}
\langle p,q\rangle
=\langle\frac{e_i-e_j}{\sqrt2},\frac{e_j-e_k}{\sqrt2}\rangle=-\frac12.
\end{aligned}
\]
The Clifford relation gives $pq+qp=2\langle p,q\rangle=-1$. Hence $qp=-1-pq$. Set $u:=pq$. Then
\[
\begin{aligned}
u^2=pqpq=p(qp)q=p(-1-pq)q=-pq-p^2q^2.
\end{aligned}
\]
Since $p$ and $q$ are unit vectors, $p^2=q^2=1$, and hence $u^2+u+1=0$. Multiplying by $u-1$, we obtain $u^3=1$. Since the image of $u$ in $\mathfrak{A}_5$ is the nontrivial $3$-cycle
$(ijk)$, we have $u\neq1$. Then, $\operatorname{ord}(u)=3$. Thus $\alpha_{ij}\alpha_{jk}$ is an order-preserving lift of the $3$-cycle $(ijk)$. Note that the other lift is $-u$. The element $-u$ has order $6$.

\bigskip

\noindent\textbf{3.2. The binary octahedral group.}\label{oct}

\bigskip

We see that an element of order $4$ in  lifts to an element of order $8$ in the binary octahedral group $2O$,
for details see \cite[\S 11.5.4]{Voight21}.

Recall that $\mathbb{H}=\mathbb{R}+\mathbb{R}i+\mathbb{R}j+\mathbb{R}k$, $i^2=j^2=k^2=ijk=-1$, is the Hamilton quaternion algebra.  For $q=a+bi+cj+dk\in\mathbb{H}$, its reduced norm is defined by $\operatorname{nrd}(q)=q\overline q=a^2+b^2+c^2+d^2$, where $\overline q=a-bi-cj-dk$ is the conjugate in $\mathbb H$.

 Let $\mathbb{H}^1=\{q\in\mathbb{H} \mid \operatorname{nrd}(q)=1\}$ denote the group of unit quaternions. Under the standard double covering $\mathbb{H}^1\cong \operatorname{Spin}(3)\longrightarrow SO(3)$, the octahedral rotation group $\mathfrak{S}_4\subset SO(3)$ pulls back to the binary octahedral group $2O\subset \mathbb{H}^1$, which has order $48$. Thus one obtains a central extension
\[
1 \longrightarrow\{\pm1\} \longrightarrow2O \longrightarrow \mathfrak{S}_4 \longrightarrow1.
\]

The elements
\[
\frac{\pm1\pm i}{\sqrt2},
\quad
\frac{\pm1\pm j}{\sqrt2},
\quad
\frac{\pm1\pm k}{\sqrt2}
\]
map to rotations of order $4$ in the octahedral group. For example, let $\alpha=\frac{1+i}{\sqrt2}\in 2O$. Since $i^2=-1$, we have $\alpha^2=\frac{(1+i)^2}{2}=i$, and hence $\alpha^4=-1$, $\alpha^8=1$. Therefore $\operatorname{ord}(\alpha)=8$. Its image in $\mathfrak{S}_4$ has order $4$, since
\[
\alpha^4=-1\in\ker(2O\to \mathfrak{S}_4), \quad \alpha^2=i\notin\{\pm1\}.
\]
Thus an element of order $4$ in $\mathfrak{S}_4$ lifts to an element of order $8$ in the binary octahedral group $2O$.

\medskip

Recall that the orientation-preserving rotation group of the cube, equivalently of the octahedron, is isomorphic to $\mathfrak{S}_4$. Thus we may regard $\mathfrak{S}_4\subset SO(3)$. Consider the standard spin double covering
\[
1
\longrightarrow\{\pm1\}
\longrightarrow\operatorname{Spin}(3)
\overset{\operatorname{pr}}{\longrightarrow}SO(3)
\longrightarrow1.
\]
Using the standard identification $\operatorname{Spin}(3)\cong SU(2)$, the inverse image of the octahedral rotation group is the binary octahedral group $2O:=\operatorname{pr}^{-1}(\mathfrak{S}_4).$

\bigskip

\begin{remark}

In the following  sections, we compute, for each of the  cases,  first  the number of $G$-invariant theta characteristics,  presenting however  our computation here only for the 7 cases whose parity is not known,  and  omitting the remaining ones.

Then we use the Hochschild-Serre spectral sequence to describe the cohomology and the invariant divisors on a surface isogenous to a product of unmixed type, using the previous result by Dolgachev for the case of curves. More precisely, we compute the group $\operatorname{Pic}(C)^G$ appearing in \eqref{eq:exact} by exploiting  its 2-torsion subgroup $\operatorname{Pic}(C)^G[2]$ and  the number of $G$-invariant theta characteristics that are computed in the forthcoming Sections and then show how this leads to the parity determination.

\end{remark}

\bigskip

\section{$G=\mathfrak{A}_5$, $g(C_1)=4$, $g(C_2)=21$, $T_1=[2,5,5]$, $T_2=[3,3,3,3]$; even fake quadrics}

We first recall  the non-split central double cover of $\mathfrak{A}_5$ (cf. \cite[Exercise 6.9.1(2)]{Weibel94}). For the trivial $\mathfrak{A}_5$-module $\mathbb{Z}/2$, the group $H^2(\mathfrak{A}_5,\mathbb{Z}/2)$ classifies central extensions of $\mathfrak{A}_5$ by $\mathbb Z/2$, and consists of two elements. The zero class corresponds to the split extension
\[
1
\longrightarrow\mathbb Z/2
\longrightarrow\mathbb Z/2\times \mathfrak{A}_5
\longrightarrow \mathfrak{A}_5
\longrightarrow1,
\]
and the nonzero class corresponds to the non-split central double cover
\[
1
\longrightarrow\mathbb Z/2
\longrightarrow2.\mathfrak{A}_5
\overset{p}{\longrightarrow}\mathfrak{A}_5
\longrightarrow1.
\]
Here $2.\mathfrak{A}_5$ is isomorphic to $\mathrm{SL}_2(\mathbb{F}_5)$, while $\mathfrak{A}_5\cong\operatorname{PSL}_2(\mathbb{F}_5)$. Thus the above extension may be identified with
\[
1
\longrightarrow\{\pm I\}
\longrightarrow\mathrm{SL}_2(\mathbb{F}_5)
\longrightarrow\operatorname{PSL}_2(\mathbb{F}_5)
\longrightarrow1.
\]
 The pullback of the non-split central extension
along $\varphi_1:\Gamma_1\to \mathfrak{A}_5$ is
\[
1
\longrightarrow\mathbb Z/2
\longrightarrow\Gamma_1\times_{\mathfrak{A}_5}2.\mathfrak{A}_5
\overset{q}\longrightarrow\Gamma_1
\longrightarrow1,
\]
where $\Gamma_1\times_{\mathfrak{A}_5}2.\mathfrak{A}_5=\{(\gamma,\widetilde g)\in \Gamma_1\times2.\mathfrak{A}_5 \mid\varphi_1(\gamma)=p(\widetilde g)\}$ and $q$ is a natural projection.

\begin{lemma}\label{A_5 involution}
Every lift of a nontrivial involution of $\mathfrak{A}_5$ to $2.\mathfrak{A}_5$ has order $4$.
\end{lemma}
\begin{proof}
Identify $2.\mathfrak{A}_5\cong\mathrm{SL}_2(\mathbb{F}_5)$, $\mathfrak{A}_5\cong\operatorname{PSL}_2(\mathbb{F}_5)$, and let $p:\mathrm{SL}_2(\mathbb{F}_5)\to\operatorname{PSL}_2(\mathbb{F}_5)$ be the natural quotient map. Choose a lift $A\in\mathrm{SL}_2(\mathbb{F}_5)$ of $h$, so that $p(A)=h$. Since $p(A^2)=p(A)^2=h^2=1$, $A^2\in\ker p=\{\pm I\}$. Suppose that $A^2=I$. Since the characteristic of $\mathbb{F}_5$ is not $2$, the polynomial $t^2-1=(t-1)(t+1)$ has distinct roots and every matrix satisfying $A^2=I$ is diagonalizable, with eigenvalues $\pm 1$. Since $A\in\mathrm{SL}_2(\mathbb{F}_5)$, we also have $\det A=1$. Therefore, $A=I$ or $A=-I$, hence $A \in\ker p=\{\pm I\}$ and $p(A)=1$, contradicting $p(A)=h\neq1$. Therefore
 $A^2=-I$. It follows that $A^4=(-I)^2=I$. Hence every lift of a nontrivial involution of $\mathfrak{A}_5$ to $2.\mathfrak{A}_5$ has order $4$. Note also that $h$ has two lifts, $A$ and $-A$.
\end{proof}

\begin{proposition}  In the case $G=\mathfrak{A}_5$, $g(C_1)=4$, $g(C_2)=21$, $T_1=[2,5,5]$, $T_2=[3,3,3,3]$, we have
$$\operatorname{Pic}(C_1)^G[2] \cong 0 \quad\text{and}\quad \#\Theta(C_1)^G=1.$$
\end{proposition}

\begin{proof}
For $C_1$, a generating vector is $(g_1,g_2,g_3)=\bigl((24)(35),(21345),(12345)\bigr)$ with orders $2,5,5$.

 For further details, we refer the reader to \cite[\S~6]{BCG08}.

We first recall the existence of an $\mathfrak{A}_5$-invariant theta
characteristic for such an $\mathfrak{A}_5$-action, i.e.,
\[
\Theta(C_1)^{\mathfrak{A}_5}\neq\varnothing.
\]
\
Indeed, for the $\mathfrak{A}_5$-action of signature $[2,5,5]$ the curve $C_1$
is called Bring's curve.

As also shown in \cite{odd},  the curve $C_1$ does not admit any $\mathfrak{A}_5$-linearized theta characteristic and  is the unique genus 4 curve with automorphism group containing $\mathfrak{A}_5$, canonically embedded as: $$ C_1 = \{ (x_1, \dots, x_5) \in \PP^4 | \sum_i x_i=0 , \sum_i x_i^2=0, \sum_i x_i^3=0  \} .$$ Burns constructed an $\mathfrak{A}_5$-invariant theta characteristic on $C_1$. Braden and Disney-Hogg  \cite[Corollary~5.10]{BradenDisneyHogg24} later proved that Bring's curve has a unique theta characteristic invariant under its full automorphism group $\mathfrak{S}_5$, and identified it with the theta characteristic constructed by Burns.

\medskip

Let
\[
\Gamma_1=\langle x_1,x_2,x_3\mid x_1^2=x_2^5=x_3^5=1,\ x_1x_2x_3=1\rangle
\]
and let $\varphi_1:\Gamma_1\twoheadrightarrow \mathfrak{A}_5$, $x_i\mapsto g_i$.  To compute $\#\Theta(C_1)^G$, by Proposition~\ref{dim} we compute
\[
\dim H^1(C,\mathbb{Z}/2)^G=\dim H^1(\Gamma,\mathbb{Z}/2)-\dim H^1(G,\mathbb{Z}/2)+\dim\ker\varphi^*.
\]
By Table~\ref{dim table}, we know that $\dim H^1(\Gamma_1,\mathbb{Z}/2)=\dim H^1(\mathfrak{A}_5,\mathbb{Z}/2)=0$. Thus,
\[
H^1(C_1,\mathbb{Z}/2)^{\mathfrak{A}_5}\cong\ker(H^2(\mathfrak{A}_5,\mathbb{Z}/2)\overset{\varphi_1^*}{\longrightarrow}H^2(\Gamma_1,\mathbb{Z}/2)).
\]
Claim that ker$\phi_1^*=0$. Suppose, for contradiction, that such a lift $\widetilde\varphi_1:\Gamma_1\to2.\mathfrak{A}_5$ exists. Set $\widetilde g_i:=\widetilde\varphi_1(x_i)$. Since $p(\widetilde g_i)=\varphi_1(x_i)=g_i$, each $\widetilde g_i$ is a lift of $g_i$. However, it is not enough to choose arbitrary preimages of the $g_i$. Since $\widetilde\varphi_1$ is a group homomorphism, the defining relations of $\Gamma_1$ must be preserved. In particular,
we must have $\widetilde g_1^2=1$, $\widetilde g_2^5=1$, $\widetilde g_3^5=1$, and $\widetilde g_1\widetilde g_2\widetilde g_3=1$. We shall see that the first relation alone is already impossible.
Since $g_1=(24)(35)$ is a nontrivial involution, by Lemma~\ref{A_5 involution}, every lift $\widetilde g_1\in2.\mathfrak{A}_5$ satisfies $\widetilde g_1^2=-1$, and therefore has order $4$. On the other hand, the corresponding generator $x_1\in\Gamma_1$ satisfies the orbifold relation $x_1^2=1$. If a lift $\widetilde\varphi_1:\Gamma_1\to 2.\mathfrak{A}_5$ existed, then $\widetilde g_1=\widetilde\varphi_1(x_1)$ would have to satisfy
\[
\widetilde g_1^2=\widetilde\varphi_1(x_1)^2=\widetilde\varphi_1(x_1^2)=\widetilde\varphi_1(1)=1.
\]
But every lift of $g_1$ satisfies $\widetilde g_1^2=-1\neq1$. This is a contradiction.

Therefore no such homomorphism $\widetilde\varphi_1:\Gamma_1\to 2.\mathfrak{A}_5$ can exist.
In other words, $\varphi_1:\Gamma_1\to \mathfrak{A}_5$ does not lift to $2.\mathfrak{A}_5$.

 It follows that $\varphi_1^*(u)\neq0$, where $u\neq 0 \in H^2(\mathfrak{A}_5,\mathbb{Z}/2)\cong\mathbb{Z}/2$.  Since $H^2(\mathfrak{A}_5,\mathbb{Z}/2)$ is one-dimensional, $\varphi_1^*$ is injective. We conclude that $H^1(C_1,\mathbb{Z}/2)^{\mathfrak{A}_5}=0$. Burns constructed an invariant theta characteristic, hence
$
\#\Theta(C_1)^{\mathfrak{A}_5}=\#H^1(C_1,\mathbb{Z}/2)^{\mathfrak{A}_5}=1.$
\end{proof}

\begin{proposition}  In the case $G=\mathfrak{A}_5$, $g(C_1)=4$, $g(C_2)=21$, $T_1=[2,5,5]$, $T_2=[3,3,3,3]$, we have
$$\operatorname{Pic}(C_2)^G[2] \cong \mathbb{Z}/2 \quad\text{and}\quad \#\Theta(C_2)^G=2.$$
\end{proposition}
\begin{proof}
For $C_2$, a generating vector is $(g_1,g_2,g_3,g_4)=\bigl((123),(345),(432),(215)\bigr)$ with orders $3,3,3,3$.  For further details, we refer the reader to \cite[\S~6]{BCG08}. Let
\[
\Gamma_2=\langle y_1,y_2,y_3,y_4\mid y_1^3=y_2^3=y_3^3=y_4^3=1,\ y_1y_2y_3y_4=1\rangle
\]
and let $\varphi_2:\Gamma_2\twoheadrightarrow \mathfrak{A}_5$, $y_i\mapsto g_i.$
 To compute  $\#\Theta(C_2)^G$, by Proposition~\ref{dim} we compute
\[
\dim H^1(C,\mathbb{Z}/2)^G=\dim H^1(\Gamma,\mathbb{Z}/2)-\dim H^1(G,\mathbb{Z}/2)+\dim\ker\varphi^*.
\]
By Table~\ref{dim table}, we know that $\dim H^1(\Gamma_2,\mathbb{Z}/2)=\dim H^1(\mathfrak{A}_5,\mathbb{Z}/2)=0$. Thus,
\[
H^1(C_2,\mathbb{Z}/2)^{\mathfrak{A}_5}\cong\ker(H^2(\mathfrak{A}_5,\mathbb{Z}/2)\overset{\varphi_2^*}{\longrightarrow}H^2(\Gamma_2,\mathbb{Z}/2)).
\]

\medskip

Now, we prove that $\varphi_2^*=0$. To prove this, we construct such a lift explicitly using the Clifford algebra(see Section~\ref{clifford}). We write $g_1=(123)=(12)(23)$, $g_2=(345)=(34)(45)$, $g_3=(432)=(43)(32)$, $g_4=(215)=(21)(15)$. Accordingly, choose
\[
\widetilde g_1=\alpha_{12}\alpha_{23},
\quad
\widetilde g_2=\alpha_{34}\alpha_{45},
\quad
\widetilde g_3=\alpha_{43}\alpha_{32},
\quad
\widetilde g_4=\alpha_{21}\alpha_{15},
\]
where $$\alpha_{ji}=\frac{e_j-e_i}{\sqrt2}=-\alpha_{ij}.$$
Each $\widetilde g_i$ is an order-$3$ lift of $g_i$. Hence $\widetilde g_i^3=1$, for $i=1,2,3,4$. Thus the four relations $x_i^3=1$ are already preserved. The only remaining defining relation of $\Gamma_2$ is $x_1x_2x_3x_4=1$. Therefore we must verify that
\[
\widetilde g_1\widetilde g_2\widetilde g_3\widetilde g_4=1.
\]
For convenience, set
\[
a:=\alpha_{12},
\quad
b:=\alpha_{23},
\quad
c:=\alpha_{34},
\quad
d:=\alpha_{45},
\quad
e:=\alpha_{15}.
\]
Since $\alpha_{ji}=-\alpha_{ij}$, we have
\[
\widetilde g_1=ab,\quad \widetilde g_2=cd,\quad \widetilde g_3=(-c)(-b)=cb,\quad \widetilde g_4=(-a)e=-ae,
\]
\[
\widetilde g_1\widetilde g_2\widetilde g_3\widetilde g_4=(ab)(cd)(cb)(-ae)=-abcdcbae.
\]
We shall simplify this expression step by step. We repeatedly use $uvu=2\langle u,v\rangle u-v$.
Since $\langle c,d\rangle=\frac12\langle e_3-e_4, e_4-e_5\rangle=-\frac12$ and
$cdc=2\langle c,d\rangle c-d=-c-d,$ we have
$$-abcdcbae=-ab(cdc)bae=-ab(-c-d)bae=ab(c+d)bae.$$
Likewise, $bcb=2\langle b,c\rangle b-c=-b-c$, $bdb=-d$, $b(c+d)b=-b-c-d$, hence
\[
ab(c+d)bae=a\bigl(b(c+d)b\bigr)ae=a(-b-c-d)ae=-a(b+c+d)ae.
\]
Also, $aba=2\langle a,b\rangle a-b=-a-b$,  $aca=-c$, $ada=-d$,
\[
a(b+c+d)a=aba+aca+ada=(-a-b)+(-c)+(-d)=-a-b-c-d,
\]
hence $$\widetilde g_1\widetilde g_2\widetilde g_3\widetilde g_4=(a+b+c+d)e=e^2=\alpha_{15}^2=1.$$
We have verified that $\widetilde g_i^3=1$ for $i=1,2,3,4$ and $\widetilde g_1\widetilde g_2\widetilde g_3\widetilde g_4=1$. Thus the elements $\widetilde g_1$, $\widetilde g_2$, $\widetilde g_3$, $\widetilde g_4$ satisfy all the defining relations in $\Gamma_2=\langle x_1,x_2,x_3,x_4\mid x_i^3=1,\;x_1x_2x_3x_4=1\rangle$. Hence, the assignment $x_i\longmapsto\widetilde g_i$ gives a group homomorphism $\widetilde\varphi_2:\Gamma_2\longrightarrow2\cdot \mathfrak{A}_5$. By construction, $p(\widetilde g_i)=g_i=\varphi_2(x_i)$ for every $i$, and hence $p\circ\widetilde\varphi_2=\varphi_2$. Then, $\varphi_2$ admits a lifting. Thus, $\varphi_2^*=0$, that is,
\[
H^1(C_2,\mathbb{Z}/2)^{\mathfrak{A}_5}\cong\ker\varphi_2^*=H^2(\mathfrak{A}_5,\mathbb{Z}/2)\cong\mathbb{Z}/2.
\]
 Proposition~\ref{admit theta} gives an $\mathfrak{A}_5$-invariant theta characteristic, thus $$\#\Theta(C_2)^G=\#H^1(C_2,\mathbb{Z}/2)^{\mathfrak{A}_5}=2.$$
\end{proof}

\medskip

 Next, we determine $\operatorname{Pic}(C_i)^G$, $i=1,2$, and the parity of the intersection form of $S=(C_1\times C_2)/G$.
\smallskip

Recall now,  by Remark~\ref{Schur}, that $\operatorname{Hom}(\mathfrak{A}_5,\mathbb{C}^*)=0$ and  $H^{2}(\mathfrak{A}_5, \mathbb{C}^*)=\mathbb{Z}/2.$

\begin{lemma} $\operatorname{Pic}(G;C_1)=  6\mathbb{Z} \oplus \mathbb{Z}/5$ and  $\operatorname{Pic}(G;C_2)=  20\mathbb{Z} \oplus (\mathbb{Z}/3)^3.$
\end{lemma}
\begin{proof}  These follow from Proposition \ref{Pic(G;C)} since the signatures are $T_1=[2,5,5]$, $T_2=[3,3,3,3]$.
\end{proof}
Now the exact sequence \eqref{eq:exact} yields the following:
\begin{subequations}\label{eq:row1-sequences}
\begin{equation}\label{eq:row1-C1-concrete}
0
\longrightarrow 0
\longrightarrow 6\mathbb{Z} \oplus \mathbb{Z}/5
\longrightarrow \operatorname{Pic}(C_{1})^{G}
\longrightarrow \mathbb{Z}/2
\longrightarrow 0 .
\end{equation}
\begin{equation}\label{eq:row1-C2-concrete}
0
\longrightarrow 0
\longrightarrow 20\mathbb{Z} \oplus (\mathbb{Z}/3)^{3}
\longrightarrow \operatorname{Pic}(C_{2})^{G}
\longrightarrow \mathbb{Z}/2
\longrightarrow 0 .
\end{equation}
\end{subequations}

Since
\[
\operatorname{Pic}(C_1)^G[2] \cong 0
\quad\text{and}\quad
\operatorname{Pic}(C_2)^G[2] \cong \mathbb{Z}/2,
\]

from the exact sequences
\eqref{eq:row1-C1-concrete} and \eqref{eq:row1-C2-concrete} we infer that
\[
\operatorname{Pic}(C_1)^G \cong 3\mathbb{Z} \oplus \mathbb{Z}/5
\quad\text{and}\quad
\operatorname{Pic}(C_2)^G \cong 20\mathbb{Z} \oplus (\mathbb{Z}/3)^3 \oplus \mathbb{Z}/2 .
\]

From \eqref{eq:row1-C1-concrete}, we can choose a $G$-invariant line bundle $L_1\in \operatorname{Pic}(C_1)^G$ of degree $3=g(C_1)-1$ which is not $G$-linearized. Since the map $\delta_2$ for $C_2$ is surjective, there exists a $G$-invariant torsion line bundle $M_1 \in \operatorname{Pic}(C_2)^G$ such that $\delta_2(M_1)=\delta_1(L_1)$.  Then the external tensor product
$L_1\boxtimes M_1$ is a $G$-linearized line bundle on $C_1\times C_2$. Since $M_1$ is a torsion, this line bundle has bidegree $(3,0)$.
On the other hand, from \eqref{eq:row1-C2-concrete} we can choose a $G$-linearized line bundle on $C_2$ of degree $20=g(C_2)-1$, whose pullback on $C_1\times C_2$ is a $G$-linearized line bundle of bidegree $(0,20)$.
Therefore, by Strategy (1) the intersection form of $S$ is even.

\bigskip

\section{$G=\mathfrak{A}_5$, $g(C_1)=13$, $g(C_2)=6$, $T_1=[5,5,5]$, $T_2=[2,2,2,3]$; odd fake quadrics}

\medskip

\begin{proposition}  In the case $G=\mathfrak{A}_5$, $g(C_1)=13$, $g(C_2)=6$, $T_1=[5,5,5]$, $T_2=[2,2,2,3]$, we have
$$\operatorname{Pic}(C_1)^G[2] \cong 0 \quad\text{and}\quad \#\Theta(C_1)^G=1.$$
\end{proposition}
\begin{proof}
For $C_1$, a generating vector is $(g_1,g_2,g_3)=\bigl((12534),(12453),(12345)\bigr)$ with orders $5,5,5.$  For further details, we refer the reader to \cite[\S~6]{BCG08}. Let
\[
\Gamma_1=\langle x_1,x_2,x_3\mid x_1^5=x_2^5=x_3^5=1,\ x_1x_2x_3=1\rangle.
\]
and let $\varphi_1:\Gamma_1 \twoheadrightarrow \mathfrak{A}_5$, $x_i\mapsto g_i$.
 To compute $\#\Theta(C_1)^G$, by Proposition~\ref{dim}, we compute
\[
\dim H^1(C,\mathbb{Z}/2)^G=\dim H^1(\Gamma,\mathbb{Z}/2)-\dim H^1(G,\mathbb{Z}/2)+\dim\ker\varphi^*.
\]
By Table~\ref{dim table}, we know that $\dim H^1(\Gamma_1,\mathbb{Z}/2)=\dim H^1(\mathfrak{A}_5,\mathbb{Z}/2)=0$. Thus,
\[
H^1(C_1,\mathbb{Z}/2)^{\mathfrak{A}_5}\cong\ker(H^2(\mathfrak{A}_5,\mathbb{Z}/2)\overset{\varphi_1^*}{\longrightarrow}H^2(\Gamma_1,\mathbb{Z}/2)).
\]
We show that the generating vector does not lift relation-preservingly to $2.\mathfrak{A}_5\cong\mathrm{SL}_2(\mathbb{F}_5)$. Suppose, to the contrary, that its three order-$5$ elements lift to order-$5$ elements $\widetilde g_1,\widetilde g_2,\widetilde g_3\in\mathrm{SL}_2(\mathbb{F}_5)$ satisfying $\widetilde g_1\widetilde g_2\widetilde g_3=I$. An order-$5$ element of $\mathrm{SL}_2(\mathbb{F}_5)$ is nontrivial unipotent, because in characteristic $5$, $X^5-1=(X-1)^5$. Thus its trace is $2$. After conjugation, write
\[
\widetilde g_1=U:=
\begin{pmatrix}
1&1\\
0&1
\end{pmatrix}.
\]
Write another order-$5$ lift as
\[
B=\begin{pmatrix}
1+a&b\\
c&1-a
\end{pmatrix}.
\]
Since $\operatorname{tr}(B)=2$ and $\det B=1$, $a^2+bc=0$. Now
\[
UB=
\begin{pmatrix}
1+a+c&b+1-a\\
c&1-a
\end{pmatrix},
\]
so
\[
\operatorname{tr}(UB)=2+c.
\]
If $\widetilde g_1\widetilde g_2\widetilde g_3=I$, then $\widetilde g_1\widetilde g_2=\widetilde g_3^{-1}$. The inverse of an order-$5$ unipotent is again order $5$ and hence has trace $2$.  Therefore $\operatorname{tr}(UB)=2$, so $c=0$. Then $a^2+bc=a^2=0$, whence $a=0$. Thus
\[
B=
\begin{pmatrix}
1&b\\
0&1
\end{pmatrix},
\]
so $U$ and $B$ lie in the same upper-unipotent subgroup
\[
\{
\begin{pmatrix}
1&t\\
0&1
\end{pmatrix}
:t\in\mathbb{F}_5
\}
\cong\mathbb{Z}/5.
\]

Their images in $\mathfrak{A}_5$ therefore lie in one cyclic $5$-subgroup and cannot generate $\mathfrak{A}_5$.  This contradicts the fact that the spherical generating vector generates $\mathfrak{A}_5$. Hence $\varphi_1^*\neq0.$ Since $H^2(\mathfrak{A}_5,\mathbb{Z}/2)$ is one-dimensional, $\varphi_1^*$ is injective. Therefore $H^1(C_1,\mathbb{Z}/2)^{\mathfrak{A}_5}=0$.

 Proposition~\ref{admit theta} gives an $\mathfrak{A}_5$-invariant theta characteristic, hence $\#\Theta(C_1)^G=1.$
\end{proof}

\begin{proposition}  In the case $G=\mathfrak{A}_5$, $g(C_1)=13$, $g(C_2)=6$, $T_1=[5,5,5]$, $T_2=[2,2,2,3]$, we have
$$\operatorname{Pic}(C_2)^G[2] \cong (\mathbb{Z}/2)^2 \quad\text{and}\quad \#\Theta(C_2)^G=4.$$
\end{proposition}

\begin{proof}
For $C_2$, a generating vector is $(g_1,g_2,g_3,g_4)=\bigl((12)(34),(24)(35),(14)(35),(234)\bigr)$ with orders $2,2,2,3$.  Let
\[
\Gamma_2=\langle y_1,y_2,y_3,y_4\mid y_1^2=y_2^2=y_3^2=y_4^3=1,\ y_1y_2y_3y_4=1\rangle
\]
and let $\varphi_2:\Gamma_2\twoheadrightarrow \mathfrak{A}_5$, $y_i\mapsto g_i$.

 To compute $\#\Theta(C_2)^G$, by Proposition~\ref{dim}, we compute
\[
\dim H^1(C,\mathbb{Z}/2)^G=\dim H^1(\Gamma,\mathbb{Z}/2)-\dim H^1(G,\mathbb{Z}/2)+\dim\ker\varphi^*.
\]
By Table~\ref{dim table}, we know that $\dim H^1(\Gamma_2,\mathbb{Z}/2)=2$ and $\dim H^1(\mathfrak{A}_5,\mathbb{Z}/2)=0$.

\medskip

The images of $y_1,y_2,y_3$ are nontrivial involutions in $\mathfrak{A}_5$. By Lemma~\ref{A_5 involution}, every lift of such an involution to $2.\mathfrak{A}_5$ has order $4$. Therefore no lift can satisfy $y_1^2=y_2^2=y_3^2=1$. Hence the map $\phi_2$ cannot lift relation-preservingly to $2.\mathfrak{A}_5$, and $\varphi_2^*\neq0$. Since $H^2(\mathfrak{A}_5,\mathbb{Z}/2)$ is one-dimensional, $\ker\varphi_2^*=0$.
The five-term exact sequence therefore gives
\[
H^1(C_2,\mathbb{Z}/2)^{\mathfrak{A}_5}\cong H^1(\Gamma_2,\mathbb{Z}/2)\cong(\mathbb{Z}/2)^2.
\]
 The  proposition below gives an invariant theta characteristic, thus $\#\Theta(C_2)^{\mathfrak{A}_5}=4.$
\end{proof}

\begin{proposition}\label{existence inv theta char}
There exists an $\mathfrak{A}_5$-invariant theta characteristic on $C_2$.
\end{proposition}

\begin{proof}
We consider the case $G=\mathfrak{A}_5$, $T_2=[2,2,2,3]$ with $g(C_2)=6$. Let $\pi:C_2\longrightarrow C_2/G\cong \mathbb{P}^1$ be the quotient map. Since the signature of the $G$-action is $[2,2,2,3]$, there are four branch points $p_1,p_2,p_3,p_4\in \mathbb{P}^1$, whose stabilizer orders are $2,2,2,3$, respectively. For each $i$, let $E_i:=\pi^{-1}(p_i)_{\mathrm{red}}$ denote the reduced fiber over $p_i$, and let $F$ denote the pullback of an ordinary point of $\mathbb{P}^1$. Thus $\pi^*(p_i)=2E_i$, $i=1,2,3$, whereas $\pi^*(p_4)=3E_4$. Since all points of $\mathbb{P}^1$ are linearly equivalent, we obtain
\[
2E_1\equiv 2E_2\equiv 2E_3\equiv 3E_4\equiv F.
\]

 Setting $A:=\operatorname{Pic}(\mathfrak{A}_5;C_2)$, $B:=\operatorname{Pic}(C_2)^{\mathfrak{A}_5}$, the exact sequence becomes
\[
0
\longrightarrow A
\longrightarrow B
\overset{\delta}{\longrightarrow} \mathbb{Z}/2
\longrightarrow 0.
\]
For the signature $[2,2,2,3]$, the $\operatorname{Pic}(\mathfrak{A}_5;C_2)$ has the presentation
\[
A=\langle F,E_1,E_2,E_3,E_4 \mid 2E_1=2E_2=2E_3=3E_4=F \rangle.
\]
We apply multiplication by $2$ to $0 \longrightarrow A \longrightarrow B \overset{\delta}{\longrightarrow} \mathbb{Z}/2 \longrightarrow 0$. This gives a commutative diagram
\[
\begin{tikzcd}
0 \arrow[r]&A \arrow[r] \arrow[d,"2"]&B \arrow[r,"\delta"] \arrow[d,"2"]&\mathbb{Z}/2 \arrow[r] \arrow[d,"0"]&0\\
0 \arrow[r]&A \arrow[r]&B \arrow[r,"\delta"]&\mathbb{Z}/2 \arrow[r]&0.
\end{tikzcd}
\]

By the snake lemma
\[
\begin{tikzpicture}[>=stealth]
  \matrix[
    matrix of math nodes,
    column sep={70pt,between origins},
    row sep={60pt,between origins},
    nodes={asymmetrical rectangle}
  ] (s)
  {
    & |[name=ka]|  A[2]
    & |[name=kb]|  B[2]
    & |[name=kc]| \mathbb{Z}/2 \\
    & |[name=A]|  A
    & |[name=B]|  B
    & |[name=C]| \mathbb{Z}/2
    & |[name=01]| 0 \\
    |[name=02]| 0
    & |[name=A']|  A
    & |[name=B']|  B
    & |[name=C']| \mathbb{Z}/2 \\
    & |[name=ca]|  A/2A
    & |[name=cb]|  B/2B
    & |[name=cc]| \mathbb{Z}/2 \\
  };

  \draw[->]
    (ka) edge (A)
    (kb) edge (B)
    (kc) edge (C)
    (A) edge (B)
    (B) edge (C)
    (C) edge (01)
    (A) edge node[auto] {$\times 2$} (A')
    (B) edge node[auto] {$\times 2$} (B')
    (C) edge node[auto] {$0$} (C')
    (02) edge (A')
    (A') edge (B')
    (B') edge (C')
    (A') edge (ca)
    (B') edge (cb)
    (C') edge (cc);

  \draw[->,gray]
    (ka) edge (kb)
    (kb) edge (kc)
    (ca) edge (cb)
    (cb) edge (cc);

  \draw[->,gray,rounded corners]
    (kc) -|
    node[auto,text=black,pos=.7] {$\partial$}
    ($(01.east)+(.5,0)$)
    |- ($(B)!.35!(B')$)
    -| ($(02.west)+(-.5,0)$)
    |- (ca);
\end{tikzpicture}
\]
we obtain the exact sequence
\[
0\longrightarrow \operatorname{Pic}(\mathfrak{A}_5;C_2)[2]
\longrightarrow \operatorname{Pic}(C_2)^{\mathfrak{A}_5}[2]
\longrightarrow \mathbb{Z}/2
\overset{\partial}{\longrightarrow}
\operatorname{Pic}(\mathfrak{A}_5;C_2)/2\operatorname{Pic}(\mathfrak{A}_5;C_2)= A/2A.
\]

The snake lemma gives a connecting homomorphism $\partial:\mathbb{Z}/2\longrightarrow A/2A$. We know explicitly how $\partial$ is defined. Indeed, let $1\in\mathbb{Z}/2$ be the nonzero element. Choose $b\in B$ such that $\delta(b)=1$. Since $\delta(2b)=2\delta(b)=0$, we have $2b\in A$. Thus, $\partial(1)=[2b]\in A/2A$. We now reduce the presentation $A=\langle F,E_1,E_2,E_3,E_4 \mid 2E_1=2E_2=2E_3=3E_4=F \rangle$ modulo $2A$, that is,
\[
A/2A=\langle \overline E_1, \overline E_2, \overline E_3 \rangle.
\]
For $i=1,2,3$, choose a point $x_i\in E_i$. Its stabilizer group is $I_i:=\operatorname{Stab}_G(x_i)=\langle g_i\rangle\cong\mathbb{Z}/2$. If $M$ is an $\mathfrak{A}_5$-linearized line bundle, then $I_i$ acts on the one-dimensional fiber $M_{x_i}$, and hence determines a character
\[
\chi_i(M): I_i\longrightarrow\mathbb{C}^*.
\]
Since $I_i\cong\mathbb{Z}/2$, this character is either trivial or nontrivial. We identify $\operatorname{Hom}(I_i,\mathbb{C}^*)\cong\mathbb{Z}/2$ by sending the trivial character to $0$ and the nontrivial character to $1$. Thus we obtain
\[
 \chi:=(\chi_1,\chi_2,\chi_3): A\longrightarrow (\mathbb{Z}/2)^3.
\]
If $L=2M$, then the stabilizer character of $L$ is the square of that of $M$. Thus, every element of $2A$ has trivial stabilizer character. Hence the map factors through
\[
\overline\chi: A/2A \longrightarrow (\mathbb{Z}/2)^3.
\]
We now compute this map on the three classes $\overline E_1,\overline E_2,\overline E_3$. Consider $\mathcal{O}_{C_2}(E_i)$, $i=1,2,3$, with its natural $\mathfrak{A}_5$-linearization. At a point $x_i\in E_i$, choose a local parameter $z$ such that $E_i=(z=0)$. Since $g_i$ generates the order-two stabilizer group  and $g_i(x_i)=x_i$, $g_i^*(z)=-z$.

A local generator of $\mathcal{O}_{C_2}(E_i)$ is $1/z$, and therefore
\[
g_i^*\left(\frac1z\right)=-\frac1z.
\]
Thus the stabilizer character at $E_i$ is nontrivial. For $i\neq j$, the divisor $E_i$ does not contain $x_j$, so stabilizer character of $\mathcal{O}_{C_2}(E_i)$ at $E_j$ is trivial. Consequently,
\[
\overline E_1\longmapsto(1,0,0), \quad \overline E_2\longmapsto(0,1,0), \quad \overline E_3\longmapsto(0,0,1).
\]
Then, we obtain an isomorphism $\overline\chi: A/2A \overset{\sim}{\longrightarrow} (\mathbb{Z}/2)^3$. Thus a class in $A/2A$ is uniquely determined by its three order-two stabilizer characters.

Let $b\in B$ be a lift of the nonzero element of $H^2(\mathfrak{A}_5,\mathbb{C}^*)\simeq\mathbb{Z}/2$. Let $L$ be an $\mathfrak{A}_5$-invariant line bundle representing $b$.
The obstruction to linearizing $L$ is the unique nonzero class of $\mathfrak{A}_5$. This class is represented by the central extension
\[
1
\longrightarrow\{\pm1\}
\longrightarrow 2.\mathfrak{A}_5
\longrightarrow \mathfrak{A}_5
\longrightarrow 1.
\]
Choose a section $s:\mathfrak{A}_5\longrightarrow 2.\mathfrak{A}_5$. Then $s(g)s(h)=c(g,h)s(gh)$ for a cocycle $c(g,h)\in\{\pm1\}$. After rescaling the chosen isomorphisms $\phi_g:g^*L\longrightarrow L$, we may assume
\[
\phi_g\circ g^*(\phi_h)=c(g,h)\phi_{hg}.
\]
Indeed, if $L$ represents $b$, then $[c_L]$ is the unique nonzero class in $H^2(\mathfrak{A}_5,\mathbb{C}^*)$. On the other hand, $c$ is a $2$-cocycle representing the same nonzero cohomology class. Therefore $[c_L]=[c]$ in $H^2(\mathfrak{A}_5,\mathbb{C}^*)$. Hence the two cocycles differ by a coboundary. More precisely, there
exists a $1$-cochain $a:\mathfrak{A}_5\longrightarrow\mathbb{C}^*$, $g\longmapsto a_g$, such that $c_L(g,h)=a_g^{-1}a_h^{-1}a_{hg}\,c(g,h)$. Equivalently, $c(g,h)=a_g a_h a_{hg}^{-1}c_L(g,h)$.
Now replace the chosen isomorphisms $\phi_g$ by $\phi'_g:=a_g\phi_g$. Then
\[
\begin{aligned}
\phi'_g\circ g^*(\phi'_h)&=a_g\phi_g\circ g^*(a_h\phi_h)=a_ga_h\bigl(\phi_g\circ g^*(\phi_h)\bigr)=a_ga_hc_L(g,h)\phi_{hg}=a_ga_ha_{hg}^{-1}c_L(g,h)\phi'_{hg}.
\end{aligned}
\]
Therefore, if $c(g,h)=a_ga_ha_{hg}^{-1}c_L(g,h)$, then $\phi'_g\circ g^*(\phi'_h)=c(g,h)\phi'_{hg}$.

\medskip

We now compute $\partial(1)=[2b]$. The line bundle $L^{\otimes2}$ represents $2b$. This means $L^{\otimes2}$ admits an $\mathfrak{A}_5$-linearization. This linearization can be seen explicitly. Namely, consider
\[
\phi_g^{\otimes2}:g^*(L^{\otimes2})
\simeq
(g^*L)^{\otimes2}
\longrightarrow
L^{\otimes2}.
\]
We obtain
\[
\begin{aligned}
\phi_g^{\otimes2}\circ g^*(\phi_h^{\otimes2})&=\bigl(\phi_g\circ g^*(\phi_h)\bigr)^{\otimes2}=c(g,h)^2\phi_{hg}^{\otimes2}=\phi_{hg}^{\otimes2},
\end{aligned}
\]
because $c(g,h)\in\{\pm1\}$. Thus $\{\phi_g^{\otimes2}\}_{g\in \mathfrak{A}_5}$ defines a $\mathfrak{A}_5$-linearization of $L^{\otimes2}$.

Fix $i\in\{1,2,3\}$. Since $g_i$ fixes $x_i$, its action on the fiber $L_{x_i}$ is multiplication by some scalar $\lambda_i\in\mathbb{C}^*$. Applying $\phi_g\circ g^*(\phi_h)=c(g,h)\phi_{hg}$ with $g=h=g_i$ gives $\lambda_i^2=c(g_i,g_i)$.

 Note that $c(g_i,g_i)=-1$, $i=1,2,3$. Indeed, as observed above, the cocycle $c(g,h)$ for the section $s:\mathfrak{A}_5\longrightarrow2.\mathfrak{A}_5$ and the cocycle $c_L(g,h)$ for the line bundle $L$ represent the same cohomology class in $H^2(\mathfrak{A}_5,\mathbb{C}^*)$. For an element $g_i \in \mathfrak{A}_5$ of order $2$, its lift $s(g_i)\in 2.\mathfrak{A}_5$ has order $4$. Hence $-1=s(g_i)^2=c(g_i,g_i)s(g_i^2)=c(g_i,g_i)$. This implies $\lambda_i^2=-1$. Therefore $g_i$ acts on $(L^{\otimes2})_{x_i}$ by $\lambda_i^2=-1$. Thus $L^{\otimes2}$ has the nontrivial stabilizer character. Consequently, $\overline\chi([2b])=(1,1,1)$. Under the isomorphism $A/2A\simeq(\mathbb{Z}/2)^3$, the vector $(1,1,1)$ corresponds precisely to $\overline E_1+\overline E_2+\overline E_3$. Hence
\[
\partial(1)=[2b]=\overline E_1+\overline E_2+\overline E_3.
\]
We now compute the canonical divisor. By the Riemann-Hurwitz formula,
\[
K_{C_2} \equiv \pi^*K_{\mathbb{P}^1}+R,
\]
where $R$ is the ramification divisor. Hence
\[
K_{C_2} \equiv -2F+E_1+E_2+E_3+2E_4.
\]
Therefore $[K_{C_2}]=\overline E_1+\overline E_2+\overline E_3=\partial(1) $in $A/2A$. Since $[K_{C_2}]=[2b]$ in $A/2A$, we have $K_{C_2}-2b \in 2A$. Therefore there exists $a\in A$ such that $K_{C_2}-2b=2a$. Equivalently, $K_{C_2}=2(a+b)$. Set $\kappa:=a+b$. Since $a\in A\subset B$ and $b\in B$, $\kappa\in B=\operatorname{Pic}(C_2)^{\mathfrak{A}_5}$.  Thus
$\kappa\in\Theta(C_2)^{\mathfrak{A}_5}.$
\end{proof}

 We now determine $\operatorname{Pic}(C_i)^G$  and the parity of the intersection form of $S$.

\medskip

In this case, by Remark~\ref{Schur} and Proposition \ref{Pic(G;C)} we have
$$\operatorname{Hom}(\mathfrak{A}_5,\mathbb{C}^*)=0,\,\,\, H^{2}(\mathfrak{A}_5, \mathbb{C}^*)=\mathbb{Z}/2$$
$$\operatorname{Pic}(G;C_1)=  12\mathbb{Z} \oplus (\mathbb{Z}/5)^2,\,\,\, \operatorname{Pic}(G;C_2)=  10\mathbb{Z} \oplus (\mathbb{Z}/2)^2.$$
The exact sequence \eqref{eq:exact} yields the following:
\begin{subequations}\label{eq:A5-row2-sequences}
\begin{equation}\label{eq:A5-row2-C1-concrete}
0
\longrightarrow 0
\longrightarrow 12\mathbb{Z} \oplus (\mathbb{Z}/5)^2
\longrightarrow \operatorname{Pic}(C_{1})^{G}
\longrightarrow \mathbb{Z}/2
\longrightarrow 0 .
\end{equation}
\begin{equation}\label{eq:A5-row2-C2-concrete}
0
\longrightarrow 0
\longrightarrow 10\mathbb{Z} \oplus (\mathbb{Z}/2)^{2}
\longrightarrow \operatorname{Pic}(C_{2})^{G}
\longrightarrow \mathbb{Z}/2
\longrightarrow 0 .
\end{equation}
\end{subequations}
  Since $\operatorname{Pic}(C_1)^G[2]=0$ and $\operatorname{Pic}(C_2)^G[2]= (\mathbb{Z}/2)^2$, we infer that
$$\operatorname{Pic}(C_1)^G=6\mathbb{Z} \oplus (\mathbb{Z}/5)^2,$$
$$\operatorname{Pic}(C_2)^G = 5\mathbb{Z}\oplus (\mathbb{Z}/2)^2 \quad\text{ or}\quad 10\mathbb{Z}\oplus \mathbb{Z}/2\oplus \mathbb{Z}/4.$$
The second possibility for $\operatorname{Pic}(C_2)^G$ can be ruled out, because there is a $G$-invariant theta characteristic on $C_2$, which has degree $5=g(C_2)-1$.
From \eqref{eq:A5-row2-C1-concrete}, we can choose a $G$-invariant line bundle $L_1\in \operatorname{Pic}(C_1)^G$ of degree $6$ which is not  $G$-linearized. Likewise, \eqref{eq:A5-row2-C2-concrete} yields a degree $5=g(C_2)-1$ line bundle $L_2\in \operatorname{Pic}(C_2)^G$ whose obstruction class in the Schur multiplier coincides with that of $L_1$. Consequently, $L_1\boxtimes L_2$ defines a $G$-linearized line bundle on $C_1\times C_2$ with bidegree $(6,5)$.

By Strategy (1) its  descent line bundle has self-intersection 1, an odd number.

\bigskip

\section{$G=\mathfrak{A}_5$, $g(C_1)=5$, $g(C_2)=16$, $T_1=[3,3,5]$, $T_2=[2,2,2,2,2]$; odd fake quadrics}

\medskip

\begin{proposition}  In the case $G=\mathfrak{A}_5$, $g(C_1)=5$, $g(C_2)=16$, $T_1=[3,3,5]$, $T_2=[2,2,2,2,2]$, we have
$$\operatorname{Pic}(C_1)^G[2] \cong 0 \quad\text{and}\quad \#\Theta(C_1)^G=1.$$
\end{proposition}
\begin{proof}
For $C_1$, a generating vector is $(g_1,g_2,g_3)=\bigl((123),(345),(54321)\bigr)$ with orders $3,3,5$.  Let
\[
\Gamma_1=\langle x_1,x_2,x_3\mid x_1^3=x_2^3=x_3^5=1,\ x_1x_2x_3=1\rangle.
\]
and let $\varphi_1:\Gamma_1\twoheadrightarrow \mathfrak{A}_5,$ $x_i\mapsto g_i.$

 To compute $\#\Theta(C_1)^G$, by Proposition~\ref{dim} we compute
\[
\dim H^1(C,\mathbb{Z}/2)^G=\dim H^1(\Gamma,\mathbb{Z}/2)-\dim H^1(G,\mathbb{Z}/2)+\dim\ker\varphi^*.
\]
By Table~\ref{dim table}, we know that $\dim H^1(\Gamma_1,\mathbb{Z}/2)=\dim H^1(\mathfrak{A}_5,\mathbb{Z}/2)=0$. Thus,
\[
H^1(C_1,\mathbb{Z}/2)^{\mathfrak{A}_5}\cong\ker(H^2(\mathfrak{A}_5,\mathbb{Z}/2)\overset{\varphi_1^*}{\longrightarrow}H^2(\Gamma_1,\mathbb{Z}/2)).
\]

Recall the definition $\alpha_{ij}=\frac{e_i-e_j}{\sqrt2}$. Each $\alpha_{ij}$ is the corresponding reflection is the transposition $(ij)$. For convenience, put
\[
a:=\alpha_{12},
\quad
b:=\alpha_{23},
\quad
c:=\alpha_{34},
\quad
d:=\alpha_{45}.
\]
Since the first two generators are $g_1=(123)=(12)(23)$ and $g_2=(345)=(34)(45)$, we may choose the lifts $\widetilde g_1=ab$ and $\widetilde g_2=cd$. Moreover, these are order-preserving lifts.

\medskip

Since the third generator is $g_3=(54321)=(45)(34)(23)(12)$, we may choose the lift $t:=dcba$. We claim that
\[
t^5=-1.
\]
We verify this directly in the Clifford algebra. Since $a^2=b^2=c^2=d^2=1$, we have $ab+ba=bc+cb=cd+dc=-1$. Also, $ac=-ca$, $ad=-da$, $bd=-db$.

\medskip

We first compute $t^2$. Since $t=dcba$, we have  $t^2=(dcba)(dcba)$. Using $ad=-da$ and $bd=-db$, we obtain $bad=dba$. Hence  $t^2=dc(bad)cba=dc(dba)cba$. Since $cd+dc=-1$, we have $cd=-1-dc$, and therefore
\[
dcd=d(-1-dc)=-d-c.
\]
It follows that
\[
 t^2=(dcd)bacba= -(c+d)bacba=-cbacba-dbacba.
\]

We now simplify the two terms separately.

\medskip

\noindent
First, since $cbc=(-1-bc)c=-c-b$ and $aba=a(-1-ab)=-a-b$, we have
$$
cbacba=-cbcaba=-(c+b)(a+b)=-ca-cb-ba-b^2=ac+(1+bc)+(1+ab)-1=1+ab+ac+bc.
$$
For the second term, we use $dba=bad$ to get $dbacba=badcba$. Since $dc=-1-cd$, $baba=ab$ and $bacdba=abcd-bd$, we have
$$
dbacba=badcba=ba(-1-cd)ba=-baba-bacdba=-ab-abcd+bd.
$$
Combining the two identities, we obtain
$$
t^2=-(1+ab+ac+bc)-(bd-ab-abcd)=-1-ac-bc-bd+abcd.
$$
Next, to compute $t^4$ we set $S:=ac+bc+bd$, $Q:=abcd$. Then $t^2=-1-S+Q$, and therefore
$$
t^4=(-1-S+Q)^2=1+S^2+Q^2+2S-2Q-(SQ+QS).
$$
For convenience, put $R:=ab+ac+ad+bc+bd+cd$. Using the Clifford relations, one obtains
\[
S^2=-2-R-2Q, \quad Q^2=1+R+Q, \quad SQ+QS=2S-2Q.
\]
Indeed, for the first identity,
\[
S^2=(ac)^2+(bc)^2+(bd)^2+(ac)(bc)+(bc)(ac)+(ac)(bd)+(bd)(ac)+(bc)(bd)+(bd)(bc).
\]
Here, $(ac)^2=acac=a(ca)c=-a^2c^2=-1$ and similarly $(bd)^2=-1$. Also,
$$
(bc)^2=bcbc=b(cb)c=b(-1-bc)c=-bc-b^2c^2=-1-bc,
$$
$$
(ac)(bc)=acbc=a(cb)c=a(-1-bc)c=-ac-ab,
$$
$$
(bc)(ac)=bcac=b(ca)c=-bac^2=-ba=1+ab,
$$

$$
(ac)(bd)=acbd=a(cb)d=a(-1-bc)d=-ad-abcd,
$$
$$
(bd)(ac)=bdac=-badc=-ba(-1-cd)=ba+bacd=-1-ab-cd-abcd,
$$

$$
(bc)(bd)=bcbd=b(cb)d=b(-1-bc)d=-bd-cd,
$$
$$
(bd)(bc)=bdbc=-bbdc=-dc=1+cd.
$$

Combining all these identities, we obtain
$$
S^2=-2-ab-ac-ad-bc-bd-cd-2abcd=-2-R-2Q.
$$

Similarly, direct reduction of $Q^2=abcdabcd$ gives $Q^2=1+R+Q$, $SQ+QS=2S-2Q$. Substituting these identities into the expression for $t^4$, we obtain
$$
t^4=1+(-2-R-2Q)+(1+R+Q)+2S-2Q-(2S-2Q)=-Q=-abcd.
$$
Therefore $t^5=(-abcd)(dcba)=-1$.  It follows that $\operatorname{ord}(t)=10$ and $\operatorname{ord}(-t)=5$. The other lift of $g_3$ is $-t$, the order-preserving lift of $g_3$.
If $\widetilde g_3=-t$, we have
$$
\widetilde g_1\widetilde g_2\widetilde g_3=(ab)(cd)(-t)=(ab)(cd)(-dcba)=-1\neq 1.
$$
Moreover, the order-preserving lift of $g_i$ is uniquely determined. The other lift of an order $3$ element has order $6$, and the other lift of the order $5$ element has order $10$. Hence there is no lifting map $\widetilde\varphi_1:\Gamma_1\longrightarrow\widetilde{\mathfrak{A}_5}$ such that $p\circ\widetilde\varphi_1=\varphi_1$.

\medskip

Since $H^2(\mathfrak{A}_5,\mathbb{Z}/2)\cong\mathbb{Z}/2$ is one-dimensional, every nonzero linear map out of it is injective. Therefore $\ker\varphi_1^*=0$. It follows that $H^1(C_1,\mathbb{Z}/2)^{\mathfrak{A}_5}=0$.  Then Proposition~\ref{admit theta} gives an invariant theta characteristic, hence $\#\Theta(C_1)^G=1$.
\end{proof}

\begin{proposition}  In the case $G=\mathfrak{A}_5$, $g(C_1)=5$, $g(C_2)=16$, $T_1=[3,3,5]$, $T_2=[2,2,2,2,2]$, we have
$$\operatorname{Pic}(C_2)^G[2] \cong (\mathbb{Z}/2)^4 \quad\text{and}\quad \#\Theta(C_2)^G=16.$$
\end{proposition}
\begin{proof}
For $C_2$, a generating vector is $(g_1,g_2,g_3,g_4,g_5)=\bigl((12)(34),(13)(24),(14)(23),(14)(25),(14)(25)\bigr)$ with orders $2,2,2,2,2.$  Let
\[
\Gamma_2=\langle y_1,y_2,y_3,y_4,y_5\mid y_1^2=y_2^2=y_3^2=y_4^2=y_5^2=1,\;y_1y_2y_3y_4y_5=1\rangle,
\]
and let $\varphi_2:\Gamma_2\twoheadrightarrow \mathfrak{A}_5$, $y_i\mapsto g_i$.

 To compute $\#\Theta(C_2)^G$, by Proposition~\ref{dim}, we compute
\[
\dim H^1(C,\mathbb{Z}/2)^G=\dim H^1(\Gamma,\mathbb{Z}/2)-\dim H^1(G,\mathbb{Z}/2)+\dim\ker\varphi^*.
\]
By Table~\ref{dim table}, we know that
\[
\dim H^1(\Gamma_2,\mathbb{Z}/2)=4,\quad \dim H^1(\mathfrak{A}_5 ,\mathbb{Z}/2)=0, \quad \dim H^2(\mathfrak{A}_5,\mathbb{Z}/2)=1.
\]

The images of $y_1,y_2,y_3,y_4,y_5$ are nontrivial involutions in $\mathfrak{A}_5$. By Lemma~\ref{A_5 involution}, every lift of such an involution to $2.\mathfrak{A}_5$ has order $4$. Therefore no lift can satisfy $y_1^2=y_2^2=y_3^2=y_4^2=y_5^2=1$. Hence the map $\phi_2$ cannot lift relation-preservingly to $2.\mathfrak{A}_5$, and $\varphi_2^*\neq0$. Since $H^2(\mathfrak{A}_5,\mathbb{Z}/2)$ is one-dimensional, $\ker\varphi_2^*=0$. It follows that
\[
\dim H^1(C_2,\mathbb{Z}/2)^{\mathfrak{A}_5}=4.
\]
 The proposition below gives an invariant theta characteristic, thus
$
\#\Theta(C_2)^{\mathfrak{A}_5}=16.$
\end{proof}

\begin{proposition}
There exists an $\mathfrak{A}_5$-invariant theta characteristic on $C_2$.
\end{proposition}

\begin{proof}
We use the same argument as in the proof of Proposition~\ref{existence inv theta char}.
\end{proof}

\medskip
 We now determine $\operatorname{Pic}(C_i)^G$ and the parity of the intersection form of $S$.

In this case, by Remark~\ref{Schur} and Proposition \ref{Pic(G;C)} we have
$$\operatorname{Hom}(\mathfrak{A}_5,\mathbb{C}^*)=0,\,\,\, H^{2}(\mathfrak{A}_5, \mathbb{C}^*)=\mathbb{Z}/2$$
$$\operatorname{Pic}(G;C_1)=  4\mathbb{Z} \oplus \mathbb{Z}/3,\,\,\, \operatorname{Pic}(G;C_2)=  30\mathbb{Z} \oplus (\mathbb{Z}/2)^4.$$
The exact sequence \eqref{eq:exact} yields
\begin{subequations}\label{eq:A5-row3-sequences}
\begin{equation}\label{eq:A5-row3-C1-concrete}
0
\longrightarrow 0
\longrightarrow 4\mathbb{Z} \oplus \mathbb{Z}/3
\longrightarrow \operatorname{Pic}(C_{1})^{G}
\longrightarrow \mathbb{Z}/2
\longrightarrow 0 .
\end{equation}
\begin{equation}\label{eq:A5-row3-C2-concrete}
0
\longrightarrow 0
\longrightarrow 30\mathbb{Z} \oplus (\mathbb{Z}/2)^{4}
\longrightarrow \operatorname{Pic}(C_{2})^{G}
\longrightarrow \mathbb{Z}/2
\longrightarrow 0 .
\end{equation}
\end{subequations}
Since $\operatorname{Pic}(C_1)^G[2]=0$ and $\operatorname{Pic}(C_2)^G[2]= (\mathbb{Z}/2)^4$,  we infer that
$$\operatorname{Pic}(C_1)^G=2\mathbb{Z} \oplus \mathbb{Z}/3,$$
$$
\operatorname{Pic}(C_2)^G = 15\mathbb{Z}\oplus (\mathbb{Z}/2)^4.$$

By \eqref{eq:A5-row3-C1-concrete} we can choose a $G$-invariant line bundle $L_1\in \operatorname{Pic}(C_1)^G$ of degree $2$ which is not $G$-linearized. Likewise, \eqref{eq:A5-row3-C2-concrete} yields a degree $15=g(C_2)-1$ line bundle $L_2\in \operatorname{Pic}(C_2)^G$ whose obstruction class in the Schur multiplier coincides with that of $L_1$. Consequently, $L_1\boxtimes L_2$ defines a $G$-linearized line bundle on $C_1\times C_2$ with bidegree $(2,15)$.   By Strategy (1) its  descent  line bundle has self-intersection 1.

\bigskip

\section{$G=\mathfrak{S}_4\times \mathbb{Z}/2$, $g(C_1)=3$, $g(C_2)=25$, $T_1=[2,4,6]$, $T_2=[2,2,2,2,2,2]$; even fake quadrics}

\medskip

\begin{proposition}  In the case $G=\mathfrak{S}_4\times \mathbb{Z}/2$, $g(C_1)=3$, $g(C_2)=25$, $T_1=[2,4,6]$, $T_2=[2,2,2,2,2,2]$, we have
$$\operatorname{Pic}(C_1)^G[2] \cong \mathbb{Z}/2 \quad\text{and}\quad  \#\Theta(C_1)^G=2.$$
\end{proposition}
\begin{proof}
For $C_1$, a generating vector is $(g_1,g_2,g_3)=(((12),0),
((1234),1),((432),1))$ with orders $2,4,6$.  Let
\[
\Gamma_1=\langle x_1,x_2,x_3\mid x_1^2=x_2^4=x_3^6=1,\ x_1x_2x_3=1\rangle
\]
and let $\varphi_1:\Gamma_1\twoheadrightarrow \mathfrak{S}_4\times \mathbb{Z}/2$, $x_i\mapsto g_i.$

 To compute $\#\Theta(C_1)^G$, by Proposition~\ref{dim} and ~\ref{kervarphi} we compute
\[
\dim H^1(C,\mathbb{Z}/2)^G=\dim H^1(\Gamma,\mathbb{Z}/2)-\dim H^1(G,\mathbb{Z}/2)+\dim H^2(G,\mathbb{Z}/2)-\rank R.
\]
By Table~\ref{dim table}, we know that
\[
\dim H^1(\Gamma_1,\mathbb{Z}/2)=2,\quad \dim H^1(\mathfrak{S}_4\times \mathbb{Z}/2 ,\mathbb{Z}/2)=2, \quad \dim H^2(\mathfrak{S}_4\times \mathbb{Z}/2,\mathbb{Z}/2)=4.
\]
Thus, it is enough to compute the rank of the matrix $R$. We choose the following basis. Define $a,b:G\longrightarrow\mathbb{Z}/2$ by
\[
a(\sigma,\epsilon)=
\begin{cases}
0,&\sigma\in \mathfrak{A}_4,\\
1,&\sigma\notin \mathfrak{A}_4,
\end{cases}
\]
and $b(\sigma,\epsilon)=\epsilon$. Thus $a$ is the mod-$2$ sign character on the $\mathfrak{S}_4$-factor, whereas $b$ is the projection onto the second factor. These two characters are linearly independent. Indeed, suppose that $\lambda a+\mu b=0$, $\lambda,\mu\in\mathbb{Z}/2$. Evaluating at $((12),0)$, we obtain
\[
0=\lambda a((12),0)+\mu b((12),0)=\lambda,
\]
and hence $\lambda=0$. Similarly, evaluating at $(id,1)$, we obtain
\[
0=\lambda a(id,1)+\mu b(id,1)=\mu.
\]
Therefore $\mu=0$. Thus $a$ and $b$ are linearly independent. Since $H^1(G,\mathbb{Z}/2)$ has dimension $2$, it follows that
\[
H^1(G,\mathbb{Z}/2)=\langle a,b\rangle_{\mathbb{Z}/2}.
\]
We compute the connecting homomorphism $\beta:H^1(G;\mathbb{Z}/2)\longrightarrow H^2(G;\mathbb{Z}/2)$, i.e. $\beta(a)$, $\beta(b).$ We now return to $G=\mathfrak{S}_4\times\mathbb Z/2$. For the first curve, let
\[
\Gamma_1=\langle x_1,x_2,x_3\mid x_1^2=x_2^4=x_3^6=1,\;x_1x_2x_3=1\rangle
\]
and let $\varphi_1:\Gamma_1\twoheadrightarrow G$ be determined by
\[
\begin{aligned}
g_1:=\varphi_1(x_1)&=((12),0),\\
g_2:=\varphi_1(x_2)&=((1234),1),\\
g_3:=\varphi_1(x_3)&=((432),1).
\end{aligned}
\]
Note that $\operatorname{ord}(g_1)=2$, $\operatorname{ord}(g_2)=4$, $\operatorname{ord}(g_3)=6$. Thus, $\langle g_1\rangle\cong \mathbb{Z}/2$, $\langle g_2\rangle\cong \mathbb{Z}/4$, $\langle g_3\rangle\cong \mathbb{Z}/6$. Since $H^2(\langle g_i\rangle;\mathbb{Z}/2)\cong\mathbb{Z}/2$ for all $i$, the restriction map is
\[
R_1:H^2(G;\mathbb{Z}/2)\longrightarrow\bigoplus_{i=1}^3H^2(\langle g_i\rangle;\mathbb{Z}/2)\cong(\mathbb{Z}/2)^3.
\]
For $g_1=((12),0)$, $g_2=((1234),1)$, $g_3=((432),1)$,
\[
(a(g_1),a(g_2),a(g_3))=(1,1,0).
\]
By Lemma~\ref{restrict} $R_1(\beta(a))=(1,0,0)$. Similarly, recall the second character $b:G\longrightarrow\mathbb{Z}/2$, $b(\sigma,\epsilon)=\epsilon$. We have $b(g_1)=0$, $b(g_2)=1$, $b(g_3)=1$. Therefore $(b(g_1),b(g_2),b(g_3))=(0,1,1)$. Thus, $R_1(\beta(b))=(0,0,1)$.

\begin{lemma}[K\"unneth formula]\cite[Proposition~6.1.13]{Weibel94}
For every group $G$ and $H$, there is a split exact sequence:
\[
0
\longrightarrow\bigoplus_{p+q=n}H_p(G;\mathbb Z)\otimes H_q(H;\mathbb Z)
\longrightarrow H_n(G\times H;\mathbb Z)
\longrightarrow \bigoplus_{p+q=n-1} \operatorname{Tor}_1^{\mathbb Z} \bigl(H_p(G;\mathbb Z),H_q(H;\mathbb Z)\bigr)
\longrightarrow 0.
\]
\end{lemma}
By the K\"unneth formula over the field $\mathbb{Z}/2$, we obtain
\[
H^2(\mathfrak{S}_4\times \mathbb{Z}/2;\mathbb{Z}/2)\cong H^2(\mathfrak{S}_4;\mathbb{Z}/2) \oplus(H^1(\mathfrak{S}_4;\mathbb{Z}/2)\otimes_{\mathbb{Z}/2} H^1(\mathbb{Z}/2;\mathbb{Z}/2))\oplus H^2(\mathbb{Z}/2;\mathbb{Z}/2).
\]
Since $\dim_{\mathbb{Z}/2}H^2(\mathfrak{S}_4;\mathbb{Z}/2)=2$, $\dim_{\mathbb{Z}/2}(H^1(\mathfrak{S}_4;\mathbb{Z}/2)\otimes H^1(\mathbb{Z}/2;\mathbb{Z}/2))=1$, and $\dim_{\mathbb{Z}/2}H^2(\mathbb{Z}/2;\mathbb{Z}/2)=1$, we obtain $\dim_{\mathbb{Z}/2}H^2(G;\mathbb{Z}/2)=2+1+1=4$.

We now identify natural generators of these three summands. We explain how the relevant classes in $H^2(G;\mathbb{Z}/2)$ arise naturally from the product structure of $G$. Let $\pi_1:G\longrightarrow \mathfrak{S}_4$, $\pi_2:G\longrightarrow \mathbb{Z}/2$ be the two projections. Since $H^1(\mathfrak{S}_4;\mathbb{Z}/2)\cong \mathbb{Z}/2$ and $H^1(\mathbb{Z}/2;\mathbb{Z}/2)\cong \mathbb{Z}/2$, choose generators $a_0\in H^1(\mathfrak{S}_4;\mathbb{Z}/2)$, $b_0\in H^1(\mathbb{Z}/2;\mathbb{Z}/2)$.  Explicitly, $a_0$ is induced by the sign homomorphism:
\[
a_0(\sigma)=
\begin{cases}
0,& \sigma\in \mathfrak{A}_4,\\
1,& \sigma\notin \mathfrak{A}_4.
\end{cases}
\]
and $b_0(\epsilon)=\epsilon$. Put $a:=\pi_1^*(a_0)$, $b:=\pi_2^*(b_0)$. Thus $H^1(G;\mathbb{Z}/2)=\langle a,b\rangle_{\mathbb{Z}/2}$.

\medskip

Let $\mathbb{Z}/2=\langle t\mid t^2=1\rangle$, and let $b_0\in H^1(\mathbb{Z}/2;\mathbb{Z}/2)\cong\operatorname{Hom}(\mathbb{Z}/2,\mathbb{Z}/2)$ be the unique nonzero class, characterized by $b_0(t)=1$. We first determine the generator of $H^2(\mathbb{Z}/2;\mathbb{Z}/2)$. Consider the connecting homomorphism
\[
\beta_{2}:H^1(\mathbb{Z}/2;\mathbb{Z}/2)\longrightarrow H^2(\mathbb{Z}/2;\mathbb{Z}/2).
\]
We claim that $\beta_{2}(b_0)\neq 0$. Indeed, if $\beta_{2}(b_0)=0$, then, by exactness of the long exact sequence, the class
$b_0$ would lie in the image of
\[
H^1(\mathbb{Z}/2;\mathbb Z/4)\longrightarrow H^1(\mathbb{Z}/2;\mathbb{Z}/2),
\]
 equivalently, would admit a lift $\widetilde b_0:\mathbb{Z}/2\longrightarrow \mathbb Z/4$ such that $\widetilde b_0(t)\equiv 1\pmod 2$. This is impossible, a lift cannot exist. Therefore $\beta_{2}(b_0)\neq0$. Since $H^2(\mathbb{Z}/2;\mathbb{Z}/2)$ is one-dimensional and $\beta_{2}(b_0)\neq0$, we obtain $H^2(\mathbb{Z}/2;\mathbb{Z}/2)=\langle\beta_{2}(b_0)\rangle_{\mathbb{Z}/2}$.

\medskip

By naturality, the following diagram commutes:
\[
\begin{tikzcd}
H^1(\mathbb{Z}/2;\mathbb{Z}/2)\arrow[r,"\beta_{2}"]\arrow[d,"\pi_2^*"']&H^2(\mathbb{Z}/2;\mathbb{Z}/2)\arrow[d,"\pi_2^*"]\\
H^1(G;\mathbb{Z}/2)\arrow[r,"\beta_G"']&H^2(G;\mathbb{Z}/2).
\end{tikzcd}
\]
Consequently, $\pi_2^*\bigl(\beta_{2}(b_0)\bigr)=\beta_G\bigl(\pi_2^*(b_0)\bigr)=\beta_G(b)$. Therefore the $H^2(\mathbb{Z}/2;\mathbb{Z}/2)$-summand in $H^2(G;\mathbb{Z}/2)$ is precisely $\langle\beta_G(b)\rangle_{\mathbb{Z}/2}$.

\medskip

Since $H^1(\mathfrak{S}_4;\mathbb{Z}/2)\otimes H^1(\mathbb{Z}/2;\mathbb{Z}/2)$ is generated by $a_0\otimes b_0$,
\[
a\smile b=\pi_1^*(a_0)\smile\pi_2^*(b_0)\in H^2(G;\mathbb{Z}/2).
\]

\medskip

 Now we determine two generators for $H^2(\mathfrak{S}_4;\mathbb{Z}/2)$. Let $a_0\in H^1(\mathfrak{S}_4;\mathbb{Z}/2)\cong\operatorname{Hom}(\mathfrak{S}_4,\mathbb{Z}/2)$ be the unique nontrivial class. Then $\ker(a_0)=\mathfrak{A}_4$. Consider the short exact sequence of trivial $\mathfrak{S}_4$-modules
\[
1
\longrightarrow \mathbb{Z}/2
\overset{\iota}{\longrightarrow}\mathbb{Z}/4
\overset{p}{\longrightarrow}\mathbb{Z}/2
\longrightarrow1,
\]
where $\iota(1)=2$, $p(x)=x \pmod 2$. The associated long exact sequence in group cohomology contains the connecting homomorphism $\beta_{\mathfrak{S}_4}: H^1(\mathfrak{S}_4;\mathbb{Z}/2)\longrightarrow H^2(\mathfrak{S}_4;\mathbb{Z}/2)$. The class $\beta_{\mathfrak{S}_4}(a_0)\in H^2(\mathfrak{S}_4;\mathbb{Z}/2)$ admits a concrete interpretation in terms of group extensions. Indeed, pull back the extension $1 \longrightarrow \mathbb{Z}/2\longrightarrow \mathbb{Z}/4\longrightarrow \mathbb{Z}/2\longrightarrow 1$ along the homomorphism $a_0:\mathfrak{S}_4\longrightarrow \mathbb{Z}/2$. The resulting group is
\[
E_{a_0}:=\mathfrak{S}_4\times_{\mathbb{Z}/2}\mathbb{Z}/4=\{(\sigma,x)\in \mathfrak{S}_4\times \mathbb{Z}/4\mid a_0(\sigma)=x\pmod 2 \}.
\]
Hence we obtain a central extension
\[
1
\longrightarrow \mathbb{Z}/2
\longrightarrow E_{a_0}
\longrightarrow \mathfrak{S}_4
\longrightarrow 1.
\]
Its extension class is precisely
\[
[E_{a_0}]=\beta_{\mathfrak{S}_4}(a_0)\in H^2(\mathfrak{S}_4;\mathbb{Z}/2).
\]
Equivalently, we have a commutative diagram
\[
\begin{tikzcd}
1 \arrow[r]&\mathbb{Z}/2 \arrow[r] \arrow[d,equal]&E_{a_0} \arrow[r] \arrow[d]&\mathfrak{S}_4 \arrow[r] \arrow[d,"a_0"]&1\\
1 \arrow[r]&\mathbb{Z}/2 \arrow[r,"\iota"]&\mathbb{Z}/4 \arrow[r,"p"]&\mathbb{Z}/2 \arrow[r]&1.
\end{tikzcd}
\]
We now restrict this extension to the subgroup $\mathfrak{A}_4\subset \mathfrak{S}_4$. Since $a_0|_{\mathfrak{A}_4}=0$, naturality gives
\[
\operatorname{res}^{\mathfrak{S}_4}_{\mathfrak{A}_4}\bigl(\beta_{\mathfrak{S}_4}(a_0)\bigr)=\beta_{\mathfrak{A}_4}\bigl(\operatorname{res}^{\mathfrak{S}_4}_{\mathfrak{A}_4}(a_0)\bigr)=\beta_{\mathfrak{A}_4}(0)=0.
\]
Thus $\operatorname{res}^{\mathfrak{S}_4}_{\mathfrak{A}_4}\bigl(\beta_{\mathfrak{S}_4}(a_0)\bigr)=0$. This vanishing can also be seen directly at the level of extensions. Indeed, over $\mathfrak{A}_4$ the defining condition of $E_{a_0}$ becomes $x\equiv0\pmod2$. Hence
\[
E_{a_0}|_{\mathfrak{A}_4}=\{(\sigma,x)\in \mathfrak{A}_4\times \mathbb{Z}/4\mid x\equiv0\pmod2\}.
\]
The even elements of $\mathbb{Z}/4$ form the subgroup $\{0,2\}\cong \mathbb{Z}/2$, so $E_{a_0}|_{\mathfrak{A}_4}\cong \mathfrak{A}_4\times \mathbb{Z}/2$. Therefore the restricted extension
\[
1
\longrightarrow \mathbb{Z}/2
\longrightarrow E_{a_0}|_{\mathfrak{A}_4}
\longrightarrow \mathfrak{A}_4
\longrightarrow1
\]
is split. Thus, we can consider another class in $H^2(\mathfrak{S}_4;\mathbb{Z}/2)$, namely one whose restriction to $\mathfrak{A}_4$ corresponds to a non-split extension. This clearly gives a different cohomology class.

\medskip

In Section~\ref{oct} we obtain a central extension
\[
1
\longrightarrow \mathbb{Z}/2
\longrightarrow 2O
\overset{q}{\longrightarrow} \mathfrak{S}_4
\longrightarrow 1.
\]
Let $u\in H^2(\mathfrak{S}_4;\mathbb{Z}/2)$ denote the corresponding extension class. We now restrict this extension to $\mathfrak{A}_4\subset \mathfrak{S}_4$. Thus restriction gives
the central extension
\[
1
\longrightarrow \mathbb{Z}/2
\longrightarrow 2.\mathfrak{A}_4
\longrightarrow \mathfrak{A}_4
\longrightarrow 1.
\]
We claim that this extension is non-split. Indeed, we can check as $\mathfrak{A}_5$ case. Therefore
\[
1 \longrightarrow \mathbb{Z}/2 \longrightarrow 2.\mathfrak{A}_4 \longrightarrow \mathfrak{A}_4 \longrightarrow 1
\]
is non-split. This implies $\operatorname{res}^{\mathfrak{S}_4}_{\mathfrak{A}_4}(u)\neq 0$ in $H^2(\mathfrak{A}_4;\mathbb{Z}/2)$.

\medskip

We have shown that $\operatorname{res}^{\mathfrak{S}_4}_{\mathfrak{A}_4}\bigl(\beta_{\mathfrak{S}_4}(a_0)\bigr)=0$, whereas
$\operatorname{res}^{\mathfrak{S}_4}_{\mathfrak{A}_4}(u)\neq0$. Therefore $u\neq\beta_{\mathfrak{S}_4}(a_0)$. Moreover, since its restriction to $\mathfrak{A}_4$ is nonzero, the class $u$ itself is nonzero.
Since $\dim_{\mathbb{Z}/2}H^2(\mathfrak{S}_4;\mathbb{Z}/2)=2$, the two distinct nonzero classes $\beta_{\mathfrak{S}_4}(a_0)$, $u$ are linearly independent over $\mathbb{Z}/2$. Consequently,
\[
H^2(\mathfrak{S}_4;\mathbb{Z}/2)=\langle\beta_{\mathfrak{S}_4}(a_0),\,u\rangle_{\mathbb{Z}/2}.
\]
Pulling them back along $\pi_1$, we obtain $\pi_1^*\bigl(\beta_{\mathfrak{S}_4}(a_0)\bigr)=\beta_G(a)$
and $s:=\pi_1^*(u)$. Therefore the $H^2(\mathfrak{S}_4;\mathbb{Z}/2)$-summand of $H^2(G;\mathbb{Z}/2)$ is represented by $\beta_G(a)$ and $s$. Combining the three summands, we obtain
\[
H^2(\mathfrak{S}_4\times \mathbb{Z}/2;\mathbb{Z}/2)=\langle\beta_G(a),\beta_G(b),a\smile b,s\rangle_{\mathbb{Z}/2}.
\]
We now compute the restrictions of the remaining two generators $a\smile b$, $s:=\pi_1^*(u)$. Put $H_i:=\langle g_i\rangle$. First consider the cup-product class $a\smile b$. By naturality of the cup product,
\[
\operatorname{res}_{H_i}(a\smile b)=\operatorname{res}_{H_i}(a)\smile\operatorname{res}_{H_i}(b).
\]
Since $(a(g_1),a(g_2),a(g_3))=(1,1,0)$ and $(b(g_1),b(g_2),b(g_3))=(0,1,1)$, we immediately obtain $\operatorname{res}_{H_1}(a\smile b)=0$, $\operatorname{res}_{H_3}(a\smile b)=0$.

 For $H_2\cong \mathbb{Z}/4$, by Lemma~\ref{restrict}, $\operatorname{res}_{H_2}(a)=0$ and $\operatorname{res}_{H_i}(b)=0$. Then $\operatorname{res}_{H_2}(a\smile b)=0$. Thus
\[
R_1(a\smile b)=(0,0,0).
\]
We next compute the restriction of $s=\pi_1^*(u)$. For each $i$, naturality gives
\[
\operatorname{res}_{H_i}(s)=\operatorname{res}_{H_i}(\pi_1^*(u))=(\pi_1|_{H_i})^*(\operatorname{res}_{\pi_1(H_i)}(u)).
\]
For $g_1=((12),0)$, the subgroup $\pi_1(H_1)=\langle(12)\rangle\cong \mathbb{Z}/2$. Every lift of an element of order $2$ to the binary octahedral group $2O$ has order $4$. Hence the restricted central extension over $\langle(12)\rangle$ is non-split, and therefore $\operatorname{res}_{H_1}(s)=1$. For $g_2=((1234),1)$, we have $\pi_1(H_2)=\langle(1234)\rangle\cong \mathbb{Z}/4$. By Section~\ref{oct}, an element of order $4$ in $\mathfrak{S}_4$ lifts to an element of order $8$ in $2O$. Hence the restricted extension over
$\langle(1234)\rangle$ is non-split, and therefore $\operatorname{res}_{H_2}(s)=1$. Finally, $\pi_1(H_3)=\langle(432)\rangle\cong \mathbb{Z}/3$. By \cite[Theorem~6.5.8]{Weibel94}, the group
$H^n(\mathbb{Z}/3;\mathbb{Z}/2)$ is annihilated by $|\mathbb{Z}/3|=3$ for every $n>0$. Since multiplication by $3$ is the identity on the $\mathbb{Z}/2$-vector space $H^n(\mathbb{Z}/3;\mathbb{Z}/2)$, it follows that $3\alpha=\alpha$ and $3\alpha=0$. Thus, $H^n(\mathbb{Z}/3;\mathbb{Z}/2)=0$. In particular, $H^2(\mathbb{Z}/3;\mathbb{Z}/2)=0$. Hence $\operatorname{res}_{\langle(432)\rangle}(u)=0$, and consequently $\operatorname{res}_{H_3}(s)=0$. Thus
\[
R_1(s)=(1,1,0).
\]
Combining this with the previous computations,
\[
R_1(\beta_G(a))=(1,0,0),
\quad
R_1(\beta_G(b))=(0,0,1),
\quad
R_1(a\smile b)=(0,0,0),
\quad
R_1(s)=(1,1,0),
\]
we obtain, with respect to the ordered basis $\mathcal B=\bigl(\beta_G(a),\,\beta_G(b),\,a\smile b,\,s\bigr)$, the restriction matrix
\[
[R_1]_{\mathcal B}=
\begin{pmatrix}
1&0&0&1\\
0&0&0&1\\
0&1&0&0
\end{pmatrix}.
\]
Therefore, $\operatorname{rank} R_1=3$. It follows that $\dim H^1(C_1,\mathbb{Z}/2)^{\mathfrak{S}_4\times\mathbb{Z}/2}=1$. By \cite[Lemma 4.1]{LeeShabalin18} there is a $G$-invariant theta characteristic on $C_1$. Thus, $\#\Theta(C_1)^G=2.$
\end{proof}

\begin{proposition}   In the case $G=\mathfrak{S}_4\times \mathbb{Z}/2$, $g(C_1)=3$, $g(C_2)=25$, $T_1=[2,4,6]$, $T_2=[2,2,2,2,2,2]$, we have
$$\operatorname{Pic}(C_2)^G[2] \cong (\mathbb{Z}/2)^5 \quad\text{and}\quad  \#\Theta(C_2)^G=32.$$
\end{proposition}

\begin{proof}
For $C_2$, a generating vector is
\[
(g_1,g_2,g_3,g_4,g_5,g_6)=\bigl(((12)(34),1),((12),1),((34),1),((14)(23),1),((23),1),((14),1)\bigr),
\]
and every $g_i$ has order $2$. Let
\[
\Gamma_2=\langle y_1,y_2,y_3,y_4,y_5,y_6\mid y_1^2=\cdots=y_6^2=1,\; y_1y_2y_3y_4y_5y_6=1\rangle
\]
and let $\varphi_2:\Gamma_2\twoheadrightarrow G=\mathfrak{S}_4\times\mathbb{Z}/2$, $y_i\longmapsto g_i$.

 To compute $\#\Theta(C_2)^G$, by Proposition~\ref{dim} and ~\ref{kervarphi}, we compute
\[
\dim H^1(C,\mathbb{Z}/2)^G=\dim H^1(\Gamma,\mathbb{Z}/2)-\dim H^1(G,\mathbb{Z}/2)+\dim H^2(G,\mathbb{Z}/2)-\rank R.
\]
By Table~\ref{dim table}, we know that
\[
\dim H^1(\Gamma_2,\mathbb{Z}/2)=5,\quad \dim H^1(\mathfrak{S}_4\times \mathbb{Z}/2 ,\mathbb{Z}/2)=2, \quad \dim H^2(\mathfrak{S}_4\times \mathbb{Z}/2,\mathbb{Z}/2)=4.
\]

It is enough to compute the rank of the matrix $R_2.$

Set $H_i:=\langle g_i\rangle\cong \mathbb{Z}/2$. Since $H^2(H_i;\mathbb{Z}/2)\cong\mathbb{Z}/2$ for every $i$, the restriction map is
\[
R_2: H^2(G;\mathbb{Z}/2) \longrightarrow \bigoplus_{i=1}^6H^2(H_i;\mathbb{Z}/2) \cong (\mathbb{Z}/2)^6.
\]

Recall that
\[
H^2(G;\mathbb{Z}/2) = \langle\beta_G(a),\beta_G(b),a\smile b,s\rangle_{\mathbb{Z}/2},
\]
where $s=\pi_1^*(u)$ and $u\in H^2(\mathfrak{S}_4;\mathbb{Z}/2)$ is the class associated to the binary octahedral extension. We first compute the restrictions. Since
\[
a(\sigma,\epsilon)=
\begin{cases}
0,&\sigma\in \mathfrak{A}_4,\\
1,&\sigma\notin \mathfrak{A}_4,
\end{cases}
\]
we have
\[
\bigl(a(g_1),a(g_2),a(g_3),a(g_4),a(g_5),a(g_6)\bigr)=(0,1,1,0,1,1).
\]
Indeed, $(12)(34)$ and $(14)(23)$ are even permutations, whereas $(12)$, $(34)$, $(23)$, and $(14)$ are odd permutations. Since every $H_i\cong \mathbb{Z}/2$, by Lemma~\ref{restrict} we obtain
\[
R_2(\beta_G(a))=(0,1,1,0,1,1).
\]
On the other hand, $b(\sigma,\epsilon)=\epsilon$, and every $g_i$ has second component equal to $1$. Hence
\[
\bigl(b(g_1),b(g_2),b(g_3),b(g_4),b(g_5),b(g_6)\bigr)=(1,1,1,1,1,1).
\]
Again by Lemma~\ref{restrict},
\[
R_2(\beta_G(b))=(1,1,1,1,1,1).
\]

\medskip

We next compute the restriction of the cup-product class
$a\smile b$. By naturality of the cup product,
\[
\operatorname{res}_{H_i}(a\smile b)=\operatorname{res}_{H_i}(a)\smile\operatorname{res}_{H_i}(b).
\]
 Since $H_i\cong \mathbb{Z}/2$ for all $i=1,\cdots, 6$,
\[
R_2(a\smile b)=(0,1,1,0,1,1)\cdot(1,1,1,1,1,1)=(0,1,1,0,1,1).
\]

\medskip

Finally, we compute the restriction of $s=\pi_1^*(u)$. For every $i$, naturality gives
\[
\operatorname{res}_{H_i}(s)=\operatorname{res}_{H_i}(\pi_1^*(u))=(\pi_1|_{H_i})^*(\operatorname{res}_{\pi_1(H_i)}(u)).
\]
The first component of each $h_i$ is a nontrivial involution of $\mathfrak{S}_4$. Hence $\pi_1(H_i)\cong \mathbb{Z}/2$. Recall that every lift of a nontrivial element of order $2$ in $\mathfrak{S}_4$ to the binary octahedral group $2O$ has order $4$. Therefore the restricted central extension
\[
1
\longrightarrow \mathbb{Z}/2
\longrightarrow q^{-1}(\pi_1(H_i))
\longrightarrow\pi_1(H_i)
\longrightarrow1
\]
is non-split. Hence $\operatorname{res}_{\pi_1(H_i)}(u)\neq0$. Since $H^2(H_i;\mathbb{Z}/2)\cong\mathbb{Z}/2$, we obtain $\operatorname{res}_{H_i}(s)=1$ for every $i$. Thus
\[
R_2(s)=(1,1,1,1,1,1).
\]
Combining the above computations,
\[
\begin{aligned}
R_2(\beta_G(a))&=(0,1,1,0,1,1),\\
R_2(\beta_G(b))&=(1,1,1,1,1,1),\\
R_2(a\smile b)&=(0,1,1,0,1,1),\\
R_2(s)&=(1,1,1,1,1,1).
\end{aligned}
\]
Therefore, with respect to the ordered basis $\mathcal B=\bigl(\beta_G(a),\beta_G(b),a\smile b,s\bigr)$, the restriction matrix is
\[
[R_2]_{\mathcal B}=
\begin{pmatrix}
0&1&0&1\\
1&1&1&1\\
1&1&1&1\\
0&1&0&1\\
1&1&1&1\\
1&1&1&1
\end{pmatrix}.
\]
The two vectors $(0,1,1,0,1,1)$, $(1,1,1,1,1,1)$ are linearly independent over $\mathbb{Z}/2$. Therefore,
\[
\operatorname{rank} R_2=2.
\]
It follows that  $\dim H^1(C_2,\mathbb{Z}/2)^{\mathfrak{S}_4\times\mathbb{Z}/2}=5$. By \cite[Lemma 4.3]{LeeShabalin18} there is a $\mathfrak{S}_4\times\mathbb{Z}/2$-invariant theta characteristic on $C_1$.
Thus, $\#\Theta(C_2)^G=32.$
\end{proof}

 We now determine $\operatorname{Pic}(C_i)^G$ and the parity of the intersection form of $S$.

\medskip

In this case, by Remark~\ref{Schur} and Proposition \ref{Pic(G;C)} we have
$$\operatorname{Hom}(\mathfrak{S}_4\times \mathbb{Z}/2,\mathbb{C}^*)=(\mathbb{Z}/2)^2,\,\,\, H^{2}(\mathfrak{S}_4\times \mathbb{Z}/2, \mathbb{C}^*)=(\mathbb{Z}/2)^2$$
$$\operatorname{Pic}(G;C_1)=  4\mathbb{Z} \oplus (\mathbb{Z}/2)^2,\,\,\, \operatorname{Pic}(G;C_2)= 24\mathbb{Z} \oplus (\mathbb{Z}/2)^5.$$
The exact sequence \eqref{eq:exact} yields
\begin{subequations}\label{eq:row4-sequences}
\begin{equation}\label{eq:row4-C1-concrete}
0
\longrightarrow (\mathbb{Z}/2)^2
\longrightarrow 4\mathbb{Z} \oplus (\mathbb{Z}/2)^2
\longrightarrow \operatorname{Pic}(C_1)^G
\longrightarrow (\mathbb{Z}/2)^2
\longrightarrow 0 .
\end{equation}
\begin{equation}\label{eq:row4-C2-concrete}
0
\longrightarrow (\mathbb{Z}/2)^2
\longrightarrow 24\mathbb{Z} \oplus (\mathbb{Z}/2)^{5}
\longrightarrow \operatorname{Pic}(C_{2})^{G}
\longrightarrow (\mathbb{Z}/2)^2
\longrightarrow 0 .
\end{equation}
\end{subequations}
 Since $\operatorname{Pic}(C_1)^G[2]= \mathbb{Z}/2$ and $\operatorname{Pic}(C_2)^G[2]= (\mathbb{Z}/2)^5$,  we infer that
$$\operatorname{Pic}(C_1)^G= 2\mathbb{Z}\oplus \mathbb{Z}/2 \quad\text{and}\quad \operatorname{Pic}(C_2)^G= 24\mathbb{Z}\oplus (\mathbb{Z}/2)^5.$$
By \eqref{eq:row4-C1-concrete} we can choose a $G$-invariant line bundle $L_1\in \operatorname{Pic}(C_1)^G$ of degree $2=g(C_1)-1$ which is not $G$-linearized, and
  by  \eqref{eq:row4-C2-concrete} we can choose a $G$-invariant torsion line bundle $\tau \in \operatorname{Pic}(C_2)^G$ whose obstruction class in the Schur multiplier coincides with that of $L_1$.

  Therefore, the external tensor product
$L_1\boxtimes \tau$
admits a $G$-linearization, that is, it is a $G$-linearized line bundle on $C_1\times C_2$. This line bundle has bidegree $(2,0)$.
On the other hand, from \eqref{eq:row4-C2-concrete} we can choose a $G$-linearized line bundle on $C_2$ of degree $24=g(C_2)-1$. Its pullback on $C_1\times C_2$ is a $G$-linearized line bundle of bidegree $(0,24)$.

Therefore, by Strategy (1) the intersection form of $S$ is even.

\bigskip

\section{$G=G(32)$, $g(C_1)=9$, $g(C_2)=5$, $T_1=[2,2,4,4]$, $T_2=[2,2,2,4]$; odd fake quadrics}

Recall $G(32)$ as $G=N\rtimes_{\Phi}\mathbb{Z}/2, N=(\mathbb{Z}/2)^4$. Write an element of $N$ as a column vector $v=(v_1,v_2,v_3,v_4)^t, v_i\in\mathbb{Z}/2$ and let the nontrivial element of $\mathbb{Z}/2$ act by
\[
\Phi_1=
\begin{pmatrix}
1&0&0&0\\
0&1&0&0\\
1&0&1&0\\
0&1&0&1
\end{pmatrix},
\quad
\text{Note $\Phi_0$ is identity.}
\]

Thus, $\Phi_1(v_1,v_2,v_3,v_4)^t=(v_1,v_2,v_1+v_3,v_2+v_4)^t$.

We write the elements of $G$ in the form $(v,\epsilon)$, $v\in N$, $\epsilon\in\mathbb{Z}/2$ with multiplication $(v,\epsilon)(w,\delta)=\bigl(v+\Phi_\epsilon(w),\epsilon+\delta\bigr)$.

Let  $e_i=((0,\ldots,\underbrace{1}_{i\text{-th}},\ldots,0)^t,0)$, $1\leq i\leq4$ and $t=((0,0,0,0)^t,1)$. Then $e_i^2=t^2=1$, $[e_i,e_j]=1$ and
\[
te_1t^{-1}=e_1e_3,\quad
te_2t^{-1}=e_2e_4,\quad
te_3t^{-1}=e_3,\quad
te_4t^{-1}=e_4.
\]
In particular, $[t,e_1]=e_3$, $[t,e_2]=e_4$.

Also, $G^{\mathrm{ab}}=G/[G,G]\cong(\mathbb{Z}/2)^3$.

\begin{proposition}\label{G32thetaC1}  In the case $G=G(32)$, $g(C_1)=9$, $g(C_2)=5$, $T_1=[2,2,4,4]$, $T_2=[2,2,2,4]$, we have
$$\operatorname{Pic}(C_1)^G[2]= (\mathbb{Z}/2)^3 \quad\text{and}\quad  \#\Theta(C_1)^G=8.$$
\end{proposition}
\begin{proof}
For $C_1$, a generating vector is given by the quadruple
\[
\begin{aligned}
g_1&=((1,1,1,0)^t,0)=e_1e_2e_3,\\
g_2&=((1,0,0,0)^t,0)=e_1,\\
g_3&=((1,1,1,0)^t,1)=e_1e_2e_3t,\\
g_4&=((1,0,1,0)^t,1)=e_1e_3t.
\end{aligned}
\]
with order $2,2,4,4$. Let
\[
\Gamma_1=\langle x_1,x_2,x_3,x_4 \mid x_1^2=x_2^2=x_3^4=x_4^4=1,\;x_1x_2x_3x_4=1\rangle
\]
and let $\varphi_1:\Gamma_1\twoheadrightarrow G=G(32)$, $x_i\longmapsto g_i$.

 To compute $\#\Theta(C_1)^G$, by Proposition~\ref{dim} and ~\ref{kervarphi}, we compute
\[
\dim H^1(C,\mathbb{Z}/2)^G=\dim H^1(\Gamma,\mathbb{Z}/2)-\dim H^1(G,\mathbb{Z}/2)+\dim H^2(G,\mathbb{Z}/2)-\rank R.
\]
By Table~\ref{dim table}, we know that
\[
\dim H^1(\Gamma_1,\mathbb{Z}/2)=3,\quad \dim H^1(G(32) ,\mathbb{Z}/2)=3, \quad \dim H^2(G(32),\mathbb{Z}/2)=7.
\]
It is enough to compute the rank of the matrix $R$. Since $G^{\mathrm{ab}}\cong(\mathbb{Z}/2)^3$, we have $H^1(G;\mathbb{Z}/2)=\operatorname{Hom}(G,\mathbb{Z}/2)\cong(\mathbb{Z}/2)^3$. We choose the following basis.

 Define characters $\chi_1,\chi_2,\chi_3:G\longrightarrow\mathbb{Z}/2$ by
\[
\chi_1(v,\epsilon)=v_1,\quad
\chi_2(v,\epsilon)=v_2,\quad
\chi_3(v,\epsilon)=\epsilon.
\]
These characters may be understood as follows. Since $[G,G]=\langle e_3,e_4\rangle$, we have $G^{\mathrm{ab}}=G/[G,G]=\langle \overline e_1,\overline e_2,\overline t\rangle\cong(\mathbb Z/2)^3$. Hence the three coordinate projections on $G^{\mathrm{ab}}$ define characters $\chi_1$, $\chi_2$, $\chi_3$.

 Their values on the four stabilizer generators are
\[
\begin{array}{c|cccc}
 &g_1&g_2&g_3&g_4\\
\hline
\chi_1&1&1&1&1\\
\chi_2&1&0&1&0\\
\chi_3&0&0&1&1.
\end{array}
\]
These three vectors are linearly independent. Indeed, if $a(1,1,1,1)+b(1,0,1,0)+c(0,0,1,1)=0$, then the second coordinate gives $a=0$, the first coordinate then gives $b=0$, and the third coordinate gives $c=0$. Thus, they form a basis, i.e.
\[
H^1(G;\mathbb{Z}/2)=\langle\chi_1,\chi_2,\chi_3\rangle_{\mathbb{Z}/2}.
\]

We compute the connecting homomorphism $\beta : H^1(G;\mathbb{Z}/2)\longrightarrow H^2(G;\mathbb{Z}/2)$, i.e. $\beta(\chi_1)$, $\beta(\chi_2)$, $\beta(\chi_3).$ Put $H_i=\langle g_i\rangle, 1\leq i\leq4$. Then $H_1\cong H_2\cong \mathbb{Z}/2$, $H_3\cong H_4\cong \mathbb{Z}/4$. Since $H^2(\mathbb{Z}/2;\mathbb{Z}/2)\cong\mathbb{Z}/2$, $H^2(\mathbb{Z}/4;\mathbb{Z}/2)\cong\mathbb{Z}/2$, we identify $\bigoplus_{i=1}^4H^2(H_i;\mathbb{Z}/2)\cong(\mathbb{Z}/2)^4$. Define
\[
R_1:H^2(G;\mathbb{Z}/2)\longrightarrow\bigoplus_{i=1}^4 H^2(\langle g_i\rangle;\mathbb{Z}/2)\cong(\mathbb{Z}/2)^4.
\]
Naturality gives $\operatorname{res}_{H_i}\beta(\chi)=\beta(\chi|_{H_i})$.  Using  Lemma~\ref{restrict} and the above table of character values, we obtain
\[
R_1(\beta(\chi_1))=(1,1,0,0),
\quad
R_1(\beta(\chi_2))=(1,0,0,0),
\quad
R_1(\beta(\chi_3))=(0,0,0,0).
\]
We therefore need more additional classes.
Let
\[
D_4=\langle r,s\mid r^4=s^2=1,\ srs=r^{-1}\rangle.
\]
Consider
\[
\widetilde D_4=\langle R,S\mid R^8=S^2=1,\ SRS=R^{-1}\rangle.
\]
There is a central extension
\[
1
\longrightarrow\langle R^4\rangle
\longrightarrow\widetilde D_4
\overset{p}{\longrightarrow}D_4
\longrightarrow1,
\]
where $p(R)=r$, $p(S)=s$, and $\langle R^4\rangle\cong \mathbb{Z}/2$. Let $u\in H^2(D_4;\mathbb{Z}/2)$ be the corresponding extension class. We first determine the restriction of $u$ to the cyclic subgroups of $D_4$.

\begin{lemma}
Let $u\in H^2(D_4;\mathbb{Z}/2)$ be the extension class defined above. Then $\operatorname{res}_{\langle r\rangle}(u)\neq0$. More generally, if $r^j$ generates the rotation subgroup
$\langle r\rangle\cong \mathbb{Z}/4$, then $\operatorname{res}_{\langle r^j\rangle}(u)\neq0$. On the other hand, for every reflection $r^ks$, $\operatorname{res}_{\langle r^ks\rangle}(u)=0$.
\end{lemma}

\begin{proof}
Recall that restriction of an extension class has a direct
interpretation in terms of pullback extensions. If $K\leq D_4$, then $\operatorname{res}_K(u)\in H^2(K;\mathbb{Z}/2)$ is represented by the restricted central extension
\[
1
\longrightarrow \mathbb{Z}/2
\longrightarrow p^{-1}(K)
\longrightarrow K
\longrightarrow 1.
\]
Consequently, $\operatorname{res}_K(u)=0$ if and only if this restricted extension splits. We first consider the rotation subgroup $K=\langle r\rangle\cong \mathbb{Z}/4$. Its inverse image under $p$ is
\[
p^{-1}(\langle r\rangle)=\langle R\rangle\cong \mathbb{Z}/8.
\]
Hence the restricted extension is
\[
1
\longrightarrow\langle R^4\rangle\cong \mathbb{Z}/2
\longrightarrow\langle R\rangle\cong \mathbb{Z}/8
\longrightarrow\langle r\rangle\cong \mathbb{Z}/4
\longrightarrow1.
\]

We claim that this extension is non-split. Suppose, to the contrary,
that it splits. Then there would exist a group homomorphic section $\sigma:\mathbb{Z}/4\longrightarrow \mathbb{Z}/8$. In particular, $\sigma(r)$ would be a lift of $r$ having order $4$. The two lifts of $r$ in $\mathbb{Z}/8=\langle R\rangle$ are $R$ and $R^5$, since $p(R)=r$ and $p(R^5)=r^5=r$. However,
\[
ord(R)=8
\quad\text{and}\quad
ord(R^5)=8.
\]
Thus there is no element of order $4$ mapping to $r$, and hence no section exists. Therefore the extension
\[
1\longrightarrow \mathbb{Z}/2
\longrightarrow \mathbb{Z}/8
\longrightarrow \mathbb{Z}/4
\longrightarrow1
\]
is non-split. It follows that $\operatorname{res}_{\langle r\rangle}(u)\neq0$. Since $r^3$ is also a generator of the same rotation subgroup, $\langle r^3\rangle=\langle r\rangle$, and therefore $\operatorname{res}_{\langle r^3\rangle}(u)\neq0$. Since $H^2(\mathbb{Z}/4;\mathbb{Z}/2)\cong\mathbb{Z}/2$, we identify its unique nonzero class with $1$. Thus we may write $\operatorname{res}_{\langle r\rangle}(u)=1$. We now consider a reflection $a=r^ks$. Since $(r^ks)^2=1$, the subgroup $\langle r^ks\rangle$ is cyclic of order $2$. The element $R^kS$ is a lift of $r^ks$, because $p(R^kS)=r^ks$. We now verify that this lift also has order $2$. From the relation $SRS=R^{-1}$, we obtain, for every integer $k$, $SR^kS=R^{-k}$. Hence
\[
\begin{aligned}
(R^kS)^2=R^kSR^kS=R^k(SR^kS)=R^kR^{-k}=1.
\end{aligned}
\]
Therefore the assignment $r^ks\longmapsto R^kS$ defines a section $\langle r^ks\rangle\longrightarrow p^{-1}(\langle r^ks\rangle)$. Thus the restricted extension over every reflection subgroup splits, and hence $\operatorname{res}_{\langle r^ks\rangle}(u)=0$.
\end{proof}

We continue the proof of Proposition ~\ref{G32thetaC1}.

The restriction behavior of $u$ can be summarized as follows:
\[
\begin{array}{c|c}
K\leq D_4 & \operatorname{res}_K(u)\\
\hline
\langle r\rangle\cong \mathbb{Z}/4 & 1\\
\langle r^3\rangle=\langle r\rangle & 1\\
\langle r^ks\rangle\cong \mathbb{Z}/2 & 0.
\end{array}
\]
We now use this class to construct an element of $H^2(G;\mathbb{Z}/2)$.
First, we construct $D_4$. Consider $N_1:=\langle e_2,e_4\rangle$. We identify $G/N_1$ with $D_4$. We choose
\[
\bar e_1\bar t\longmapsto r
\quad
\bar e_1\longmapsto rs,
\quad
\bar e_3\longmapsto r^2,
\quad
\bar t\longmapsto s.
\]
Thus, these three elements give the quotient $G/\langle e_2,e_4\rangle \cong D_4$. Put $N_2:=\langle e_1,e_3\rangle$. We identify $G/N_2$ with $D_4$ by
\[
\bar e_2\bar t\longmapsto r
\quad
\bar e_2\longmapsto rs,
\quad
\bar e_4\longmapsto r^2,
\quad
\bar t\longmapsto s.
\]
We now define two homomorphisms $q_1,q_2:G\longrightarrow D_4$ by
\[
\begin{array}{c|ccccc}
&e_1&e_2&e_3&e_4&t\\
\hline
q_1&rs&1&r^2&1&s\\
q_2&1&rs&1&r^2&s.
\end{array}
\]
Define $s_1:=q_1^*(u)\in H^2(G;\mathbb{Z}/2)$. Let $H_i=\langle g_i\rangle$, $1\leq i\leq4$. By naturality of restriction and pullback, $\operatorname{res}_{H_i}(s_1)=\operatorname{res}_{H_i}(q_1^*(u))=(q_1|_{H_i})^*(u)$. Therefore, in order to determine $\operatorname{res}_{H_i}(s_1)$, it suffices to determine whether $q_1(H_i)$ is a reflection subgroup or the rotation subgroup of $D_4$. We compute the four stabilizer generators separately.

For $g_1=e_1e_2e_3$,
\[
\begin{aligned}
q_1(g_1)=q_1(e_1)q_1(e_2)q_1(e_3)=(rs)\cdot 1\cdot r^2=r^3s.
\end{aligned}
\]
Therefore $q_1(H_1)=\langle r^3s\rangle$ is a reflection subgroup. By the preceding lemma, $\operatorname{res}_{q_1(H_1)}(u)=0$. Consequently, $\operatorname{res}_{H_1}(s_1)=0.$

For $g_2=e_1$, $q_1(g_2)=rs$. Since $rs$ is a reflection, $q_1(H_2)=\langle rs\rangle$ is a reflection subgroup. Hence $\operatorname{res}_{H_2}(s_1)=0$.

For $g_3=e_1e_2e_3t$,
\[
\begin{aligned}
q_1(g_3)=q_1(e_1)q_1(e_2)q_1(e_3)q_1(t)=(rs)\cdot1\cdot r^2\cdot s=r^3.
\end{aligned}
\]
 Since the element $r^3$ has order $4$, $q_1(H_3)=\langle r^3 \rangle$ is the  rotation subgroup of $D_4$.

  Hence $\operatorname{res}_{\langle r\rangle}(u)\neq0$. Pulling a nonzero cohomology class back along the isomorphism $q_1|_{H_3}$ gives a nonzero class. Therefore $\operatorname{res}_{H_3}(s_1)\neq0$.

\medskip

For $g_4=e_1e_3t$,
\[
\begin{aligned}
q_1(g_4)=q_1(e_1)q_1(e_3)q_1(t)=(rs)r^2s=r^3.
\end{aligned}
\]
Hence
\[
q_1(H_4)=\langle r^3\rangle=\langle r\rangle\cong \mathbb{Z}/4.
\]
Since $H_4\cong \mathbb{Z}/4$, the map $q_1|_{H_4}:H_4\longrightarrow\langle r\rangle$ is again an isomorphism. Therefore the same argument as for $H_3$ gives $\operatorname{res}_{H_4}(s_1)\neq0.$
Combining the four calculations,
\[
\begin{aligned}
R_1(s_1)=\bigl(\operatorname{res}_{H_1}(s_1),\operatorname{res}_{H_2}(s_1),\operatorname{res}_{H_3}(s_1),\operatorname{res}_{H_4}(s_1)\bigr)=(0,0,1,1).
\end{aligned}
\]
We now perform the corresponding calculation for $q_2$. Define $s_2:=q_2^*(u)$. Similarly,
\[
\begin{aligned}
q_2(g_1)&=q_2(e_1e_2e_3)=rs,\\
q_2(g_2)&=1,\\
q_2(g_3)&=q_2(e_1e_2e_3t)=(rs)s=r,\\
q_2(g_4)&=q_2(e_1e_3t)=s.
\end{aligned}
\]
The first and the fourth are reflections, the second is trivial, and
the third generates $\langle r\rangle\cong \mathbb{Z}/4$. Hence
\[
R_1(s_2)=(0,0,1,0).
\]
We have therefore obtained
\[
\begin{array}{c|c}
\text{class} & R_1(\text{class})\\
\hline
\beta(\chi_1) & (1,1,0,0)\\
\beta(\chi_2) & (1,0,0,0)\\
\beta(\chi_3) & (0,0,0,0)\\
s_1 & (0,0,1,1)\\
s_2 & (0,0,1,0).
\end{array}
\]
We claim that $\beta(\chi_1)$, $\beta(\chi_2)$, $\beta(\chi_3)$, $s_1$, $s_2$ are linearly independent in $H^2(G;\mathbb{Z}/2)$.
If $a_1\beta(\chi_1)+a_2\beta(\chi_2)+a_3\beta(\chi_3)+a_4s_1+a_5s_2=0$,
then, by applying $R_1$, we obtain
$a_1(1,1,0,0)+a_2(1,0,0,0)+a_4(0,0,1,1)+a_5(0,0,1,0)=0,
$
thus we get $a_1=0$, $a_2=0$, $a_4=0$, $a_5=0$ and $a_3\beta(\chi_3)=0$.

Note that in this case the Bockstein homomorphism $\beta:H^1(G;\mathbb{Z}/2)\longrightarrow H^2(G;\mathbb{Z}/2)$ is injective. Indeed, since the coefficient modules are trivial, we have $H^1(G;\mathbb{Z}/4)=\operatorname{Hom}(G,\mathbb{Z}/4)$, $H^1(G;\mathbb{Z}/2)=\operatorname{Hom}(G,\mathbb{Z}/2)$. Moreover, every homomorphism $f:G\longrightarrow\mathbb{Z}/4$ factors through the abelianization:
\[
G
\longrightarrow G^{\mathrm{ab}}\cong (\mathbb{Z}/2)^3
\longrightarrow \mathbb{Z}/4.
\]
Since $(\mathbb{Z}/2)^3$ is $2$-elementary, the image of $f$ must be contained in the unique subgroup $\{0,2\}\subset\mathbb{Z}/4$. Therefore, after reduction modulo $2$, the map $\rho^*: H^1(G;\mathbb{Z}/4)\longrightarrow H^1(G;\mathbb{Z}/2)$ is the zero map. Consequently, Im$\rho^*=\ker\beta=0$, Since $\chi_3\neq0$, we have $\beta(\chi_3)\neq0$. Therefore $a_3=0$. Hence these are linearly independent. Since $\dim_{\mathbb{Z}/2}H^2(G;\mathbb{Z}/2)=7$, we may choose two additional classes $v_6,v_7\in H^2(G;\mathbb{Z}/2)$ such that
\[
\mathcal B=\bigl(\beta(\chi_1),\beta(\chi_2),\beta(\chi_3),s_1,s_2,v_6,v_7\bigr)
\]
is a basis of $H^2(G;\mathbb{Z}/2)$. Thus, the restriction map $R_1:H^2(G;\mathbb{Z}/2)\longrightarrow(\mathbb{Z}/2)^4$ is represented by a matrix of the form
\[
[R_1]_{\mathcal B}=
\begin{pmatrix}
1&1&0&0&0&*&*\\
1&0&0&0&0&*&*\\
0&0&0&1&1&*&*\\
0&0&0&1&0&*&*
\end{pmatrix}.
\]
Therefore
$$
\operatorname{rank} R_1=4.
$$
It follows that $\dim H^1(C_1;\mathbb{Z}/2)^G=3$. Proposition ~\ref{G32C1} gives a $G(32)$-invariant theta characteristic on $C_1$. Thus $\#\Theta(C_1)^G=8.$
\end{proof}

\begin{proposition}\label{G32C1}
There exists a $G(32)$-invariant theta characteristic on $C_1$.
\end{proposition}
This follows from Lemmas ~\ref{pull line} - ~\ref{G32C1-last}.

\begin{lemma}\label{pull line}
Let $C$ be a smooth projective curve with a $G$-action, and let
$H\triangleleft G$ act freely on $C$. Set $q:C\longrightarrow B:=C/H$. If $\theta\in\Theta(B)^{G/H}$, then $q^*\theta\in\Theta(C)^G$.
\end{lemma}

\begin{proof}
Since $H$ acts freely, $q$ is an \'etale covering. Hence $K_C\cong q^*K_B$. Therefore, if $\theta^{\otimes2}\cong K_B$, then $(q^*\theta)^{\otimes2}\cong q^*(\theta^{\otimes2})\cong q^*K_B\cong K_C$. Thus $q^*\theta$ is a theta characteristic on $C$. Since $H\triangleleft G$, the $G$-action on $C$ descends to a
$G/H$-action on $B$. For $g\in G$, denoting its image in $G/H$ by $\bar g$, we have $q\circ g=\bar g\circ q$. Consequently, $g^*(q^*\theta)=q^*(\bar g^*\theta)\cong q^*\theta$, because $\theta$ is $G/H$-invariant. Therefore $q^*\theta\in\Theta(C)^G$.
\end{proof}

We first pass to a suitable \'etale quotient.

\begin{lemma}\label{lem:H1-free}
The normal subgroup $N_1:=\langle e_2,e_4\rangle\cong(\mathbb{Z}/2)^2$ acts freely on $C_1$.
\end{lemma}

\begin{proof}
Let $\pi:C_1\longrightarrow C_1/G\simeq\mathbb{P}^1$ be the quotient map, and let $(g_1,g_2,g_3,g_4)$ be the stabilizer generators corresponding to the four branch points. Their orders are $ord(g_1)=ord(g_2)=2$, $ord(g_3)=ord(g_4)=4$. Note that $N_1$ acts freely on $C_1$ iff $N_1\cap\langle g_i\rangle=\{1\}$, $i=1,\ldots,4$. We now use the quotient map $q:G\longrightarrow G/N_1$. Restricting $q$ to the cyclic subgroup $\langle g_i\rangle$, we obtain $q|_{\langle g_i\rangle}:\langle g_i\rangle\longrightarrow G/N_1$. Its kernel is
\[
\begin{aligned}
\ker\left(q|_{\langle g_i\rangle}\right)=\{g\in\langle g_i\rangle\mid q(g)=1\}=\{g\in\langle g_i\rangle\mid g\in N_1\}=N_1\cap\langle g_i\rangle.
\end{aligned}
\]
Therefore $N_1\cap\langle g_i\rangle=\{1\}$ iff $q|_{\langle g_i\rangle}$ is injective. Since $\langle g_i\rangle$ is cyclic, injectivity is equivalent to saying
that the order of its generator is preserved under the quotient map. Indeed, by the first isomorphism theorem, $\frac{\langle g_i\rangle}{N_1\cap\langle g_i\rangle}\cong\langle\bar g_i\rangle$, where $\bar g_i:=q(g_i)$. Taking orders, we obtain $ord(\bar g_i)=\frac{ord(g_i)}{\left|N_1\cap\langle g_i\rangle\right|}$. Thus $N_1\cap\langle g_i\rangle=\{1\}$ is equivalent to $ord(\bar g_i)=ord(g_i)$. We now identify the quotient group. As shown above, $G/N_1\simeq D_4=\langle r,s\mid r^4=s^2=1,\ srs=r^{-1}\rangle,$ and the images of the four stabilizer generators are $\bar g_1=r^3s$, $\bar g_2=rs$, $\bar g_3=r^3$, $\bar g_4=r^3$. Therefore $ord(\bar g_1)=2=ord(g_1)$, $ord(\bar g_2)=2=ord(g_2)$, and $ord(\bar g_3)=4=ord(g_3)$, $ord(\bar g_4)=4=ord(g_4)$. This implies $N_1\cap\langle g_i\rangle=\{1\}$ for all $i=1,\ldots,4$. Therefore $N_1$ acts freely on $C_1$.
\end{proof}

Set $B_1:=C_1/N_1$ and let $q_1:C_1\longrightarrow B_1$ be the quotient map. Since $q_1$ is \'etale of degree $4$, Riemann--Hurwitz gives $2g(C_1)-2=4\bigl(2g(B_1)-2\bigr)$. Since $g(C_1)=9$, we have $16=4\bigl(2g(B_1)-2\bigr)$, and hence $g(B_1)=3$. We have $G/N_1\cong D_4$, where $D_4=\langle r,s\mid r^4=s^2=1,\ srs=r^{-1}\rangle$. Under this identification, the generating vector of the $D_4$-action on $B_1$ is $(r^3s,rs,r^3,r^3)$.
We will construct explicitly a $D_4$-invariant theta characteristic on $B_1$.

\begin{lemma}
The curve $B_1$ is a hyperelliptic curve.
\end{lemma}
\begin{proof}
Let $R:=\langle r\rangle\cong\mathbb{Z}/4$. Since $R\triangleleft D_4$, the quotient map $B_1\longrightarrow B_1/D_4$ factors as $B_1\overset{\rho}{\longrightarrow}B_1/R\longrightarrow B_1/D_4\cong\mathbb{P}^1$. We first determine the ramification of the cyclic covering $\rho:B_1\longrightarrow B_1/R$. Consider one of the two order-$4$ branch points $P\in B_1/D_4$. Its stabilizer is $R=\langle r\rangle$. Since $|D_4|=8$ and $|R|=4$, the fibre of $B_1\longrightarrow B_1/D_4$ over $P$ consists of $|D_4|/|R|=2$ points. Let these points be $p$ and $p'$. Since the stabilizer is $R$, both $p$ and $p'$ are fixed by $R$. Therefore they remain two distinct points in $B_1/R$. Since there are two order-$4$ branch points on $B_1/D_4$, the covering $\rho$ has exactly $2\cdot2=4$ ramification points, each with ramification index $4$. Applying Riemann--Hurwitz to the degree-$4$ map $\rho$, we obtain $2g(B_1)-2=4\bigl(2g(B_1/R)-2\bigr)+4(4-1)$. Since $g(B_1)=3$, this becomes $4=4\bigl(2g(B_1/R)-2\bigr)+12$, and hence $g(B_1/R)=0$. Therefore $B_1/R\cong\mathbb{P}^1$.

\medskip

We next determine the local monodromies of the four branch points of $\rho$. Since $\rho$ is an $R$-Galois cover and $R=\langle r\rangle\cong\mathbb{Z}/4$, the local monodromy at a point with ramification index $4$ must have order $4$. The only elements of order $4$ in $R$ are $r$ and $r^3=r^{-1}$. There are two points $a,a'\in B_1/R$ lying over $P$. Since $D_4/R\cong\mathbb{Z}/2$ and the nontrivial element is represented by $s$, the element $s$ exchanges these two points. Thus, the four local monodromies are $r,r^{-1},r,r^{-1}$. After renumbering the four branch points, we may write them as $r,r,r^{-1},r^{-1}$.

\medskip

 Denote the two branch points with local monodromy $r$ by $a_1,a_2$ and the two branch points with local monodromy $r^{-1}$ by $b_1,b_2$.

  Since $B_1/R\cong\mathbb{P}^1$, choose a local coordinate $x$ on $B_1/R\cong \mathbb{P}^1$. Then its function field is $\mathbb{C}(x)$. Let $L:=\mathbb{C}(B_1)$ be the function field of $B_1$. The quotient map $\rho:B_1\longrightarrow B_1/R$ induces an inclusion of function fields $\mathbb{C}(x)\subset L$. Since $\rho$ is an $R$-Galois covering and $R=\langle r\rangle\cong\mathbb{Z}/4$, the extension $L/\mathbb{C}(x)$ is a cyclic Galois extension of degree $4$ with $\operatorname{Gal}(L/\mathbb{C}(x))=\langle r\rangle$. We regard $r$ as a $\mathbb{C}(x)$-linear automorphism of the $4$-dimensional $\mathbb{C}(x)$-vector space $L$. Since $r$ has order $4$, we have $r^4=1$, so every eigenvalue of $r$ is a fourth root of unity. Replacing $r$ by $r^{-1}$ if necessary, we may choose a nonzero element $y\in L$ such that $r(y)=iy$. Then, we have $r(y^4)=r(y)^4=(iy)^4=y^4$. Thus $y^4$ is invariant under the action of $R$. Since the fixed field of $R$ is $\mathbb{C}(B_1)^R=\mathbb{C}(B_1/R)=\mathbb{C}(x)$, it follows that $y^4\in\mathbb{C}(x)$. Therefore there exists a rational function $f(x)\in\mathbb{C}(x)$ such that $y^4=f(x)$.

\medskip

To determine $f(x)$, recall that if locally $f(x)=(x-a)^m u(x)$ with $u(a)\neq0$ gives local monodromy $r^m$. Hence exponent $1$ gives monodromy $r$, while exponent $-1$ gives monodromy $r^{-1}$. We may therefore choose $f(x)$ with simple zeros at $a_1,a_2$ and simple poles at $b_1,b_2$. Thus the covering can be written as $y^4=\frac{(x-a_1)(x-a_2)}{(x-b_1)(x-b_2)}$, that is, $(y^4-1)x^2+\bigl((a_1+a_2)-y^4(b_1+b_2)\bigr)x+\bigl(y^4b_1b_2-a_1a_2\bigr)=0$. Thus $x$ is quadratic over $\mathbb{C}(y)$, so $[L:\mathbb C(y)]\leq2$. The degree cannot be $1$, since otherwise $L=\mathbb C(y)$ and $B_1$ would be rational, contrary to $g(B_1)=3$. Therefore $[L:\mathbb C(y)]=2$. and hence the function $y$ induces a morphism $\eta:B_1\longrightarrow\mathbb P^1$ of degree $2$. Therefore $B_1$ is hyperelliptic.
\end{proof}

\begin{lemma}\label{G32C1-last}
The curve $B_1$ admits a $D_4$-invariant theta characteristic.
\end{lemma}

\begin{proof}
Let $\eta:B_1\longrightarrow \mathbb P^1$ be the hyperelliptic covering and let $A$ denote the divisor class of a fiber of $\eta$. Since $g(B_1)=3$, the canonical divisor  $K_{B_1}\equiv (g(B_1)-1)A=2A$; see \cite[(3.2)]{KallelSjerve10}. Hence $\theta_1:=\mathcal O_{B_1}(A)$ is a theta characteristic. Let $\iota$ be the hyperelliptic involution. Since $\iota$ is central in $\operatorname{Aut}(B_1)$, every automorphism $\varphi\in \operatorname{Aut}(B_1)$ descends to an automorphism $\bar\varphi$ of $B_1/\langle \iota\rangle\cong\mathbb{P}^1$, so that $\eta\circ\varphi=\bar\varphi\circ\eta$; see \cite[Section~3.2]{KallelSjerve10}. Thus $\varphi$ sends fibers of $\eta$ to fibers of $\eta$. Consequently, $\varphi^*\theta_1\cong\mathcal O_{B_1}(\varphi^*A)\cong\mathcal O_{B_1}(A)=\theta_1$. In particular, $\theta_1\in\Theta(B_1)^{D_4}$.
\end{proof}

\begin{proposition}  In the case $G=G(32)$, $g(C_1)=9$, $g(C_2)=5$, $T_1=[2,2,4,4]$, $T_2=[2,2,2,4]$, we have
$$\operatorname{Pic}(C_2)^G[2]= (\mathbb{Z}/2)^3 \quad\text{and}\quad  \#\Theta(C_2)^G=8.$$
\end{proposition}

\begin{proof}

For $C_2$, we use the generating vector given by the quadruple
\[
\begin{aligned}
g_1&=((0,0,1,1)^t,1)=e_3e_4t,\\
g_2&=((1,1,1,1)^t,0)=e_1e_2e_3e_4,\\
g_3&=((1,0,1,1)^t,0)=e_1e_3e_4,\\
g_4&=((0,1,1,1)^t,1)=e_2e_3e_4t,
\end{aligned}
\]
with orders $2,2,2,4$. Let
\[
\Gamma_2=\langle y_1,y_2,y_3,y_4 \mid y_1^2=y_2^2=y_3^2=y_4^4=1,\;y_1y_2y_3y_4=1\rangle,
\]
and let $\varphi_2:\Gamma_2\twoheadrightarrow G$, $y_i\longmapsto g_i$.

 To compute $\#\Theta(C_2)^G$, by Proposition~\ref{dim} and ~\ref{kervarphi}, we compute
\[
\dim H^1(C,\mathbb{Z}/2)^G=\dim H^1(\Gamma,\mathbb{Z}/2)-\dim H^1(G,\mathbb{Z}/2)+\dim H^2(G,\mathbb{Z}/2)-\rank R.
\]
By Table~\ref{dim table}, we know that
\[
\dim H^1(\Gamma_2,\mathbb{Z}/2)=3,\quad \dim H^1(G(32) ,\mathbb{Z}/2)=3, \quad \dim H^2(G(32),\mathbb{Z}/2)=7.
\]

Thus, it is enough to compute the rank of the matrix $R$. We use the same basis $\chi_1(v,\epsilon)=v_1$, $\chi_2(v,\epsilon)=v_2$, $\chi_3(v,\epsilon)=\epsilon$. Their values on the four  generators of the stabilizers are
\[
\begin{array}{c|cccc}
 & g_1 & g_2 & g_3 & g_4\\
\hline
\chi_1 & 0&1&1&0\\
\chi_2 & 0&1&0&1\\
\chi_3 & 1&0&0&1.
\end{array}
\]
These three vectors are linearly independent.Thus $H^1(G;\mathbb{Z}/2)=\langle\chi_1,\chi_2,\chi_3\rangle_{\mathbb{Z}/2}$. Put $H_i:=\langle g_i\rangle$, $1\leq i\leq4$. Then $H_1\cong H_2\cong H_3\cong \mathbb{Z}/2$, $H_4\cong \mathbb{Z}/4$. Define
\[
R_2:H^2(G;\mathbb{Z}/2)\longrightarrow\bigoplus_{i=1}^4H^2(H_i;\mathbb{Z}/2)\cong(\mathbb{Z}/2)^4.
\]
We first compute the restrictions of the classes $\beta(\chi_1)$, $\beta(\chi_2)$, $\beta(\chi_3)$.  By naturality, $\operatorname{res}_{H_i}\beta(\chi)=\beta(\chi|_{H_i})$.  Using Lemma~\ref{restrict} and the above table of character values, we obtain
\[
R_2(\beta(\chi_1))=(0,1,1,0),
\quad
R_2(\beta(\chi_2))=(0,1,0,0),
\quad
R_2(\beta(\chi_3))=(1,0,0,0).
\]
But we need more additional classes. We use the same central extension as in the proof for $C_1$. Recall the central extension
\[
1
\longrightarrow \langle R^4\rangle
\longrightarrow \widetilde D_4
\overset{p}{\longrightarrow} D_4
\longrightarrow 1.
\]
Let $u\in H^2(D_4;\mathbb{Z}/2)$ be the corresponding extension class. Recall also the homomorphism $q_2:G\longrightarrow D_4$ defined by
\[
\begin{array}{c|ccccc}
 & e_1&e_2&e_3&e_4&t\\
\hline
q_2&1&rs&1&r^2&s.
\end{array}
\]
Put $s_2:=q_2^*(u)\in H^2(G;\mathbb{Z}/2)$. We shall use the restriction of $u$ to the subgroup $\langle r^2\rangle\cong \mathbb{Z}/2$. Its inverse image in $\widetilde D_4$ is $p^{-1}(\langle r^2\rangle)=\langle R^2\rangle\cong \mathbb{Z}/4$. Hence the restricted extension is
\[
1
\longrightarrow \langle R^4\rangle\cong \mathbb{Z}/2
\longrightarrow \langle R^2\rangle\cong \mathbb{Z}/4
\longrightarrow \langle r^2\rangle\cong \mathbb{Z}/2
\longrightarrow 1.
\]
This extension is non-split. Indeed, the two lifts of $r^2$ are $R^2$ and $R^6$, and both have order $4$. Therefore $\operatorname{res}_{\langle r^2\rangle}(u)\neq0$. Under the identification $H^2(\mathbb{Z}/2;\mathbb{Z}/2)\cong\mathbb{Z}/2$, we write $\operatorname{res}_{\langle r^2\rangle}(u)=1$. We now compute the images of the stabilizer generators under $q_2$.

\medskip

For $g_1=e_3e_4t$,
\[
\begin{aligned}
q_2(g_1)=q_2(e_3)q_2(e_4)q_2(t)=1\cdot r^2\cdot s=r^2s.
\end{aligned}
\]
Thus,
\[
\operatorname{res}_{H_1}(s_2)=0.
\]

For $g_2=e_1e_2e_3e_4$,
\[
\begin{aligned}
q_2(g_2)=q_2(e_1)q_2(e_2)q_2(e_3)q_2(e_4)=1\cdot(rs)\cdot1\cdot r^2=(rs)r^2=r^3s.
\end{aligned}
\]
Thus,
\[
\operatorname{res}_{H_2}(s_2)=0.
\]

For $g_3=e_1e_3e_4$,
\[
\begin{aligned}
q_2(g_3)=q_2(e_1)q_2(e_3)q_2(e_4)=1\cdot1\cdot r^2=r^2.
\end{aligned}
\]
Thus,
\[
q_2(H_3)=\langle r^2\rangle\cong \mathbb{Z}/2.
\]
As shown above, the extension restricts nontrivially to
$\langle r^2\rangle$. Since $q_2|_{H_3}:H_3\longrightarrow\langle r^2\rangle$ is an isomorphism, we obtain
\[
\operatorname{res}_{H_3}(s_2)=1.
\]

Finally, for $g_4=e_2e_3e_4t$,
\[
\begin{aligned}
q_2(g_4)=q_2(e_2)q_2(e_3)q_2(e_4)q_2(t)=(rs)\cdot1\cdot r^2\cdot s=r^3.
\end{aligned}
\]
Thus $q_2(g_4)$ generates the rotation subgroup $\langle r\rangle\cong \mathbb{Z}/4$. Since the extension
\[
1
\longrightarrow \mathbb{Z}/2
\longrightarrow \mathbb{Z}/8
\longrightarrow \mathbb{Z}/4
\longrightarrow1
\]
is non-split, we obtain
\[
\operatorname{res}_{H_4}(s_2)=1.
\]
Consequently,
\[
R_2(s_2)=(0,0,1,1).
\]

We obtain
\[
\begin{array}{c|c}
\text{class} & R_2(\text{class})\\
\hline
\beta(\chi_1)&(0,1,1,0)\\
\beta(\chi_2)&(0,1,0,0)\\
\beta(\chi_3)&(1,0,0,0)\\
s_2&(0,0,1,1).
\end{array}
\]

These four classes are linearly independent in $H^2(G;\mathbb{Z}/2)$.
Indeed, suppose
\[
a_1\beta(\chi_1)+a_2\beta(\chi_2)+a_3\beta(\chi_3)+a_4s_2=0.
\]
Applying $R_2$, we obtain
\[
a_1(0,1,1,0)+a_2(0,1,0,0)+a_3(1,0,0,0)+a_4(0,0,1,1)=0.
\]
If we calculate this equation $(a_3,\;a_1+a_2,\;a_1+a_4,\;a_4)=(0,0,0,0)$, then we get the linearly independency. Since $\dim_{\mathbb{Z}/2}H^2(G;\mathbb{Z}/2)=7$, we may choose $v_5,v_6,v_7\in H^2(G;\mathbb{Z}/2)$ such that
\[
\mathcal B=\bigl(\beta(\chi_1),\beta(\chi_2),\beta(\chi_3),s_2,v_5,v_6,v_7\bigr)
\]
is a basis of $H^2(G;\mathbb{Z}/2)$. With respect to this basis and the standard basis of $(\mathbb{Z}/2)^4$, the restriction matrix has the form
\[
[R_2]_{\mathcal B}=
\begin{pmatrix}
0&0&1&0&*&*&*\\
1&1&0&0&*&*&*\\
1&0&0&1&*&*&*\\
0&0&0&1&*&*&*
\end{pmatrix}.
\]
Therefore
\[
\operatorname{rank}_{\mathbb{Z}/2}R_2=4.
\]
It follows that $\dim_{\mathbb{Z}/2}H^1(\mathbb{Z}/2;\mathbb{Z}/2)^G=3$. By \cite[Lemma 3.1]{KimKimLee17} there is a $G(32)$-invariant theta characteristic on $C_2$. Thus $\#\Theta(C_2)^G=8$.
\end{proof}

 We now determine $\operatorname{Pic}(C_i)^G$.

\medskip

In this case, by Remark~\ref{Schur} and Proposition \ref{Pic(G;C)},  we have
$$\operatorname{Hom}(G(32),\mathbb{C}^*)=(\mathbb{Z}/2)^3,\,\,\, H^{2}(G(32), \mathbb{C}^*)=(\mathbb{Z}/2)^4$$
$$\operatorname{Pic}(G;C_1)=  8\mathbb{Z} \oplus (\mathbb{Z}/2)^2 \oplus \mathbb{Z}/4,\,\,\, \operatorname{Pic}(G;C_2)= 8\mathbb{Z} \oplus (\mathbb{Z}/2)^3.$$
The exact sequence \eqref{eq:exact} yields

\begin{subequations}\label{eq:row5-sequences}
\begin{equation}\label{eq:row5-C1-concrete}
0
\longrightarrow (\mathbb{Z}/2)^3
\longrightarrow 8\mathbb{Z} \oplus (\mathbb{Z}/2)^2\oplus \mathbb{Z}/4
\longrightarrow \operatorname{Pic}(C_{1})^{G}
\longrightarrow (\mathbb{Z}/2)^4
\longrightarrow 0 .
\end{equation}
\begin{equation}\label{eq:row5-C2-concrete}
0
\longrightarrow (\mathbb{Z}/2)^3
\longrightarrow 8\mathbb{Z} \oplus (\mathbb{Z}/2)^{3}
\longrightarrow \operatorname{Pic}(C_{2})^{G}
\longrightarrow (\mathbb{Z}/2)^4
\longrightarrow 0 .
\end{equation}
\end{subequations}

 Since $\operatorname{Pic}(C_1)^G[2]= (\mathbb{Z}/2)^3$ and  $\operatorname{Pic}(C_2)^G[2]= (\mathbb{Z}/2)^3$,  we infer that
$$\operatorname{Pic}(C_1)^G= 4\mathbb{Z}\oplus (\mathbb{Z}/2)^2\oplus \mathbb{Z}/4 \quad\text{and}\quad
\operatorname{Pic}(C_2)^G= 4\mathbb{Z}\oplus (\mathbb{Z}/2)^3.$$

Next, we determine the intersection form. Let
\[
\delta_1 \colon \operatorname{Pic}(C_1)^G \longrightarrow H^2(G,\mathbb{C}^*),
\quad
\delta_2 \colon \operatorname{Pic}(C_2)^G \longrightarrow H^2(G,\mathbb{C}^*)
\]
be the obstruction maps. Both maps are surjective, since both curves are smooth projective of genus $\ge 2$.
Note that the restriction of $\delta_i$ to the torsion part of $\operatorname{Pic}(C_i)^G$ cannot be surjective.

Choose generators $L_1$, $\tau_2, \tau_3$, $\tau_4$ of $\operatorname{Pic}(C_1)^G= 4\mathbb{Z}\oplus (\mathbb{Z}/2)^2\oplus \mathbb{Z}/4$, where $\deg L_1=4$, $\tau_2$ and $\tau_3$ are 2-torsions, and $\tau_4$ a 4-torsion. The images $e_1=\delta_1(L_1)$ and $e_i=\delta_1(\tau_i)$ $(i=2,3,4)$ are generators of $ H^{2}(G(32), \mathbb{C}^*)=(\mathbb{Z}/2)^4$.

\begin{itemize}
\item Suppose that the preimage $\delta_2^{-1}(e_1)$ in $\operatorname{Pic}(C_2)^G$ consists of only torsion elements, (which implies that $\delta_2^{-1}(<e_2, e_3, e_4>)$ contains a
 degree 4 element $L_2\in \operatorname{Pic}(C_2)^G$.)
Pick such a torsion, say $M_2$, then  $L_1\boxtimes M_2$ defines a $G$-linearized line bundle on $C_1\times C_2$ of bidegree $(4,0)$.
On the other hand, since $\delta_2(L_2)\in <e_2, e_3, e_4>$, there is a torsion $M_1\in \operatorname{Pic}(C_1)^G$ such that $\delta_2(L_2)=\delta_1(M_1)$, then  $M_1\boxtimes L_2$ defines a $G$-linearized line bundle on $C_1\times C_2$ of bidegree $(0, 4)$.
Their descended line bundles have intersection number
\[
\frac{ (4, 0)\cdot(0, 4)}{|G|}=\frac{16}{|G|}=\frac{1}{2}\notin \mathbb{Z},
\]
a contradiction.
\item Suppose that the preimage $\delta_2^{-1}(e_1)$ in $\operatorname{Pic}(C_2)^G$ contains a
 degree 4 element $L_2\in \operatorname{Pic}(C_2)^G$.
 Then  $L_1\boxtimes L_2$ defines a $G$-linearized line bundle on $C_1\times C_2$ with bidegree $(4, 4)$.  By Strategy (1) its descended line bundle has self-intersection

\[
\frac{(4,4)\cdot(4,4)}{|G|}=\frac{32}{32}=1,
\]
an odd number. Thus the intersection form of $S$ is odd.
\end{itemize}

\bigskip

\section{$G=D_4\times \mathbb{Z}/2$, $g(C_1)=3$, $g(C_2)=9$, $T_1=[2,2,2,4]$, $T_2=[2,2,2,2,2,2]$; even fake quadrics}

\noindent

From \cite[Lemma 3.1, Lemma 3.2, Lemma 7.3]{LeeShabalin18}, we see that there is a $G$-linearized line bundle of bidegree $(g(C_1)-1,0)=(2,0)$ on $C_1 \times C_2$.

 On the other hand, by Proposition \ref{Pic(G;C)} the free part of $\operatorname{Pic}(G;C_{2})$ is generated by a degree $8$ line bundle.
Its pullback to $C_1 \times C_2$ is a $G$-linearized line bundle of bidegree $(0, g(C_2)-1)=(0,8)$.
Therefore, it follows from the above results  that  $S$ is an even fake quadric. \\

 Another proof has been given by a different method in \cite[Theorem 16]{CataneseFrapporti26}.\\

We include now a proof using the method of the present article.

In this case, we see that
$$\operatorname{Hom}(D_4\times \mathbb{Z}/2,\mathbb{C}^*)=(\mathbb{Z}/2)^3,\,\,\, H^{2}(D_4\times \mathbb{Z}/2, \mathbb{C}^*)=(\mathbb{Z}/2)^3$$
$$\operatorname{Pic}(G;C_1)=  4\mathbb{Z} \oplus (\mathbb{Z}/2)^3,\,\,\, \operatorname{Pic}(G;C_2)= 8\mathbb{Z} \oplus (\mathbb{Z}/2)^5.$$
The exact sequence \eqref{eq:exact} yields

\begin{subequations}\label{eq:row9-sequences}
\begin{equation}\label{eq:row9-C1-concrete}
0 \longrightarrow (\mathbb{Z}/2)^3
\longrightarrow 4\mathbb{Z} \oplus (\mathbb{Z}/2)^3
\longrightarrow \operatorname{Pic}(C_{1})^{G}
\longrightarrow (\mathbb{Z}/2)^3
\longrightarrow 0 .
\end{equation}
\begin{equation}\label{eq:row9-C2-concrete}
0 \longrightarrow (\mathbb{Z}/2)^3
\longrightarrow 8\mathbb{Z} \oplus (\mathbb{Z}/2)^{5}
\longrightarrow \operatorname{Pic}(C_{2})^{G}
\longrightarrow (\mathbb{Z}/2)^3
\longrightarrow 0 .
\end{equation}
\end{subequations}
We have also obtained that $\operatorname{Pic}(C_1)^G[2]= (\mathbb{Z}/2)^2$ and $\operatorname{Pic}(C_2)^G[2]= (\mathbb{Z}/2)^5$ and computed  $\#\Theta(C_i)^G$, but omit the proof. We infer that
$$\operatorname{Pic}(C_1)^G= 2\mathbb{Z}\oplus (\mathbb{Z}/2)^2 \quad\text{and}\quad
 \operatorname{Pic}(C_2)^G= 8 \mathbb{Z}\oplus(\mathbb{Z}/2)^5.$$
Now, from \eqref{eq:row9-C1-concrete} we can choose a $G$-invariant line bundle $L_1\in \operatorname{Pic}(C_1)^G$ of degree $2$ which is not $G$-linearized. By the exact sequence \eqref{eq:row9-C2-concrete}, there exists a torsion line bundle $\tau\in \operatorname{Pic}(C_2)^G$  whose obstruction class in the Schur multiplier coincides with that of $L_1$. Consequently, $L_1\boxtimes \tau$ defines a $G$-linearized line bundle on $C_1\times C_2$ of bidegree $(2,0)$.
On the other hand, there exists a $G$-linearized line bundle on $C_1\times C_2$ of bidegree
$(0,8)$.
By Strategy (1), the intersection form of $S$ is even.

\bigskip

\section{$G=(\mathbb{Z}/2)^3$, $g(C_1)=3$, $g(C_2)=5$, $T_1=[2,2,2,2,2]$, $T_2=[2,2,2,2,2,2]$; even fake quadrics}

Let $e_1, e_2, e_3$ be the generators of $G$.

\begin{proposition}  In the case $G=(\mathbb{Z}/2)^3$, $g(C_1)=3$, $g(C_2)=5$, $T_1=[2,2,2,2,2]$, $T_2=[2,2,2,2,2,2]$, we have
$$\operatorname{Pic}(C_1)^G[2]= (\mathbb{Z}/2)^3 \quad\text{and}\quad
\#\Theta(C_1)^G=8.$$
\end{proposition}
\begin{proof}
For $C_1$, a generating vector is
\[
(g_1,g_2,g_3,g_4,g_5)=(e_1,e_2,e_3,e_1,e_2+e_3),
\]
with orders $2,2,2,2,2$.

Let
\[
\Gamma_1=\langle x_1,x_2,x_3,x_4,x_5 \mid x_1^2=x_2^2=x_3^2=x_4^2=x_5^2=1,\; x_1x_2x_3x_4x_5=1\rangle,
\]
and let $\varphi_1:\Gamma_1\twoheadrightarrow G$, $\varphi_1(x_i)=g_i$.

 To compute $\#\Theta(C_1)^G$, by Proposition~\ref{dim} and ~\ref{kervarphi}, we compute
\[
\dim H^1(C,\mathbb{Z}/2)^G=\dim H^1(\Gamma,\mathbb{Z}/2)-\dim H^1(G,\mathbb{Z}/2)+\dim H^2(G,\mathbb{Z}/2)-\rank R.
\]
By Table~\ref{dim table}, we know that
\[
\dim H^1(\Gamma_1,\mathbb{Z}/2)=4,\quad \dim H^1((\mathbb{Z}/2)^3 ,\mathbb{Z}/2)=3, \quad \dim H^2((\mathbb{Z}/2)^3,\mathbb{Z}/2)=6.
\]

Thus, it is enough to compute the rank of the matrix $R$.
We choose the basis defined by
\[
\begin{aligned}
a(g)&=\epsilon_1,\\
b(g)&=\epsilon_2,\\
c(g)&=\epsilon_3,
\end{aligned}
\qquad
\text{for }
g=\epsilon_1e_1+\epsilon_2e_2+\epsilon_3e_3.
\]
Thus, $H^1(G;\mathbb{Z}/2)=\langle a,b,c\rangle_{\mathbb{Z}/2}$. Also, by the K\"unneth formula,
\[
H^2((\mathbb{Z}/2)^3;\mathbb{Z}/2)\cong\bigoplus_{i=1}^3 H^2(\mathbb{Z}/2;\mathbb{Z}/2)\oplus\bigoplus_{1\leq i<j\leq3}H^1(\mathbb{Z}/2;\mathbb{Z}/2)\otimes H^1(\mathbb{Z}/2;\mathbb{Z}/2).
\]
Since $\dim_{\mathbb{Z}/2}H^1(\mathbb{Z}/2;\mathbb{Z}/2)=\dim_{\mathbb{Z}/2}H^2(\mathbb{Z}/2;\mathbb{Z}/2)=1$, we have $\dim_{\mathbb{Z}/2}H^2(G;\mathbb{Z}/2)=3+3=6$. Now we compute the basis for each summands. The three $H^2(\mathbb{Z}/2;\mathbb{Z}/2)$-summands are generated by $\beta(a)$, $\beta(b)$, $\beta(c)$, where $\beta:H^1(G;\mathbb{Z}/2)\longrightarrow H^2(G;\mathbb{Z}/2)$ is the Bockstein homomorphism associated to $0\longrightarrow\mathbb{Z}/2\longrightarrow\mathbb{Z}/4\longrightarrow\mathbb{Z}/2\longrightarrow0$. Also, the three $\otimes$ terms are generated by $a\smile b$, $a\smile c$, $b\smile c$. Therefore
\[
\mathcal B=\bigl(\beta(a),\beta(b),\beta(c),a\smile b,a\smile c,b\smile c\bigr)
\]
is a basis of $H^2(G;\mathbb{Z}/2)$.

\medskip

Put $H_i:=\langle g_i\rangle$. Then $H_i\cong \mathbb{Z}/2$, $1\leq i\leq5$. Define
\[
R_1:H^2(G;\mathbb{Z}/2)\longrightarrow\bigoplus_{i=1}^5H^2(H_i;\mathbb{Z}/2)\cong(\mathbb{Z}/2)^5.
\]
We first evaluate the characters $a,b,c$ on the stabilizer generators:
\[
\begin{array}{c|ccc|c}
 & a&b&c&ord(g_i)\\
\hline
g_1=e_1       &1&0&0&2\\
g_2=e_2       &0&1&0&2\\
g_3=e_3       &0&0&1&2\\
g_4=e_1       &1&0&0&2\\
g_5=e_2+e_3   &0&1&1&2
\end{array}
\]
By naturality, $\operatorname{res}_{H_i}\beta(\chi)=\beta(\chi|_{H_i})$. Since every $H_i\cong \mathbb{Z}/2$, by Lemma~\ref{restrict},
\[
R_1(\beta(a))=(1,0,0,1,0),
\quad
R_1(\beta(b))=(0,1,0,0,1),
\quad
R_1(\beta(c))=(0,0,1,0,1).
\]
We next compute the cup-product classes. By naturality,
\[
\operatorname{res}_{H_i}(a\smile b)=\operatorname{res}_{H_i}(a)\smile\operatorname{res}_{H_i}(b).
\]
For every $i$, at least one of $a(g_i)$ and $b(g_i)$ is zero.
Therefore
\[
R_1(a\smile b)=(0,0,0,0,0).
\]
Similarly,
\[
\operatorname{res}_{H_i}(a\smile c)=\operatorname{res}_{H_i}(a)\smile\operatorname{res}_{H_i}(c).
\]
For every $i$, at least one of $a(g_i)$ and $c(g_i)$ is zero.
Therefore
\[
R_1(a\smile c)=(0,0,0,0,0).
\]
Finally,
\[
\operatorname{res}_{H_i}(b\smile c)=\operatorname{res}_{H_i}(b)\smile\operatorname{res}_{H_i}(c).
\]
Both characters $b$ and $c$ are nonzero only at $g_5=e_2+e_3$. Consequently,
\[
R_1(b\smile c)=(0,0,0,0,1).
\]
Collecting the computations, we obtain
\[
\begin{array}{c|c}
\text{class}&R_1(\text{class})\\
\hline
\beta(a)&(1,0,0,1,0)\\
\beta(b)&(0,1,0,0,1)\\
\beta(c)&(0,0,1,0,1)\\
a\smile b&(0,0,0,0,0)\\
a\smile c&(0,0,0,0,0)\\
b\smile c&(0,0,0,0,1).
\end{array}
\]
With respect to the ordered basis $\mathcal B=\bigl(\beta(a),\beta(b),\beta(c),a\smile b,a\smile c,b\smile c\bigr)$, the restriction matrix is
\[
[R_1]_{\mathcal B}=
\begin{pmatrix}
1&0&0&0&0&0\\
0&1&0&0&0&0\\
0&0&1&0&0&0\\
1&0&0&0&0&0\\
0&1&1&0&0&1
\end{pmatrix}.
\]
Therefore,
\[
\operatorname{rank} R_1=4.
\]
It follows that $\dim H^1(C_1;\mathbb{Z}/2)^G=3$. By \cite[Lemma 5.1]{Lee16} there is a $G$-invariant theta characteristic on $C_1$. Thus,
$\#\Theta(C_1)^G=8.$
\end{proof}

\begin{proposition} In the case $G=(\mathbb{Z}/2)^3$, $g(C_1)=3$, $g(C_2)=5$, $T_1=[2,2,2,2,2]$, $T_2=[2,2,2,2,2,2]$, we have
$$\operatorname{Pic}(C_2)^G[2]= (\mathbb{Z}/2)^5 \quad\text{and}\quad
\#\Theta(C_2)^G=32.$$
\end{proposition}

\begin{proof}
For $C_2$, a generating vector is given by the sextuple
\begin{align*}
g_1&=e_1+e_2,\\
g_2&=e_1+e_3,\\
g_3&=e_1+e_2+e_3,\\
g_4&=e_1+e_2,\\
g_5&=e_1+e_3,\\
g_6&=e_1+e_2+e_3.
\end{align*}
with orders $2,2,2,2,2,2$.
Let
\[
\Gamma_2=\langle y_1,y_2,y_3,y_4,y_5,y_6 \mid y_1^2=y_2^2=y_3^2=y_4^2=y_5^2=y_6^2=1,\;y_1y_2y_3y_4y_5y_6=1\rangle,
\]
and let $\varphi_2:\Gamma_2\twoheadrightarrow G$, $\varphi_2(y_i)=g_i$.
 To compute $\#\Theta(C_2)^G$, by Proposition~\ref{dim} and ~\ref{kervarphi}, we compute
\[
\dim H^1(C,\mathbb{Z}/2)^G=\dim H^1(\Gamma,\mathbb{Z}/2)-\dim H^1(G,\mathbb{Z}/2)+\dim H^2(G,\mathbb{Z}/2)-\rank R.
\]
By Table~\ref{dim table}, we know that
\[
\dim H^1(\Gamma_2,\mathbb{Z}/2)=5,\quad \dim H^1((\mathbb{Z}/2)^3 ,\mathbb{Z}/2)=3, \quad \dim H^2((\mathbb{Z}/2)^3,\mathbb{Z}/2)=6.
\]
Thus, it is enough to compute the rank of the matrix $R$.

\medskip

Recall that $
H^1(G;\mathbb{Z}/2)$ has the dual  basis $a,b,c$ with

\[
\begin{aligned}
a(g)&=\epsilon_1,\\
b(g)&=\epsilon_2,\\
c(g)&=\epsilon_3,
\end{aligned}
\qquad
\text{for }
g=\epsilon_1e_1+\epsilon_2e_2+\epsilon_3e_3.
\]
Also, from the proof for $C_1$, we have $\dim_{\mathbb{Z}/2}H^2(G;\mathbb{Z}/2)=6$ and
$$\mathcal B=\bigl(\beta(a),\beta(b),\beta(c),a\smile b,a\smile c,b\smile c\bigr)$$
 is a basis of $H^2(G;\mathbb{Z}/2).$

\medskip

Put $H_i:=\langle g_i\rangle$. Since every $g_i$ has order $2$, $H_i\cong \mathbb{Z}/2$, $1\leq i\leq6$. Define
\[
R_2: H^2(G;\mathbb{Z}/2) \longrightarrow \bigoplus_{i=1}^6H^2(H_i;\mathbb{Z}/2) \cong (\mathbb{Z}/2)^6.
\]
We first evaluate the characters $a,b,c$ on the stabilizer generators:
\[
\begin{array}{c|ccc|c}
 & a&b&c& ord(g_i)\\
\hline
g_1=e_1+e_2           &1&1&0&2\\
g_2=e_1+e_3           &1&0&1&2\\
g_3=e_1+e_2+e_3       &1&1&1&2\\
g_4=e_1+e_2           &1&1&0&2\\
g_5=e_1+e_3           &1&0&1&2\\
g_6=e_1+e_2+e_3       &1&1&1&2
\end{array}
\]
By naturality, $\operatorname{res}_{H_i}\beta(\chi)=\beta(\chi|_{H_i})$. Since every $H_i\cong \mathbb{Z}/2$, by Lemma~\ref{restrict},
\[
R_2(\beta(a))=(1,1,1,1,1,1),
\quad
R_2(\beta(b))=(1,0,1,1,0,1),
\quad
R_2(\beta(c))=(0,1,1,0,1,1).
\]
We next compute the cup-product classes. By naturality,
\[
\operatorname{res}_{H_i}(a\smile b)=\operatorname{res}_{H_i}(a)\smile\operatorname{res}_{H_i}(b).
\]
Consequently,
\[
R_2(a\smile b)=(1,0,1,1,0,1).
\]
Similarly,
\[
\operatorname{res}_{H_i}(a\smile c)=\operatorname{res}_{H_i}(a)\smile\operatorname{res}_{H_i}(c).
\]
Therefore
\[
R_2(a\smile c)=(0,1,1,0,1,1).
\]
Finally,
\[
\operatorname{res}_{H_i}(b\smile c)=\operatorname{res}_{H_i}(b)\smile\operatorname{res}_{H_i}(c).
\]
Therefore
\[
R_2(b\smile c)=(0,0,1,0,0,1).
\]
Collecting the computations, we obtain
\[
\begin{array}{c|c}
\text{class}&R_2(\text{class})\\
\hline
\beta(a)&(1,1,1,1,1,1)\\
\beta(b)&(1,0,1,1,0,1)\\
\beta(c)&(0,1,1,0,1,1)\\
a\smile b&(1,0,1,1,0,1)\\
a\smile c&(0,1,1,0,1,1)\\
b\smile c&(0,0,1,0,0,1).
\end{array}
\]
With respect to the ordered basis $\mathcal B=\bigl(\beta(a),\beta(b),\beta(c),a\smile b,a\smile c,b\smile c\bigr)$, the restriction matrix is
\[
[R_2]_{\mathcal B}=
\begin{pmatrix}
1&1&0&1&0&0\\
1&0&1&0&1&0\\
1&1&1&1&1&1\\
1&1&0&1&0&0\\
1&0&1&0&1&0\\
1&1&1&1&1&1
\end{pmatrix}.
\]
Therefore,
\[
\operatorname{rank} R_2=3.
\]
It follows that $\dim H^1(\mathbb{Z}/2;\mathbb{Z}/2)^G=5$. Example below gives an invariant theta characteristic on $C_2$, hence $\#\Theta(C_2)^G=32.$
\end{proof}

\begin{example}(A $G$-invariant theta characteristic for $\pi_2:C_2\longrightarrow C_2/G\cong\mathbb{P}^1$)\\
Let $F_1,\ldots,F_6$ denote the reduced fibers over the six branch points. Let $F=\pi_2^*(p)$ be a general fiber. Since every ramification index is $2$, $2F_i\equiv F$ for every $i$. By the Hurwitz formula,
\[
K_{C_2} \equiv \pi_2^*K_{\mathbb{P}^1} + F_1+\cdots+F_6 \equiv -2F+F_1+\cdots+F_6.
\]
By \cite[Lemma 5.5]{Lee16},
\[
F_1\equiv F_4, \quad F_2\equiv F_5, \quad F_3\equiv F_6.
\]
Thus, we obtain
\[
\begin{aligned}
K_{C_2}\equiv -2F+2F_1+2F_2+2F_3=2(-F+F_1+F_2+F_3).
\end{aligned}
\]
Therefore
\[
\kappa_{C_2}:=\mathcal{O}_{C_2}(-F+F_1+F_2+F_3)
\]
is a theta characteristic. Since each $F_i$ is $G$-invariant,
so is $\kappa_{C_2}$.
\end{example}

 We now determine $\operatorname{Pic}(C_i)^G$ and the parity of the intersection form of $S$.

\medskip

In this case, by Remark~\ref{Schur} and Proposition~\ref{Pic(G;C)}  we have
$$\operatorname{Hom}((\mathbb{Z}/2)^3,\mathbb{C}^*)=(\mathbb{Z}/2)^3,\,\,\, H^{2}((\mathbb{Z}/2)^3, \mathbb{C}^*)=(\mathbb{Z}/2)^3$$
$$\operatorname{Pic}(G;C_1)=  4\mathbb{Z} \oplus (\mathbb{Z}/2)^4,\,\,\, \operatorname{Pic}(G;C_2)= 4\mathbb{Z} \oplus (\mathbb{Z}/2)^5.$$
The exact sequence \eqref{eq:exact} yields
\begin{subequations}\label{eq:row12-sequences}
\begin{equation}\label{eq:row12-C1-concrete}
0
\longrightarrow (\mathbb{Z}/2)^3
\longrightarrow 4\mathbb{Z} \oplus (\mathbb{Z}/2)^4
\longrightarrow \operatorname{Pic}(C_{1})^{G}
\longrightarrow (\mathbb{Z}/2)^3
\longrightarrow 0 .
\end{equation}
\begin{equation}\label{eq:row12-C2-concrete}
0
\longrightarrow (\mathbb{Z}/2)^3
\longrightarrow 4\mathbb{Z} \oplus (\mathbb{Z}/2)^{5}
\longrightarrow \operatorname{Pic}(C_{2})^{G}
\longrightarrow (\mathbb{Z}/2)^3
\longrightarrow 0 .
\end{equation}
\end{subequations}
 Since $\operatorname{Pic}(C_1)^G[2]= (\mathbb{Z}/2)^3$ and  $\operatorname{Pic}(C_2)^G[2]= (\mathbb{Z}/2)^5$,  from the above exact sequences  we infer that
$$\operatorname{Pic}(C_1)^G= 2\mathbb{Z}\oplus (\mathbb{Z}/2)^3 \quad\text{and}\quad \operatorname{Pic}(C_2)^G= 4\mathbb{Z}\oplus(\mathbb{Z}/2)^5.$$
By  \eqref{eq:row12-C1-concrete}, we can choose a $G$-invariant line bundle $L_1\in \operatorname{Pic}(C_1)^G$ of degree $2$, and  by the exact sequence \eqref{eq:row12-C2-concrete}, there exists a torsion line bundle $\tau\in \operatorname{Pic}(C_2)^G$  whose obstruction class in the Schur multiplier coincides with that of $L_1$. Consequently, $L_1\boxtimes \tau$ defines a $G$-linearized line bundle on $C_1\times C_2$ of bidegree $(2,0)$.
On the other hand, by the exact sequence \eqref{eq:row12-C2-concrete} there exists a $G$-linearized line bundle of degree 4 on $C_2$, whose pullback  on $C_1\times C_2$ gives a $G$-linearized line bundle of bidegree
$(0,4)$.

By Strategy (1), the intersection form of $S$ is even.

\medskip

We have completed the proof of Theorem \ref{main}.

\bigskip

\noindent{\bf Acknowledgments} Kim and Lee were partially supported by the National Research Foundation (NRF) of Korea grant funded by the Korea government (MSIT)(RS-2024-00394189) and the POSCO Science Fellowship of POSCO TJ Park Foundation. Kim was also partially supported by the Basic Research Laboratory Program of NRF (RS-2026-25508638). Lee was also partially supported by Basic Science Research Institute Fund, whose NRF grant number is RS-2021-NR060139.

\bigskip

\noindent{\bf Tool Use} The authors conceived and developed all the mathematical arguments and wrote the proofs and the manuscript. They are solely responsible for the final text. Generative AI tools (ChatGPT and Claude) were used only for minor English editing and stylistic improvements.

\end{document}